\newif\ifarxiv\arxivfalse
\arxivtrue
\ifarxiv
	\documentclass[10pt, twoside, reqno]{amsart}
\else
	\RequirePackage{amsmath}
	\documentclass[smallextended]{svjour3}
\fi

\PassOptionsToPackage{noend}{algpseudocode}
\PassOptionsToPackage{dvipsnames}{xcolor}
\PassOptionsToPackage{capitalize}{cleveref}
\PassOptionsToPackage{bookmarksdepth=4}{hyperref}

\ifarxiv
	\usepackage[foot]{amsaddr}
	\newcommand{\review}[1]{#1}
\else
	\newcommand{\review}[1]{{\color{orange}#1}}
	\makeatletter
		\let\cl@chapter\undefined
	\makeatother

	\smartqed
\fi

\usepackage{longtable}
\usepackage[%
	noload={subcaption,todonotes},
	eqreset=section,
	thmreset=section,
]{myPreamble}

\allowdisplaybreaks
\ifarxiv\else
	
	\makeatletter
		\def\operator@font{\rm}
	\makeatother
\fi

\ifarxiv
	\BeforeBeginEnvironment{table}{\bgroup\footnotesize}%
	\AfterEndEnvironment{table}{\egroup}%
	\BeforeBeginEnvironment{algorithmic}{\bgroup\footnotesize}%
	\AfterEndEnvironment{algorithmic}{\egroup}%
\fi

	\begingroup\lccode`~=`&
	\lowercase{\endgroup\let~}&
	\catcode`&=\active
	\AfterEndPreamble{%
		\def\noampersand{\def&{}}
		\ifarxiv
			\newenvironment{Align}{\begin{equation}\noampersand}{\end{equation}}
			\newenvironment{Align*}{\csname equation*\endcsname\noampersand}{\csname endequation*\endcsname}
		\else
			
			\newenvironment{Align*}{\csname align*\endcsname}{\csname endalign*\endcsname}
		\fi
	}

\DeclarePairedDelimiter\ceil{\lceil}{\rceil}

\DeclareMathOperator{\symm}{Sym}

\newcommand{\smallinv}{{\scriptscriptstyle -1}}
\newcommand{\smallinf}{{\scriptscriptstyle\infty}}
\newcommand{\submatrix}[1]{_{\text{\tiny{\(#1\)}}}}
\newcommand{\supmatrix}[1]{^{\text{\tiny{\(#1\)}}}}

\renewcommand{\I}{\mathcal I}
\newcommand{\Id}{I} 
\newcommand{\J}{\mathcal J}
\newcommand{\K}{\mathcal K}
\newcommand{\M}{\mathcal M}
\newcommand{\LL}{\mathcal L}
\newcommand{\LDL}{LD\trans L}
\newcommand{\U}{\mathcal U}

\newcommand{\x}{\text{\sf\textsc x}}
\newcommand{\y}{\text{\sf\textsc y}}
\newcommand{\X}{\mathcal X}

\makeatletter
	\def\Sx{%
		\Sigma%
		\@ifnextchar_{\Sx@}{_\x}%
	}
	\def\Sx@_#1{%
		_{\x,#1}%
	}
	\def\Sy{%
		\Sigma%
		\@ifnextchar_{\Sy@}{_\y}%
	}
	\def\Sy@_#1{%
		_{\y,#1}%
	}
\makeatother

\crefname{figure}{Figure}{Figures}

	\makeatletter
		\newcommand\OrigAlglinenumber[1]{{
			\color{black}\footnotesize
			\textbf{\thealgorithm}.%
			\fillwidthof[l]{\oldstylenums{\oldstylenums{88}:}}{{\arabic{ALG@line}}:}%
		}}
		\newcommand\NewAlglinenumber[1]{{
			\color{black}\footnotesize
			\textbf{\thealgorithm}.%
			\fillwidthof[l]{\oldstylenums{\oldstylenums{88}:}}{{\arabic{ALG@line}}:\smash{\rlap{\(^\star\)}}}%
		}}
		\let\alglinenumber\OrigAlglinenumber
		\renewcommand\theALG@line{\thealgorithm.\oldstylenums{\arabic{ALG@line}}}
		\providecommand\theHALG@line{\thealgorithm.\arabic{ALG@line}}
	\makeatother
	
	\AtBeginEnvironment{algorithmic}{\linespread{1.75}\selectfont}

	\makeatletter
		\let\OLDState\State
		\let\OLDIf\If
		\let\OLDElse\Else
		\newcommand{\@State}{\OLDState\color{black}}
		\newcommand{\@@State}{%
			\renewcommand\alglinenumber[1]{\NewAlglinenumber{##1}}%
			\OLDState
			\renewcommand\alglinenumber[1]{\OrigAlglinenumber{##1}}%
		}
		\renewcommand{\State}{\@ifstar\@@State\@State}
		\newcommand{\@If}[1]{\OLDIf{#1}\color{black}}
		\newcommand{\@@If}[1]{%
			\renewcommand\alglinenumber[1]{\NewAlglinenumber{##1}}%
			\OLDIf{#1}%
			\renewcommand\alglinenumber[1]{\OrigAlglinenumber{##1}}%
		}
		\renewcommand{\If}{\@ifstar\@@If\@If}
		\newcommand{\@Else}{\OLDElse\color{black}}
		\newcommand{\@@Else}{%
			\renewcommand\alglinenumber[1]{\NewAlglinenumber{##1}}%
			\OLDElse
			\renewcommand\alglinenumber[1]{\OrigAlglinenumber{##1}}%
		}
		\renewcommand{\Else}{\@ifstar\@@Else\@Else}
	\makeatother

	\definecolor{amber}{rgb}{1.0, 0.49, 0.0}%
	\pgfplotsset{%
		myline/.style  = {line width = 1pt},
		QPALM/.style   = {myline, color = blue},
		qpOASES/.style = {myline, color = green!75!black},
		OSQP/.style    = {myline, color = red},
		IPOPT/.style   = {myline, color = cyan},
		Gurobi/.style  = {myline, color = black},
		HPIPM/.style   = {myline, color = amber},
		myaxis/.style = {%
			title style  = {font = \bfseries},
			xlabel style = {font = \color{white!15!black}},
			ylabel style = {font = \color{white!15!black}},
			axis background/.style = {fill = white},
			minor grid style = {dotted, black!20},
			major grid style = {dotted, black!50},
		},
		PPaxis/.style = {%
			myaxis,
			xmin = 0,
			ymin = 0,
			ymax = 1,
			yminorticks = true,
			xlabel = {\(\log_{10}(f)\)},
			ylabel = {fraction of solver within \(f\) of best},
			legend style = {mylegend, bottomright},
		},
		showgrids/.style = {%
			xmajorgrids,
			ymajorgrids,
			yminorgrids,
			minor grid style={dotted,black!20},
			major grid style={dotted,black!50},
		},
		mylegend/.style = {%
			rounded corners,
			align = left,
			draw = black,
			text opacity = 1,
			fill = white,
			fill opacity = 0.6,
			legend cell align = left,
		},
		topright/.style    = {at = {(0.975,0.975)}, anchor = north east},
		topleft/.style     = {at = {(0.025,0.975)}, anchor = north west},
		bottomright/.style = {at = {(0.975,0.025)}, anchor = south east},
		bottomleft/.style  = {at = {(0.025,0.025)}, anchor = south west},
	}
\showgrayfalse
\ifarxiv
	\renewcommand{\includetikz}[2][]{\includegraphics[#1]{Pics/Tikz/#2.pdf}}
\fi

\newcommand{\TheTitle}{QPALM: A Proximal Augmented Lagrangian Method for Nonconvex Quadratic Programs}
\newcommand{\TheShortTitle}{\TheTitle}
\newcommand{\TheShortAuthor}{B. Hermans, A. Themelis and P. Patrinos}
\newcommand{\TheFunding}{%
	Ben Hermans is with the MECO Research Team, Department of Mechanical Engineering, KU Leuven, and Flanders Make - DMMS\textunderscore M, Leuven, Belgium.
	His research benefits from KU Leuven-BOF PFV/10/002 Centre of Excellence: Optimization in Engineering (OPTEC), from project G0C4515N of the Research Foundation - Flanders (FWO - Flanders), from Flanders Make ICON: Avoidance of collisions and obstacles in narrow lanes, and from the KU Leuven Research project C14/15/067: B-spline based certificates of positivity with applications in engineering.%
	\newline
	Andreas Themelis is with the Faculty of Information Science and Electrical Engineering, Kyushu University, 744 Motooka, Nishi-ku, Fukuoka  819-0395 Japan.
	This work was supported by the JSPS KAKENHI grant number JP21K17710.
	\newline
	Panagiotis Patrinos is with the \TheAddressKU.
	This work was supported by the \emph{Research Foundation Flanders (FWO)} research projects G086518N, G086318N, and G0A0920N;
	\emph{Research Council KU Leuven} C1 project No. C14/18/068;
	\emph{Fonds de la Recherche Scientifique --- FNRS and the Fonds Wetenschappelijk Onderzoek --- Vlaanderen} under EOS project no 30468160 (SeLMA)%
}
\newcommand{\TheKeywords}{%
	Nonconvex QPs,
	proximal augmented Lagrangian,
	semismooth Newton method,
	exact linesearch,
	factorization updates%
}
\newcommand{\TheSubjclass}{%
	90C05, 
	90C20, 
	90C26, 
	49J53, 
	49M15.
}
\newcommand{\TheAbstract}{%
	We propose QPALM, a nonconvex quadratic programming (QP) solver based on the proximal augmented Lagrangian method.
	This method solves a sequence of inner subproblems which can be enforced to be strongly convex and which therefore \review{admit a} unique solution.
	The resulting steps are shown to be equivalent to inexact proximal point iterations on the extended-real-valued cost function\review{, which allows for a fairly simple analysis where convergence to a stationary point at an \(R\)-linear rate is shown.}
	The QPALM algorithm solves the subproblems iteratively using semismooth Newton directions and an exact linesearch.
	\review{The former can be computed efficiently in most iterations by making use of suitable factorization update routines, while the latter requires the zero of a monotone, one-dimensional, piecewise affine function.}
	QPALM is implemented in open-source C code, with tailored linear algebra routines for the factorization in a self-written package LADEL.
	The resulting implementation is shown to be extremely robust in numerical simulations, solving all of the Maros-Meszaros problems and finding a stationary point for most of the nonconvex QPs in the Cutest test set.
	Furthermore, it is shown to be competitive against state-of-the-art convex QP solvers in typical QPs arising from application domains such as portfolio optimization and model predictive control.
	As such, QPALM strikes a unique balance between solving both easy and hard problems efficiently.%
}

\ifarxiv
	\title[\TheShortTitle]{\LARGE\sc\TheTitle}

	\author[B. Hermans]{Ben Hermans}
	\author[A. Themelis]{Andreas Themelis}
	\author[P. Patrinos]{Panagiotis Patrinos}

	\thanks{\TheFunding}

	\email[B. Hermans]{ben.hermans2@kuleuven.be}
	\email[A. Themelis]{andreas.themelis@ees.kyushu-u.ac.jp}
	\email[P. Patrinos]{panos.patrinos@esat.kuleuven.be}

	\subjclass{\TheSubjclass}

\begin{document}
		\vspace*{-2cm}%
		\begin{abstract}
			\TheAbstract
		\end{abstract}

		\maketitle

\else
	\begin{document}
		\journalname{Mathematical Programming Computation}

		\title{\TheTitle\thanks{\TheFunding.}}
		\titlerunning{\TheShortTitle}

		\author{%
			Ben Hermans\and
			Andreas Themelis\and
			Panagiotis Patrinos%
		}
		\authorrunning{\TheShortAuthor}

		\institute{%
			B. Hermans
			\at
			Tel.: +32 (0)16 373659\\
			\email{ben.hermans2@kuleuven.be}\\
			Department of Mechanical Engineering --
			KU Leuven, Celestijnenlaan 300, 3001 Leuven, Belgium
		\and
			A. Themelis
			\at
			\email{andreas.themelis@ees.kyushu-u.ac.jp}\\
			Faculty of Information Science and Electrical Engineering, Kyushu University, 744 Motooka, Nishi-ku, Fukuoka 819-0395, Japan
		\and
			P. Patrinos
		\at
			Tel.: +32 (0)16 374445\\
			\email{panos.patrinos@esat.kuleuven.be}\\
			\TheAddressKU
		}

		\date{Received: date / Accepted: date}

		\maketitle

		\begin{abstract}
			\TheAbstract
			\subclass{\TheSubjclass}%
		\end{abstract}
\fi


	\section{Introduction}

		This paper considers QPs, namely
		\[\tag{QP}\label{eq:P}
			\minimize_{x\in\R^n}\tfrac12\trans xQx+\trans qx
		\quad\stt{}
			Ax\in C,
		\]
		where \review{\(Q\in\R^{n\times n}\) is symmetric,} \(A\in\R^{m\times n}\)\review{,} and \(C=\set{z\in\R^m}[\ell\leq z\leq u]\) for some vectors \(\ell,u\in\R^m\) is a box.
		\ifarxiv\else
			\par
		\fi
		Convex QPs, in which \review{\(Q\) is positive semidefinite}, are ubiquitous in numerical optimization, as they arise in many applications domains such as portfolio optimization, support vector machines, sparse regressor selection, linear \review{model predictive control} (MPC), etc.
		The solution of a QP is also required in the general nonlinear optimization technique known as sequential quadratic programming (SQP).
		Therefore, substantial research has been performed to develop robust and efficient QP solvers.
		State-of-the-art algorithms to solve convex QPs typically fall into one of three categories: active-set methods, interior-point methods or first-order methods.
		
		Active-set methods iteratively determine a working set of active constraints and require the solution of a linear system every time this set changes.
		Because the change is small however, typically restricted to \review{one or two} constraints, the linear system changes only slightly and low-rank factorization updates can be used to make this method efficient.
		An advantage is that active-set methods can easily make use of an initial guess, also known as a warm-start, which is very useful when solving a series of related QPs, such as in SQP or in MPC.
		The biggest drawback of active-set methods, however, is that a large number of iterations can be required to converge to the right active set, as the number of possible sets grows exponentially with the number of constraints.
		Popular active-set\review{-}based QP solvers include the open-source solver qpOASES \cite{ferreau2014qpoases} and the QPA module in the open-source software GALAHAD \cite{gould2003galahad}.
		
		Interior point methods typically require fewer but more expensive iterations than active-set methods.
		Their iterations involve the solution of a new linear system at every iteration.
		Interior-point methods are generally efficient, but suffer from not having warm-starting capabilities.
		Examples of state-of-the-art QP solvers using an interior-point method are the commercial solvers Gurobi \cite{gurobi2018gurobi} and MOSEK \cite{mosek2019mosek}, the closed-source solver BPMPD \cite{meszaros1999bpmpd} and the open-source solver OOQP \cite{gertz2003object}.
		
		First-order methods rely only on first-order information of the problem.
		Particularly popular among first-order methods are the proximal algorithms, also known as operator splitting methods.
		Such methods can typically be described in terms of simple operations, and their iterations are relatively inexpensive.
		They may, however, exhibit slow asymptotic convergence for ill-conditioned problems.
		The recently proposed OSQP solver \cite{stellato2020osqp}, based on the alternating direction method of multipliers (ADMM), addresses this crucial issue somewhat by means of preconditioning.
		
		It is generally difficult to extend the aforementioned methods to be able to find stationary points of nonconvex QPs, without additional assumptions.
		Augmented-Lagrangian-based algorithms such as the ADMM, for example, would require surjectivity of the constraint matrix \(A\) \cite{li2015global,bolte2018nonconvex,bot2020proximal,themelis2020douglas}.
		Some proposals have been made for interior-\review{point} methods to solve nonconvex QPs \cite{ye1992affine,absil2007newton}, but these methods were found to often exhibit numerical issues in our benchmarks.
		\review{The active-set solvers SQIC \cite{gill2015methods} and qpOASES \cite{ferreau2014qpoases} are also able to find critical points of nonconvex QPs and include checks for second-order sufficient conditions to identify local minima.
		However, the former is not publicly available and the latter is tailored to small-to-medium-scale problems.}
		Finally, global optimization of nonconvex QPs has been the topic of a large amount of research, see for example \cite{sherali1995reformulation,burer2008finite,chen2012globally}, but this is a separate issue and will not be discussed further here, as we are only interested in finding a stationary point, \review{characterized by the first-order necessary conditions for optimality.
		
		As mentioned in \cite[Result 2.1]{gill2015methods} for instance, the second-order necessary and sufficient condition for optimality requires the positive-definiteness of \(Q\) over all feasible directions orthogonal to the local gradient.
		However, verifying this condition requires finding the global minimizer of an indefinite quadratic form over a cone, which is an NP-hard problem \cite{cottle1970classes}.
		The authors of \cite{gill2015methods} propose additionally a necessary (but not sufficient) second-order criterion verifiable in polynomial time, which consists of verifying the positive-definiteness of \(Q\) on the nullspace of \(A_\I\), where \(\I\) is the index set of active constraints.
		This paper restricts itself to considering only first-order conditions, with the exception of using the criterion above a posteriori in the simulations on nonconvex QPs in \cref{sec:Simulations}.}
		
		In this paper we show that the proximal augmented Lagrangian method (P-ALM), up to a simple modification, still enjoys convergence guarantees without convexity or surjectivity assumptions.
		In particular, this allows us to extend the recently proposed convex QP solver QPALM \cite{hermans2019qpalm} to nonconvex QPs.


			P-ALM when applied to convex problems has been shown to be equivalent to resolvent iterations on the monotone operator encoding the KKT optimality conditions \cite{rockafellar1976augmented}.
			While this interpretation is still valid for \eqref{eq:P} under our working assumptions, the resulting KKT system lacks the monotonicity requirement that is needed for guaranteeing convergence of the iterates.
			In fact, while \(\M\) is hypo-monotone, in the sense that it can be made monotone by adding a suitably large multiple of the identity mapping, the same cannot be said about its inverse \(\M^{-1}\) whence recent advancements in the nonconvex literature would apply, see \cite{iusem2003inexact,combettes2004proximal}.
			For this reason, we here propose a different interpretation of \emph{a P-ALM step as an inexact proximal point iteration on the extended-real-valued cost}
			\[
				\varphi(x)
			{}\coloneqq{}
				\tfrac12\trans xQx+\trans qx+\indicator_C(Ax),
			\]
			where \(\indicator_C\) is the indicator function of set \(C\), namely \(\indicator_C(x)=0\) if \(x\in C\) and \(\infty\) otherwise.
			\review{As will be better detailed in \cref{sec:PALM}, the proximal point (PP) subproblems are} addressed by means of an ALM method where the hard constraint \(Ax\in C\) is replaced by a quadratic penalty\review{, overall resulting in the proposed \emph{proximal-ALM method for quadratic programs QPALM}}.
			
			Some recent papers \cite{kong2019complexity,lin2019inexact} developed and analyzed the iteration complexity of a closely related three-layer algorithm, not specifically for QPs, which involves solving a series of quadratically penalized subproblems using inexact proximal point iterations\review{.}
			\review{Although dealing} with nonconvexity in the objective through the proximal penalty \review{as well, the approach} is quite different from ours in that it uses a pure quadratic penalty instead of ALM\review{, and thus} requires the penalty parameters to go to infinity and is prone to exhibit \review{slower convergence rates} \cite[\S2.2.5]{bertsekas1982constrained}.
			Furthermore, the inner subproblems are solved using \review{an accelerated (first-order)} composite gradient method, whereas QPALM uses a semismooth Newton method.
			Finally, \review{to the best of our knowledge}, no code for their algorithm was provided.

		\subsection{Contributions}

			Our contributions can be summarized as follows.
			\begin{enumerate}
			\item
				We show the equivalence between P-ALM and \review{PP} iterations, and more specifically the relation between the inexactness in both methods.
				As such, we can make use of the convergence results of \cite[\S4.1]{sun2017convergence}.
				In addition, we show that inexact PP on possibly nonconvex QPs is globally convergent and, in fact, with \(R\)-linear rates, thus complementing \cite[Prop. 3.3]{luo1993error} that covers the exact case.
			\item
				We modify the QPALM algorithm introduced in a previous paper for convex QPs \cite{hermans2019qpalm}, such that it can now also deal with nonconvex QPs.
				We highlight the (minimal) changes required, and add a self-written C version of the LO(B)PCG algorithm, used to find the minimum eigenvalue, to the QPALM source code.
			\item
				We outline the necessary linear algebra routines and present a standalone C package LADEL that implements them.
				Therefore, differently from our previous version which relied on CHOLMOD \cite{chen2008algorithm}, \review{QPALM is now a standalone software package.}
				Furthermore, all the details of the parameter selection and initialization routines are outlined here.
			\item
				We provide extensive benchmarking results, not only for nonconvex QPs which we obtain from the Cutest test set \cite{gould2015cutest}, but also for convex QPs.
				Here, we vastly extend the limited results presented in \cite{hermans2019qpalm} by solving the Maros-Meszaros problems, many of which are \review{large scale and ill conditioned}, alongside some QPs arising from specific application domains, such as portfolio optimization and model predictive control.
			\end{enumerate}

		\subsection{Notation}\label{subsec:Notation}

			The following notation is used throughout the paper.
			We denote the extended real line by \(\Rinf\coloneqq\R\cup\set{\infty}\).
			The scalar product on \(\R^n\) is denoted by \(\innprod{{}\cdot{}}{{}\cdot{}}\).
			With $[x]_+\coloneqq\max\set{x,0}$ we indicate the positive part of vector \(x\in\R^n\), meant in a componentwise sense.
			A sequence of vectors \(\seq{x^k}\) is said to be summable if \(\sum_{k\in\N}\|x^k\|<\infty\).
			
			With \(\symm(\R^n)\) we indicate the set of symmetric \(\R^{n\times n}\) matrices, while \(\symm_+(\R^n)\) and \(\symm_{++}(\R^n)\) denote the subsets of those which are positive semidefinite and positive definite, respectively.
			\review{For two matrices \(A,B\in\symm(\R^n)\) we write \(A\succeq B\) (resp. \(A\succ B\)) to indicate that \(A-B\) is positive semidefinite (resp. positive definite), and given} \(\Sigma\in\symm_{++}(\R^n)\) we indicate with
			\(
				\|{}\cdot{}\|_\Sigma
			\)
			the norm on \(\R^n\) induced by \(\Sigma\), namely
			\(
				\|x\|_\Sigma
			{}\coloneqq{}
				\sqrt{\innprod{x}{\Sigma x}}
			\).
			\review{With \(\Id_n\) we denote the \(n\times n\)-identity matrix, and we simply write \(\Id\) when \(n\) is clear from context.}
			
			Given a nonempty closed convex set \(C\subseteq\R^n\), with \(\proj_C(x)\) we indicate the projection of a point \(x\in\R^n\) onto \(C\), namely \(\proj_C(x)=\argmin_{y\in C}\|y-x\|\) or, equivalently, the unique point \(z\in C\) satisfying the inclusion
			\begin{equation}\label{eq:proj}
				x-z\in\ncone_C(z),
			\end{equation}
			where
			\(
				\ncone_C(z)
			{}\coloneqq{}
				\set{v\in\R^n}[\innprod{v}{z-z'}\leq 0~\forall z'\in C]
			\)
			is the normal cone of the set \(C\) at \(z\).
			\(\dist(x,C)\) and \(\dist_\Sigma(x,C)\) denote the distance from \(x\) to set \(C\) in the Euclidean norm and in that induced by \(\Sigma\), respectively, while \(\indicator_C\) is the indicator function of set \(C\), namely
			\(
				\indicator_C(x)=0
			\)
			if \(x\in C\) and \(\infty\) otherwise.
			For a set of natural numbers \(\I\in\N\) we let \(|\I|\) denote its cardinality, whereas for a (sparse) matrix \(A\in\R^{m\times n}\) we let \(|A|\) denote the number of nonzero elements in \(A\).
			The element of \(A\) in the \(i\)-th row and \(j\)-th column is denoted as \(A_{ij} \in \R\).
			For an index \(i \in [1,m]\), let \(A_{i\cdot}\) denote the \(i-\)th row of \(A\).
			Similarly, for a set of indices \(\mathcal I \subseteq [1,m]\), let \(A\submatrix{\I\cdot} \in \R^{|\I|\times n}\) be the submatrix comprised of all the rows \(i \in \mathcal I\) of \(A\).
			Analogously, for \(j \in [1,n]\), and \(\mathcal J \subseteq [1,n]\), let \(A_{\cdot j}\) denote the \(j\)-th column of \(A\), and \(A\submatrix{\cdot\J}\) \review{the submatrix comprised of all the columns \(j \in \mathcal J\) of \(A\)}.
			Combined, let \(A\submatrix{\I\J} \in \R^{|\I|\times |\J|}\) denote the submatrix comprised of all the rows \(i \in \mathcal I\) and all the columns \(j \in \mathcal J\) of \(A\).
			Finally, let us denote the matrix \(A\supmatrix{\I\cdot} \in \R^{m\times n}\) as the matrix with the corresponding rows from \(A\) and \(\mathbf{0}\) elsewhere, i.e.
			\[ \label{eq:supmatrix}
				A_{i\cdot} \supmatrix{\I} = 
				\begin{ifcases}
					A_{i\cdot} &  i\in \I \\
					\mathbf{0} \otherwise,	
				\end{ifcases}
			\]
			and similarly \(A\supmatrix{\I\J} \in \R^{m\times n}\) the matrix with elements
			\[ 
				A_{ij} \supmatrix{\I\J} = 
				\begin{ifcases}
					A_{ij} &  i\in \I \textrm{ and } j\in\J \\
					0 \otherwise.
				\end{ifcases}
			\]
			
			\review{A mapping \(\func{F}{\R^n}{\R^m}\) is Lipschitz continuous on \(\Omega\subseteq\R^n\) if there exists \(L\geq0\) such that \(\|F(x)-F(y)\|\leq L\|x-y\|\) holds for all \(x,y\in\Omega\).
			The smallest such constant \(L\) is the Lipschitz modulus of \(F\) on \(\Omega\), denoted as \(\lip_\Omega F\) or simply \(\lip F\) when \(\Omega=\R^n\).
			We say that a real-valued function \(\func h{\R^n}{\R}\) is Lipschitz differentiable if \(h\) is continuosly differentiable and its gradient \(\nabla h\) is Lipschitz continuous on \(\R^n\).
			We may also say that \(h\) is \(L_h\)-smooth as a shorthand notation to indicate that \(h\) is Lipschitz differentiable with \(\lip\nabla h=L_h\).
			
			For an extended-real-valued function \(\func{f}{\R^n}{\Rinf}\) and \(\alpha\in\R\) we indicate with
			\(
				\lev_{\leq\alpha}f
			{}\coloneqq{}
				\set{x\in\R^n}[f(x)\leq\alpha]
			\)
			the \(\alpha\)-sublevel set of \(f\), and we say that \(f\) is level bounded if \(\lev_{\leq\alpha}f\) is a bounded set for any \(\alpha\in\R\), this condition being equivalent to \(\lim_{\|x\|\to\infty}f(x)=\infty\).}

		\subsection{Nonconvex subdifferential and proximal mapping}

			The \DEF{regular subdifferential} of \(\func{f}{\R^n}{\Rinf}\) at \(x\in\R^n\) is the set \(\hat\partial f(x)\), where
			\[
				v\in\hat\partial f(x)
			\quad\text{iff}\quad
				\liminf_{x'\to x}\frac{f(x')-f(x)-\innprod{v}{x'-x}}{\|x'-x\|}
			{}\geq{}
				0,
			\]
			whereas the \DEF{(limiting) subdifferential} of \(f\) at \(x\) is \(\review{\partial f(x)}=\emptyset\) if \(x\notin\dom f\), and
			\[
				\partial f(x)
			{}\coloneqq{}
				\set{v\in\R^n}[
					\exists (x^k,v^k)\to (x,v)
				~\text{such that}~
					f(x^k)\to f(x)
				~\text{and}~
					v^k\in\hat\partial f(x^k)
				~
					\forall k
				]
			\]
			otherwise.
			Notice that \(\hat\partial f(x)\subseteq\partial f(x)\) for any \(x\in\R^n\), and that the inclusion \(0\in\hat\partial f(x)\) is a necessary condition for local minimality of \(x\) for \(f\) \cite[Thm.s 8.6 and 10.1]{rockafellar2011variational}.
			A point \(x\) satisfying this inclusion is said to be \DEF{stationary} (for \(f\)).
			If \(f\) is proper lower semicontinuous (lsc) and convex, then
			\[
				\hat\partial f(x)
			{}={}
				\partial f(x)
			{}={}
				\set{v\in\R^n}[
					f(x')\geq f(x)+\innprod{v}{x'-x}
					~
					\forall x'\in\R^n
				],
			\]
			and stationarity of \(x\) for \(f\) is a necessary and sufficient condition for global minimality \cite[Prop. 8.12 and Thm. 10.1]{rockafellar2011variational}.
			For a proper lsc function \(f\) and \(\Sx\in\symm_{++}(\R^n)\), the \DEF{proximal mapping} of \(f\) with (matrix) stepsize \(\Sx\) is the set-valued mapping \(\ffunc{\prox_f^{\Sx}}{\R^n}{\R^n}\) given by
			\begin{align*}
				\prox_f^{\Sx}(x)
			{}\coloneqq{} &
				\argmin_{w\in\R^n}\set{
					f(w)
					{}+{}
					\tfrac12\|w-x\|_{\Sx^{-1}}^2
				},
			\shortintertext{%
				and the corresponding \DEF{Moreau envelope} is \(\func{f^{\Sx}}{\R^n}{\R}\) defined as
			}
				f^{\Sx}(x)
			{}\coloneqq{} &
				\min_{w\in\R^n}\set{
					f(w)
					{}+{}
					\tfrac12\|w-x\|_{\Sx^{-1}}^2
				}.
			\end{align*}
			It follows from the definition that \(\bar x\in\prox_f^{\Sx}(x)\) iff
			\begin{equation}\label{eq:proxIneq}
				f^{\Sx}(x)
			{}={}
				f(\bar x)
				{}+{}
				\tfrac12\|x-\bar x\|_{\Sx^{-1}}^2
			{}\leq{}
				f(x')
				{}+{}
				\tfrac12\|x-x'\|_{\Sx^{-1}}^2
			\quad
				\forall x'\in\R^n.
			\end{equation}
			Moreover, for every \(x\in\R^n\) it holds that
			\begin{equation}\label{eq:subgradMoreau}
				\Sx^{-1}(x-\bar x)
			{}\in{}
				\hat\partial f(\bar x),
			\quad\text{where}\quad
				\bar x
			{}\in{}
				\prox_f^{\Sx}(x).
			\end{equation}
			If \(\Sx\) is such that \(f+\tfrac12\|{}\cdot{}\|_{\Sx^{-1}}^2\) is strongly convex (in which case \(f\) is said to be \DEF{hypoconvex}), then \(\prox_f^{\Sx}\) is (single-valued and) Lipschitz continuous and \(f^{\Sx}\) Lipschitz differentiable, and \eqref{eq:subgradMoreau} can be strenghtened to
			\begin{equation}\label{eq:gradMoreau}
				\nabla f^{\Sx}(x)
			{}={}
				\Sx^{-1}(x-\bar x)
			{}\in{}
				\hat\partial f(\bar x),
			\quad\text{where}\quad
				\bar x
			{}={}
				\prox_f^{\Sx}(x).
			\end{equation}

		\subsection{Paper outline}

			The remainder of the paper is outlined as follows.
			\Cref{sec:PALM} discusses the theoretical convergence of inexact proximal point iterations on the extended-real-valued cost of \eqref{eq:P}, and shows equivalence between these iterations and the inexact proximal augmented Lagrangian method.
			\review{\Cref{alg:PALM} therein illustrates the proposed (modification of) proximal ALM applied to QPs, thus providing a snapshot of the main steps of the proposed QPALM algorithm pruned of all the implementation details, which will instead be covered in the subsequent sections.
			Specifically, \cref{sec:Subproblem} deals with the inner minimization procedure required at \cref{state:PALM:x}, that is, the semismooth Newton method with exact linesearch of \cite{hermans2019qpalm}.
			\Cref{sec:Code} covers the required heavy duty linear algebra routines, including factorizations and factorization updates, and the LO(B)PCG algorithm used to compute the minimum eigenvalue of \(Q\).
			\Cref{sec:Parameter} lays out in detail the parameters used in QPALM, how they are initialized and updated.
			It furthermore discusses preconditioning of the problem data, as well as termination criteria and infeasibility detection routines.
			After a brief recap on the material of the previous sections, \cref{sec:FullQPALM} presents the fully detailed implementation of QPALM in the dedicated \cref{alg:FullQPALM}, together with a comprehensive overview of all the algorithmic steps therein.}
			\Cref{sec:Simulations} then presents numerical results obtained by comparing this implementation against state-of-the-art solvers.
			Finally, \Cref{sec:Conclusion} draws the concluding remarks of the paper.

	\section{Proximal ALM}\label{sec:PALM}%

		\review{As mentioned in the introduction, our methodology revolves around the interpretation of a P-ALM step on \eqref{eq:P} as an inexact proximal point iteration on the extended-real-valued cost
		\begin{equation}
			\varphi(x)
		{}\coloneqq{}
			\tfrac12\trans xQx+\trans qx+\indicator_C(Ax),
		\end{equation}
		namely
		\begin{equation}\label{eq:prox}
			\hat x^{k+1}
		{}\approx{}
			\prox_\varphi^{\Sx}(\hat x^k)
		{}\coloneqq{}
			\argmin_{x\in\R^n}\set{
				\varphi(x)+\tfrac12\|x-\hat x^k\|_{\Sx^{-1}}^2
			}.
		\end{equation}
		Although differing from the original \eqref{eq:P} only by a quadratic term, similarly to what suggested in \cite{bertsekas1979convexification} by selecting a suitably small weight \(\Sx\in\symm_{++}(\R^n)\) this minimization subproblem can be made strongly convex and addressed by means of an ALM method where the hard constraint \(Ax\in C\) is replaced by a quadratic penalty.
		Note that, in order to make the subproblem strongly convex, a diagonal matrix \(\Sx\) can be chosen based on the minimum eigenvalue \(\lambda_{\rm min}\) of the objective Hessian matrix \(Q\).
		If \(\lambda_{\rm min}\geq0\), corresponding to a convex QP, then any positive definite \(\Sx\) works; otherwise, it is easy to verify that it is enough to select the diagonal elements as \(0<\Sx_{ii}< -\nicefrac{1}{\lambda_{\rm min}}\), and that \(\lambda_{\rm min}\) can be replaced by any (under-)estimation carried out at initialization.
		
		Using ALM to solve the inner PP subproblems gives rise to a modification of the proximal ALM scheme, the difference being that the proximity point \(\hat x^k\) is kept constant for some subsequent iterations until a suitably accurate solution \(\hat x^{k+1}\) of \eqref{eq:prox} has been found.
		This subproblem amounts to the nonsmooth composite minimization}
		\begin{equation}\label{eq:primal}
			\minimize_{x\in\R^n,z\in\R^m} f(x)+g(z)+\tfrac12\|x-\hat x^{k}\|_{\Sx^{-1}}^2
		\quad\stt{}
			Ax-z=0,
		\end{equation}
		where \(f=\tfrac12\innprod{{}\cdot{}}{Q{}\cdot{}}+\innprod{q}{{}\cdot{}}\) and \(g=\indicator_C\), and thus itself requires an iterative procedure.
		\begin{subequations}\label{eq:xyz}%
			Starting from a vector \(y^k\in\R^m\) and for a given dual weight matrix \(\Sy_k\in\symm_{++}(\R^m)\), one iteration of ALM applied to \eqref{eq:primal} produces a triplet \((x^{k+1},z^{k+1},y^{k+1})\) according to the following update rule:
			\begin{equation}
				\begin{cases}[r @{{}={}} l]
					(x^{k+1},z^{k+1})
				&
					\argmin_{x,z}\LL_{\hat x^{k},\Sx,\Sy_k}(x,z,y^k)
				\\[3pt]
					y^{k+1}
				&
					y^{k}
					{}+{}
					\Sy_k(Ax^{k+1}-z^{k+1}),
				\end{cases}
			\end{equation}
			where
			\begin{equation}
				\LL_{\hat x^{k},\Sx,\Sy}(x,z,y)
			{}\coloneqq{}
				f(x)
				{}+{}
				\tfrac12\|x-\hat x^{k}\|_{\Sx^{-1}}^2
				{}+{}
				g(z)
				{}+{}
				\innprod{y}{Ax-z}
				{}+{}
				\tfrac12\|Ax-z\|_{\Sy}^2
			\end{equation}
			is the \(\Sy\)-augmented Lagrangian associated to \eqref{eq:primal}.
		\end{subequations}
		Notice that, by first minimizing with respect to \(z\), apparently \(x^{k+1}\) and \(z^{k+1}\) are given by
		\begin{equation}\label{eq:xyz_explicit}
			\begin{cases}[r @{{}={}} l]
				x^{k+1}
			&
				\argmin_{x\in\R^n}\set{
					f(x)
					{}+{}
					g^{\Sy_k^{-1}}(Ax+\Sy_k^{-1}y^{k})
					{}+{}
					\tfrac12\|x-\hat x^{k}\|_{\Sx^{-1}}^2
				}
			\\
				z^{k+1}
			&
				\prox_g^{\Sy_k^{-1}}(Ax^{k+1}+\Sy_k^{-1}y^{k})
				{}={}
				\proj_C(Ax^{k+1}+\Sy_k^{-1}y^{k}),
			\end{cases}
		\end{equation}
		\review{where the second equality in the \(z\)-update owes to the fact that \(C\) is a box and \(\Sy_k\) is diagonal, so that the projections onto \(C\) with respect to \(\|{}\cdot{}\|\) and \(\|{}\cdot{}\|_{\Sy_k^{-1}}\) coincide.}
		
		\begin{rem}[Proximal ALM vs plain ALM]%
			A major advantage of proximal ALM over plain ALM when applied to a nonconvex QP is that by suitably selecting the proximal weights each subproblem is guaranteed to have solutions.
			An illustrative example showing how ALM may not be applicable is given by the nonconvex QP
			\[
				\minimize_{x\in\R^2}x_1x_2
			\quad
				\stt x_1=0,
			\]
			which is clearly lower bounded and with minimizers given by \(\set{x\in\R^2}[x_1=0]\).
			For a fixed penalty \(\beta>0\) and a Lagrangian multiplier \(y\in\R\), the \(x\)-minimization step prescribed by ALM is
			\[
				x_{\rm ALM}^+
			{}\in{}
				\argmin_{w\in\R^2}\set{
					w_1w_2
					{}+{}
					\innprod{y}{w_1}
					{}+{}
					\tfrac\beta2\|w_1\|^2
				}
			{}={}
				\emptyset,
			\]
			owing to lower unboundedness of the augmented Lagrangian (take, \eg, \(w^k=(1,-k)\) for \(k\to\infty\)).
			The problem is readily solved by proximal ALM, as long as the proximal weight \(\Sx\in\symm_{++}(\R^2)\) satisfies \(\Sx\prec\review{\Id_2}\) (the \(2\times2\)-identity matrix).
			In fact, the P-ALM update step results in
			\begin{Align*}
				x_{\text{\rm P-ALM}}^+
			{}\in{} &
				\argmin_{w\in\R^2}\set{
					w_1w_2
					{}+{}
					\innprod{y}{w_1}
					{}+{}
					\tfrac\beta2\|w_1\|^2
					{}+{}
					\tfrac12\|w-x\|_{\Sx^{-1}}^2
				}
			\\
			{}={} &
				\set{\textstyle
					\Bigl[
						\Sx^{-1}
						{}+{}
						\binom{\beta~~1}{1~~\hphantom\beta}
					\Bigr]^{-1}
			\ifarxiv
				\hspace*{-2pt}
			\fi
					\Bigl(
						\Sx x
						{}-{}
						\binom y0
					\Bigr)
				},
			\end{Align*}
			which is well defined \review{for any Lagrange multiplier \(y\) and penalty \(\beta>0\)}.
		\end{rem}
		
		\review{The resulting proximal ALM is outlined in \cref{alg:PALM}.
		The iterate \(x^{k+1}\) retrieved at \cref{state:PALM:x} corresponds to an approximate solution of the smooth \(x\)-subproblem in \eqref{eq:xyz_explicit} which, in principle, can be carried out by any smooth minimization technique terminating when the norm of the gradient falls within the prescribed tolerance \(\delta_k\).
		The proposed QPALM solver, whose full implementation is outlined in \cref{alg:FullQPALM} in \cref{sec:FullQPALM}, will ultimately address this step with a semismooth Newton method with exact linesearch.
		For the moment being, however, we shall regard this update as a black box to merely focus on the convergence analysis of the outer proximal ALM.
		
		\begin{algorithm}
			\caption{Proximal augmented Lagrangian method for nonconvex QPs}%
			\label{alg:PALM}%

		\review{	
		\begin{algorithmic}[1]
		\Require
			\begin{tabular}[t]{@{}l@{}}
				\((\hat x^0,y^0)\in\R^n\times\R^m\);~~
				\(\delta_0,\varepsilon_0>0\);~~
				\(\rho\in(0,1)\);~~
				\(\Sy_0\succeq\Sy_{\rm min}\in\symm_{++}(\R^m)\)
			\\
				\(\Sx\in\symm_{++}(\R^n)\)
				~such that~
				\(Q+\Sx^{-1}\succ0\)
			\end{tabular}
		\For{ \(k=0,1,\ldots\) }
			\State\label{state:PALM:x}\label{state:PALM:z}%
		\review{%
				\begin{tabular}[t]{@{}l@{}}
					Let \(x^{k+1}\) be such that~
					\(
						\Bigl\|~
							\smash{
								\overbrace*{
									Qx^{k+1}
									{}+{}
									q
								}^{\nabla f(x^{k+1})}
								{}+{}
								\overbrace*{
									\Sy_k(Ax^{k+1}-z^{k+1})+y^k
								}^{\nabla g^{\Sy_k^{-1}}(x^{k+1})}
								{}+{}
								\Sx^{-1}(x^{k+1}-\hat x^k)
							}
						~\Bigr\|
					{}\leq{}
						\delta_k
					\)
				\\
					where
					\(
						z^{k+1}
					{}={}
						\proj_C\bigl(Ax^{k+1}+\Sy_k^{-1}y^k\bigr)
					\)
				\end{tabular}
			\State\label{state:PALM:y}%
		\review{%
				\(
					y^{k+1}
				{}={}
					y^k
					{}+{}
					\Sy_k(Ax^{k+1}-z^{k+1})
				\)
		}%
		\review{%
			\If{ \(\|Ax^{k+1}-z^{k+1}\|_{\Sy_k}\leq\varepsilon_k\) }\label{step:innerTermination}%
				\Comment
					{{\color{ForestGreen}\sf\footnotesize [Quit ALM inner loop]}}
				\State\label{state:y}%
		\review{%
					Update \(\hat x^{k+1}=x^{k+1}\)
					~and choose~
					\(\Sy_{k+1}\succeq\Sy_{\rm min}\)
					and \(\varepsilon_{k+1}\leq\rho\varepsilon_k\)
		}%
		\review{%
			\Else{}
				\State\label{state:PALM:inner}%
		\review{%
					Set~
					\(\hat x^{k+1}=\hat x^k\) and
					\(\varepsilon_{k+1}=\varepsilon_k\),
					~and choose~
					\(\Sy_{k+1}\succeq\Sy_k\)
		}%
		}%
			\EndIf{}
		}%
			\State\label{state:PALM:delta}%
		\review{%
				Choose \(\delta_{k+1}\leq\rho\delta_k\)
		}%
		}%
		\EndFor{}
		\end{algorithmic}
		}
		\end{algorithm}}%

		\subsection{Inexact proximal point}

			We now summarize a key result that was shown in \cite[\S4.1]{sun2017convergence} in the more general setting of proximal gradient iterations.
			Given that \cite{sun2017convergence} has not yet been peer-reviewed and also for the sake of self-containedness, we provide a proof tailored to our simplified setting in the dedicated \cref{proof:thm:PP}.
			
			\begin{thm}[Inexact nonconvex PP {\cite[\S4.1]{sun2017convergence}}]\label{thm:PP}%
				Let \(\func{\varphi}{\R^n}{\Rinf}\) be a proper, lsc and lower bounded function.
				Starting from \(x^0\in\R^n\) and given a sequence \(\seq{e^k}\subset\R^n\) such that \(\sum_{k\in\N}\|e^k\|<\infty\), consider the inexact PP iterations
				\[
					x^{k+1}
				{}\in{}
					\prox_\varphi^{\Sx}(x^k+e^k)
				\]
				for some \(\Sx\in\symm_{++}(\R^n)\).
				Then, the following hold:
				\begin{enumerate}
				\item\label{thm:PP:cost}%
					the real-valued sequence \(\seq{\varphi(x^{k+1})}\) converges to a finite value;
				\item\label{thm:PP:summable}%
					the sequence \(\seq{\|x^{k+1}-x^k\|^2}\) has finite sum, and in particular \(\min_{i\leq k}\|x^{k+1}-x^i\|\leq o(\nicefrac{1}{\sqrt k})\);
				\item\label{thm:PP:omega}%
					\(\varphi\) is constant and equals the limit of \(\seq{\varphi(x^{k+1})}\) on the set of cluster points of \(\seq{x^k}\), which is made of stationary points for \(\varphi\);
				\item\label{thm:PP:bounded}%
					if \(\varphi\) is coercive, then \(\seq{x^k}\) is bounded.
				\end{enumerate}
			\grayout{%
				If moreover \(\varphi\) is semialgebraic, \(\seq{x^k}\) remains bounded, and
				\(
					\smash{
						\sum_{k\in\N}\!\sqrt{\sum_{j\geq k}\|e^j\|^2}
						{}<{}
						\infty
					}
				\)
				(as is the case when \(\|e^k\|\sim O(\rho^k)\) for some \(\rho\in(0,1)\)), then
				\begin{enumerate}[resume]
				\item\label{thm:PP:global}
					 the entire sequence \(\seq{x^k}\) converges to a (stationary) point \(x_\star\).
				\end{enumerate}
			}%
				\begin{proof}
					See \cref{proof:thm:PP}.
				\end{proof}
			\end{thm}
			
			We remark that with trivial modifications of the proof the arguments also apply to time-varying proximal weights \(\Sx_k\), \(k\in\N\), as long as there exist \(\Sx_{\rm min},\Sx_{\rm max}\in\symm_{++}(\R^n)\) such that \(\Sx_{\rm min}\preceq\Sx_k\preceq\Sx_{\rm max}\) holds for all \(k\).
			\Cref{thm:PP:bounded} indicates that coerciveness of the cost function is a sufficient condition for inferring boundedness of the iterates.
			In \eqref{eq:P}, however, the cost function \(\varphi\) may fail to be coercive even if lower bounded on the feasible set \(\set{x}[Ax\in C]\).
			This happens when there is a feasible direction for which the objective is constant, \ie when \(\lim_{\tau\to\infty}\indicator_C(A(x+\tau d)) = 0\), \(Qd = 0\) and \(\trans q d = 0\) hold for some \(x,d\in\R^n\) with \(d\neq 0\).
			Nevertheless, it has been shown in
			\cite{luo1993error}
			that (exact) proximal point iterations on a lower bounded nonconvex quadratic program remain bounded and, in fact, converge to a stationary point.
			We next show that this remains true even for inexact proximal point iterations, provided that the inexactness vanishes at linear rate.
			The proof hinges on the (exact) proximal gradient error bound analysis of \cite{luo1993error} and on the close relation existing among proximal point and proximal gradient iterations for this kind of problems.
			Before showing the result in \cref{thm:QP}, we present a simple technical lemma that will be needed in the proof.
			
			\begin{lem}\label{thm:Lip}%
				Let \(h\) be a lower bounded and \(L_h\)-smooth function.
				Then, for every \(\alpha>0\) it holds that
				\[
					\lip_{\lev_{\leq\inf h+\alpha}h}h
				{}\leq{}
					\tfrac{1+\sqrt2}{2}
					\sqrt{2\alpha L_h}.
				\]
				\begin{proof}
					Without loss of generality we may assume that \(\inf h=0\).
					Let \(\alpha>0\) be fixed, and consider \(x,y\in\lev_{\leq\alpha}h\) with \(x\neq y\).
					From the quadratic upper bound of Lipschitz differentiable functions (see \eg \cite[Prop. A.24]{bertsekas2016nonlinear}) and the fact that \(0\leq h(x),h(y)\leq\alpha\) we have that
					\begin{Align*}
					\numberthis\label{eq:Lip}
						\frac{|h(y)-h(x)|}{\|y-x\|}
					{}\leq{} &
						\min\set{
							\tfrac{\alpha}{\|y-x\|},\,
							\|\nabla h(x)\|+\tfrac{L_h}{2}\|y-x\|
						}
					\\
					{}\leq{} &
						\min\set{
							\tfrac{\alpha}{\|y-x\|},\,
							\sqrt{2\alpha L_h}+\tfrac{L_h}{2}\|y-x\|
						},
					\end{Align*}
					where the last inequality follows from the fact that
					\[
						\alpha
					{}\geq{}
						h(x)-h(x-\tfrac{1}{L_h}\nabla h(x))
					{}\geq{}
						\tfrac{1}{2L_h}\|\nabla h(x)\|^2,
					\]
					where the last equality again uses the quadratic lower bound of \cite[Prop. A.24]{bertsekas2016nonlinear}.
					By solving a second-order equation in \(\|y-x\|\), we see that
					\begin{align*}
						\min\set{
							\tfrac{\alpha}{\|y-x\|},
							\sqrt{2\alpha L_h}+\tfrac{L_h}{2}\|y-x\|
						}
					{}={} &
						\begin{ifcases}
							\tfrac{\alpha}{\|y-x\|}
						&
							\|y-x\|
							{}\geq{}
							(2-\sqrt2)\sqrt{\frac{\alpha}{L_h}}
						\\
							\sqrt{2\alpha L_h}+\tfrac{L_h}{2}\|y-x\|
						\otherwise,
						\end{ifcases}
					\\
					{}\leq{} &
						\tfrac{\alpha}{2-\sqrt2}
						\sqrt{\tfrac{L_h}{\alpha}}
					{}={}
						\tfrac{1+\sqrt2}{2}\sqrt{2\alpha L_h},
					\end{align*}
					resulting in the claimed bound.
				\end{proof}
			\end{lem}
			
			\begin{rem}
				By discarding the term \(\frac{\alpha}{\|y-x\|}\) in \eqref{eq:Lip} and letting \((y,x)\to(\bar x,\bar x)\) with \(x\neq y\), one obtains that the \emph{pointwise} Lipschitz constant of a lower bounded and \(L_h\)-smooth function \(h\) can be estimated as \(\lip h(\bar x)\leq\sqrt{2(h(\bar x)-\inf h)L_h}\).
				Therefore, if \(h\) is also (quasi-)convex \cref{thm:Lip} can be tightened to \(\lip_{\lev_{\leq\inf h+\alpha}h}h\leq\sqrt{2\alpha L_h}\), owing to convexity of the sublevel set together with \cite[Thm. 9.2]{rockafellar2011variational}.
			\end{rem}
			
			\begin{thm}[Linear convergence of inexact PP on nonconvex quadratic programs]\label{thm:QP}%
				Let \(\varphi=f+\indicator_\Omega\), where \(\Omega\subseteq\R^n\) is a nonempty polyhedral set and \(\func{f}{\R^n}{\R}\) is a (possibly nonconvex) quadratic function which is lower bounded on \(\Omega\).
				Starting from \(x^0\in\R^n\) and given a sequence \(\seq{e^k}\subset\R^n\) such that \(\|e^k\|\in O(\rho^k)\) for some \(\rho\in(0,1)\), consider the inexact PP iterations
				\[
					x^{k+1}
				{}\in{}
					\prox_\varphi^{\Sx}(x^k+e^k)
				\]
				for some \(\Sx\in\symm_{++}(\R^n)\).
				Then, the sequence \(\seq{x^k}\) converges at \(R\)-linear rate to a stationary point of \(\varphi\).
				\begin{proof}
					Let \(x_\star^k\) be a projection of \(x^k\) onto the set \(\zer\partial\varphi\) of stationary points for \(\varphi\).
					Such a point exists for every \(k\) owing to nonemptiness and closedness of \(\zer\partial\varphi\), the former condition holding by assumption and the latter holding because of closedness of \(\graph\partial\varphi\), cf. \cite[Prop. 8.7 and Thm. 5.7(a)]{rockafellar2011variational}.
					From \cite[Eq.s (2.1) and (A.3)]{luo1993error}, which can be invoked owing to \cite[Thm. 2.1(b)]{luo1993error}, it follows that there exists \(\tau>0\) such that
					\begin{equation}\label{eq:subreg}
						\varphi(x_\star^k)=\varphi_\star
					\quad\text{and}\quad
						\dist(x^k,\zer\partial\varphi)
					{}\leq{}
						\tau\|x^k-\proj_\Omega^{\Sx}(x^k-\Sx\nabla f(x^k))\|_{\Sx^{-1}}
					\end{equation}
					hold for \(k\) large enough, where \(\proj_\Omega^{\Sx}=\prox_{\indicator_\Omega}^{\Sx}\) is the projection with respect to the distance \(\|{}\cdot{}\|_{\Sx^{-1}}\).
					Let \(L_f\) and \(L_{\varphi^{\Sx}}\) be Lipschitz constants for \(\nabla f\) and \(\nabla\varphi^{\Sx}\), respectively.
					Note that stationarity of \(x_\star^k\) implies that \(\varphi(x_\star^k)=\varphi^{\Sx}(x_\star^k)\) and \(\nabla\varphi^{\Sx}(x_\star^k)=0\).
					We have
					\begin{align*}
						\varphi^{\Sx}(x^k)
						{}-{}
						\varphi_\star
					{}={} &
						\varphi^{\Sx}(x^k)
						{}-{}
						\varphi^{\Sx}(x_\star^k)
					{}\leq{}
						\tfrac{L_{\varphi^{\Sx}}}{2}\|x^k-x_\star^k\|^2
					{}={}
						\tfrac{L_{\varphi^{\Sx}}}{2}\dist(x^k,\zer\partial\varphi)^2
					\\
					\numberthis\label{eq:QPcostk}
					{}\overrel[\leq]{\eqref{eq:subreg}}{} &
						\tfrac{L_{\varphi^{\Sx}}\tau^2}{2}
						\|x^k-\proj_\Omega^{\Sx}(x^k-\Sx\nabla f(x^k))\|_{\Sx^{-1}}^2.
					\end{align*}
					Next, observe that
					\begin{align*}
						x^{k+1}
					{}={}
						\prox_\varphi^{\Sx}(x^k+e^k)
					~\Leftrightarrow~ &
						\Sx^{-1}(x^k+e^k-x^{k+1})
					{}\in{}
						\partial\varphi(x^{k+1})
					{}={}
						\nabla f(x^{k+1})+\partial\indicator_\Omega(x^{k+1})
					\\
					~\Leftrightarrow~ &
						\Sx^{-1}(x^k+e^k-\Sx\nabla f(x^{k+1})-x^{k+1})
					{}\in{}
						\partial\indicator_\Omega(x^{k+1})
					\\
					~\Leftrightarrow~ &
						x^{k+1}
					{}={}
						\prox_{\indicator_\Omega}^{\Sx}[x^k+e^k-\Sx\nabla f(x^{k+1})]
					\\
					\numberthis\label{eq:PP=proj}
					~\Leftrightarrow~ &
						x^{k+1}
					{}={}
						\proj_\Omega^{\Sx}[x^k+e^k-\Sx\nabla f(x^{k+1})].
					\end{align*}
					Denoting \(c\coloneqq L_{\varphi^{\Sx}}\tau^2\), we obtain that 
					\begin{align*}
						\varphi^{\Sx}(x^k)
						{}-{}
						\varphi_\star
					{}\overrel[\leq]{\eqref{eq:QPcostk}}{} &
						\mathtight[0.01]
						c
						\|x^k-x^{k+1}\|_{\Sx^{-1}}^2
						{}+{}
						c
						\bigl\|
							\overbracket*{
								\proj_\Omega^{\Sx}[x^k+e^k-\Sx\nabla f(x^{k+1})]
							}^{x^{k+1}}
							{-}
							\proj_\Omega^{\Sx}[x^k-\Sx\nabla f(x^k)]
						\bigr\|_{\Sx^{-1}}^2
					\shortintertext{%
						which, using 1-Lipschitz continuity of \(\proj^{\Sx}\) in the norm \(\|{}\cdot{}\|_{\Sx^{-1}}\),
					}
					{}\leq{} &
						c
						\|x^k-x^{k+1}\|_{\Sx^{-1}}^2
						{}+{}
						c
						\bigl\|
							[x^k+e^k-\Sx\nabla f(x^{k+1})]
							{}-{}
							[x^k-\Sx\nabla f(x^k)]
						\bigr\|_{\Sx^{-1}}^2
					\\
					{}\leq{} &
						c
						\|x^k-x^{k+1}\|_{\Sx^{-1}}^2
						{}+{}
						2c
						\|e^k\|_{\Sx^{-1}}^2
						{}+{}
						2c
						\bigl\|
							\nabla f(x^k)
							{}-{}
							\nabla f(x^{k+1})
						\bigr\|_{\Sx}^2
					\\
					{}\leq{} &
						(c+2cL_f^2\|\Sx\|)
						\|x^k-x^{k+1}\|_{\Sx^{-1}}^2
						{}+{}
						2c\|e_k\|_{\Sx^{-1}}^2
					\\
					\numberthis\label{eq:inexactEB}
					{}\leq{} &
						c_1
						\|x^k-x^{k+1}\|_{\Sx^{-1}}^2
						{}+{}
						c_2\rho^{2k}
					\end{align*}
					for some constants \(c_1,c_2>0\).
					Observe that
					\begin{align*}
						\varphi^{\Sx}(x^{k+1})
					{}\leq{}
						\varphi(x^{k+1})
					{}={} &
						\varphi^{\Sx}(x^k+e^k)
						{}-{}
						\tfrac12\|x^{k+1}-x^k-e^k\|_{\Sx^{-1}}^2
					\\
					{}\leq{} &
						\varphi^{\Sx}(x^k)
						{}+{}
						L\|e^k\|
						{}-{}
						\tfrac14\|x^{k+1}-x^k\|_{\Sx^{-1}}^2
						{}+{}
						\tfrac12\|e^k\|_{\Sx^{-1}}^2
					\\
					\numberthis\label{eq:inexactSDMoreau}
					{}\leq{} &
						\varphi^{\Sx}(x^k)
						{}-{}
						\tfrac14\|x^{k+1}-x^k\|_{\Sx^{-1}}^2
						{}+{}
						c_3\rho^k
					\end{align*}
					for some constant \(c_3>0\), where in the second inequality \(L\) denotes a Lipschitz constant of the smooth function \(\varphi^{\Sx}\) on a sublevel set that contains all iterates; the existence of such an \(L\) is guaranteed by \cref{thm:PP:cost,thm:Lip}, since \(-\infty<\inf\varphi\leq\varphi^{\Sx}\leq\varphi\).
					Therefore,
					\[
						\bigl(
							\varphi^{\Sx}(x^k)
							{}-{}
							\varphi_\star
						\bigr)
						{}-{}
						\bigl(
							\varphi^{\Sx}(x^{k+1})
							{}-{}
							\varphi_\star
						\bigr)
					{}\overrel[\geq]{\eqref{eq:inexactSDMoreau}}{}
						\tfrac14\|x^{k+1}-x^k\|_{\Sx^{-1}}^2
						{}-{}
						c_3\rho^k
					{}\overrel[\geq]{\eqref{eq:inexactEB}}{}
						\tfrac{1}{4c_1}
						\bigl(
							\varphi^{\Sx}(x^k)
							{}-{}
							\varphi_\star
						\bigr)
						{}-{}
						c_4\rho^k
					\]
					holds for some constant \(c_4>0\).
					By possibly enlarging \(c_1\) we may assume without loss of generality that \(\rho\geq1-\nicefrac{1}{4c_1}\), so that
					\begin{align*}
						\bigl(
							\varphi^{\Sx}(x^{k+1})
							{}-{}
							\varphi_\star
						\bigr)
					{}\leq{} &
						\rho
						\bigl(
							\varphi^{\Sx}(x^k)
							{}-{}
							\varphi_\star
						\bigr)
						{}+{}
						c_4\rho^k
					\\
					{}\leq{} &
						\rho^{k+1}
						\bigl(
							\varphi^{\Sx}(x^0)
							{}-{}
							\varphi_\star
						\bigr)
						{}+{}
						c_4\sum_{j=0}^k{
							\rho^{k-j}
							\rho^j
						}
					\\
					\numberthis\label{eq:costRlinear0}
					{}={} &
						\bigl(
							\rho(\varphi^{\Sx}(x^0)-\varphi_\star)
							{}+{}
							c_4(k+1)
						\bigr)
						\rho^k
					{}\leq{}
						c_5(\sqrt{\rho})^k,
					\end{align*}
					where \(c_5\) is any such that
					\(
						\bigl(
							\rho(\varphi^{\Sx}(x^0)-\varphi_\star)
							{}+{}
							c_4(k+1)
						\bigr)
						(\sqrt{\rho})^k
					{}\leq{}
						c_5
					\)
					holds for every \(k\in\N\).
					Next, denoting
					\(
						\varphi_k
					{}\coloneqq{}
						\varphi^{\Sx}(x^k)
						{}+{}
						\tfrac{c_3}{1-\rho}\rho^k
					\)
					observe that
					\begin{equation}\label{eq:exactSDMoreau}
						\varphi_{k+1}
					{}\leq{}
						\varphi_k
						{}-{}
						\tfrac14\|x^{k+1}-x^k\|_{\Sx^{-1}}^2
					\end{equation}
					as it follows from \eqref{eq:inexactSDMoreau}, and that \(\varphi_\star<\varphi_k\to\varphi_\star\) as \(k\to\infty\).
					In fact, \eqref{eq:costRlinear0} implies that
					\begin{equation}\label{eq:costRlinear}
						0
					{}\leq{}
						\varphi_k-\varphi_\star
					{}\leq{}
						c_6\rho^{\nicefrac k2}
					\end{equation}
					holds for some \(c_6>0\) and all \(k\in\N\).
					Therefore,
					\begin{align*}
						\sum_{j\geq k}\|x^{j+1}-x^j\|_{\Sx^{-1}}
					{}\overrel[\leq]{\eqref{eq:exactSDMoreau}}{} &
						2\sum_{j\geq k}\sqrt{
							\varphi_j-\varphi_{j+1}
						}
					{}\overrel[\leq]{\eqref{eq:costRlinear}}{}
						2
						\sum_{j\geq k}\sqrt{
							\varphi_j-\varphi_\star
						}
					\\
					{}\leq{} &
						2c_6^{\nicefrac12}
						\sum_{j\geq k}{
							\rho^{(j-1)/4}
						}
					{}\leq{}
						c_7
						\rho^{k/4}
					\end{align*}
					for some constant \(c_7>0\).
					In particular, the sequence \(\seq{x^k}\) has finite length and thus converges to a point \(x^\star\), which is stationary for \(\varphi\) owing to \cref{thm:PP:omega}.
					In turn, the claimed \(R\)-linear convergence follows from the inequality
					\(
						\|x^k-x^\star\|_{\Sx^{-1}}
					{}\leq{}
						\sum_{j\geq k}\|x^{j+1}-x^j\|_{\Sx^{-1}}
					\).
				\end{proof}
			\end{thm}

		\subsection{Convergence of \texorpdfstring{\Cref{alg:PALM}}{Algorithm \ref*{alg:PALM}} for nonconvex QPs}%

			\begin{thm}\label{thm:QPALM}
				Suppose that problem \eqref{eq:P} is lower bounded, and consider the iterates generated by \cref{alg:PALM} with \(f(x)=\tfrac12\trans xQx+\trans qx\) and \(g(z)=\indicator_C(z)\).
				Then, the following hold:
				\begin{enumerate}
				\item\label{thm:QPALM:inner}%
					The triplet \((x^{k+1},y^{k+1},z^{k+1})\) produced at the \(k\)-th iteration satisfies
					\[
						\|\nabla f(x^{k+1})+\trans Ay^{k+1}\|\leq\delta_k+\|\Sx^{-1}(x^{k+1}-\hat x^{k})\|
					\quad\text{and}\quad
						y^{k+1}\in\partial g(z^{k+1}).
					\]
				\item\label{thm:QPALM:termination}%
					The condition at \cref{step:innerTermination} is satisfied infinitely often, and \(\|\hat x^{k+1}-\hat x^{k}\|\to0\) as \(k\to\infty\).
					In particular, for every primal-dual tolerances \(\epsilon_{\rm p},\epsilon_{\rm d}>0\), the termination criteria
					\[
						\|\nabla f(x^{k+1})+\trans Ay^{k+1}\|\leq\epsilon_{\rm d}
					\qquad
						y^{k+1}\in\partial g(z^{k+1})
					\qquad
						\|Ax^{k+1}-z^{k+1}\|\leq\epsilon_{\rm p}
					\]
					are satisfied in a finite number of iterations.
				\item\label{thm:QPALM:global}%
					The sequence \(\seq{\hat x^k}\) converges to a stationary point of problem \eqref{eq:P}; in fact, denoting \(\seq{k_i}[i\in\N]\) as the (infinite) set of those indices at which the condition at \cref{step:innerTermination} is satisfied, the sequence \(\seq{x^{k_i+1}}[i\in\N]\) converges at \(R\)-linear rate.
				\end{enumerate}
				\begin{proof}
					The definition of \(z^{k+1}\) at \cref{state:PALM:z} and the characterization of \(\prox_g^{\smash{\Sy_k^{-1}}}\) yield
					\[
						\partial g(z^{k+1})
					{}\ni{}
						\nabla g^{\Sy_k^{-1}}(Ax^{k+1}+\Sy_k^{-1}y^{k})
					{}={}
						\Sy_k\bigl(
							Ax^{k+1}+\Sy_k^{-1}y^{k}
							{}-{}
							z^{k+1}
						\bigr)
					{}={}
						y^{k+1}.
					\]
					By expanding the gradient appearing in the norm at \cref{state:PALM:x}, we thus have
					\begin{equation} \label{eq:grad phi_j}
						\delta_k
					{}\geq{}
						\|
							\nabla f(x^{k+1})
							{}+{}
							\trans Ay^{k+1}
							{}+{}
							\Sx(x^{k+1}-\hat x^k)
						\|,
					\end{equation}
					and assertion \ref{thm:QPALM:inner} follows from the triangular inequality.
					Next, observe that whenever the condition at \cref{step:innerTermination} is not satisfied the variable \(\hat x^{k+1}\) is not updated (cf. \cref{state:PALM:inner}), and thus \cref{state:PALM:x,state:PALM:y,state:PALM:z} amount to ALM iterations applied to the \review{convex problem}
					\[
						\minimize_{x\in\R^n,z\in\R^m}f(x)+\tfrac12\|x-\hat x^{k}\|_{\Sx^{-1}}^2+g(z)
					\quad\stt{}
						Ax=z
					\]
					with a summable inexactness in the computation of the \(x\)-minimization step.
					The existence of dual solutions entailed by the strong duality of convex QPs guarantees through \cite[Thm. 4 and \S6]{rockafellar1976augmented} that the feasibility residual vanishes, hence that eventually \cref{step:innerTermination} holds.
					
					Let
					\(
						d^k
					{}\coloneqq{}
						\nabla f(x^{k+1})
						{}+{}
						\trans Ay^{k+1}
						{}+{}
						\Sx(x^{k+1}-\hat x^k)
					\)
					be the gradient appearing in the norm at \cref{state:PALM:x} and let \(e^k\coloneqq Ax^{k+1}-z^{k+1}\).
					Let \(\seq{k_i}[i\in\N]\) be the (infinite) set of all indices at which the condition at \cref{step:innerTermination} is satisfied, so that \(\hat x^{k_i+1}=x^{k_i+1}\) and \(\|e^{k_i}\|\leq\varepsilon_{k_i}\leq\rho^i\varepsilon_0\).
					Then, for every \(i\in\N\)
					\[
						\begin{cases}[r>{{}}c<{{}}l]
							0 & = & \nabla f(x^{k_i+1})+\Sx^{-1}\bigl(x^{k_i+1} - (x^{k_{i-1}+1}+\Sx d^{k_i})\bigr)+\trans Ay^{k_i+1}\\
							0 &\in& \partial g(z^{k_i+1})-y^{k_i+1}\\
							0 & = & Ax^{k_i+1}-z^{k_i+1}-e^{k_i}.
						\end{cases}
					\]
					In particular, \((x^{k_i+1},z^{k_i+1},y^{k_i+1})\) is a primal-dual solution of
					\[
						\minimize_{x\in\R^n,z\in\R^m}f(x)+g(z)+\tfrac12\|x-(x^{k_{i-1}+1}+\Sx d^{k_i})\|_{\Sx^{-1}}^2
					\quad\stt{}
						Ax-z=e^{k_i}.
					\]
					Therefore, denoting \(\func{\mathcal X}{\dom\mathcal X\subseteq\R^n\times\R^m}{\R^n}\) as the operator
					\[
						\mathcal X(u,v)
					{}={}
						\argmin_{x\in\R^n}\set{
							\tfrac12\trans xQx+\trans qx+\tfrac12\|x-u\|_{\Sx^{-1}}^2
						}[Ax-v\in C],
					\]
					we have that
					\(
						x^{k_i+1}
					{}={}
						\mathcal X(\hat x^{k_{i-1}+1}+\Sx d^{k_i}, e^{k_i})
					\).
					Notice further that \(\mathcal X(u,0)=\prox_\varphi^{\Sx}(u)\) for the QP function \(\varphi(x)=\tfrac12\trans xQx+\trans qx+\indicator_C(Ax)\).
					As shown in \cite[Thm. 1]{patrinos2010new}, \(\mathcal X\) is a polyhedral mapping, and as it is at most single valued (owing to strong convexity of the QP) we deduce from \cite[Cor. 3D.5]{dontchev2009implicit} that it is globally Lipschitz continuous on its (polyhedral) domain with constant, say, \(L\).
					Therefore,
					\begin{align*}
						\|x^{k_i+1}-\prox_\varphi^{\Sx}(x^{k_{i-1}+1})\|^2
					{}={} &
						\|\mathcal X(x^{k_{i-1}+1}+\Sx d^{k_i}, e^{k_i})-\mathcal X(\hat x^{k_{i-1}+1}, 0)\|^2
					\\
					{}\leq{} &
						L^2\left(
							\|\Sx d^{k_i}\|^2+\|e^{k_i}\|^2
						\right)
					\\
					{}\leq{} &
						L^2\|\Sx\|^2\delta_{k_i}^2
						{}+{}
						L^2\|\Sy_{\rm min}\|^{-1}\varepsilon_{k_i}^2
					\\
					{}\leq{} &
						c\rho^i
					\end{align*}
					for some constant \(c>0\) that only depends on the problem and on the algorithm initialization.
					Denoting \(\xi^i\coloneqq x^{k_i+1}\) as the \(i\)-th ``outer'' iterate, this shows that \(\seq{\xi^i}[i\in\N]\) is generated by an inexact proximal point algorithm on function \(\varphi\) with error \(\|e^i\|\leq O(\rho^i)\), namely,
					\[
						\xi^{i+1}
					{}={}
						\prox_\varphi^{\Sx}(\xi^i+e^i).
					\]
					In particular, all the assertions follow from \cref{thm:PP,thm:QP}.
			\grayout{%
					\color{red}%
					
					We now prove linear convergence of the outer sequence \(\seq{z^i}[i\in\N]\).
					Let \(z_\star^i\) be a projection of \(z^i\) onto the set \(\zer\partial\varphi\) of stationary points for \(\varphi\).
					From \cite[Eq.s (2.1) and (A.3)]{luo1993error}, which can be invoked owing to \cite[Thm. 2.1(b)]{luo1993error}, it follows that there exists \(\tau>0\) such that
					\begin{equation}\label{eq:subreg}
						\varphi(z_\star^i)=\varphi_\star
					\quad\text{and}\quad
						\dist(z^i,\zer\partial\varphi)
					{}\leq{}
						\tau\|z^i-\proj_\Omega^{\Sx}(z^i-\Sx\nabla f(z^i))\|_{\Sx^{-1}}
					\end{equation}
					hold for \(k\) large enough, where \(\proj_\Omega^{\Sx}=\prox_{\indicator_\Omega}^{\Sx}\) is the projection with respect to the distance \(\|{}\cdot{}\|_{\Sx^{-1}}\).
					Let \(L_f\) and \(L_{\varphi^{\Sx}}\) be Lipschitz constants for \(\nabla f\) and \(\nabla\varphi^{\Sx}\), respectively.
					Note that stationarity of \(z_\star^i\) implies that \(\varphi(z_\star^i)=\varphi^{\Sx}(z_\star^i)\) and \(\nabla\varphi^{\Sx}(z_\star^i)=0\).
					We have
					\begin{align*}
						\varphi^{\Sx}(z^i)
						{}-{}
						\varphi_\star
					{}={} &
						\varphi^{\Sx}(z^i)
						{}-{}
						\varphi^{\Sx}(x_\star^i)
					{}\leq{}
						\tfrac{L_{\varphi^{\Sx}}}{2}\|z^i-z_\star^i\|^2
					{}={}
						\tfrac{L_{\varphi^{\Sx}}}{2}\dist(z^i,\zer\partial\varphi)^2
					\\
					{}\overrel[\leq]{\eqref{eq:subreg}}{} &
						\tfrac{L_{\varphi^{\Sx}}\tau^2}{2}
						\|z^i-\proj_\Omega^{\Sx}(z^i-\Sx\nabla f(z^i))\|_{\Sx^{-1}}^2.
					\end{align*}
					Next, observe that \(z^{i+1}\) is characterized by
					\(
						z^{i+1}
					{}={}
						\proj_\Omega^{\Sx}[z^i+e^i-\Sx\nabla f(z^{i+1})]
					\),
					and from Lipschitz continuity of the projection we obtain
					\begin{align*}
						\varphi^{\Sx}(z^i)
						{}-{}
						\varphi_\star
					{}\leq{} &
						c
						\|z^i-z^{i+1}\|_{\Sx^{-1}}^2
						{}+{}
						c
						\bigl\|
							[z^i+e^i-\Sx\nabla f(z^{i+1})]
							{}-{}
							[z^i-\Sx\nabla f(z^i)]
						\bigr\|_{\Sx^{-1}}^2
					\\
					{}\leq{} &
						c
						\|z^i-z^{i+1}\|_{\Sx^{-1}}^2
						{}+{}
						2c
						\|e_i\|_{\Sx^{-1}}^2
						{}+{}
						2c
						\bigl\|
							\nabla f(z^i)
							{}-{}
							\nabla f(z^{i+1})
						\bigr\|_{\Sx}^2
					\\
					{}\leq{} &
						(c+2cL_f^2\|\Sx\|)
						\|z^i-z^{i+1}\|_{\Sx^{-1}}^2
						{}+{}
						2c\|e_i\|_{\Sx^{-1}}^2
					\\
					\numberthis\label{eq:inexactEB}
					{}\leq{} &
						c_1
						\|z^i-z^{i+1}\|_{\Sx^{-1}}^2
						{}+{}
						c_2\rho^{2i}.
					\end{align*}
					Observe that
					\begin{align*}
						\varphi^{\Sx}(z^{i+1})
					{}\leq{} &
						\varphi(z^{i+1})
					{}={}
						\varphi^{\Sx}(z^i+e^i)
						{}-{}
						\tfrac12\|z^{i+1}-z^i-e^i\|_{\Sx^{-1}}^2
					\\
					{}\leq{} &
						\varphi^{\Sx}(z^i)
						{}+{}
						L\|e^i\|
						{}-{}
						\tfrac14\|z^{i+1}-z^i\|_{\Sx^{-1}}^2
						{}+{}
						\tfrac12\|e^i\|_{\Sx^{-1}}^2
					\\
					\numberthis\label{eq:inexactSDMoreau}
					{}\leq{} &
						\varphi^{\Sx}(z^i)
						{}-{}
						\tfrac14\|z^{i+1}-z^i\|_{\Sx^{-1}}^2
						{}+{}
						c_3\rho^i,
					\end{align*}
					where in the second inequality \(L\) denotes a Lipschitz constant of the smooth function \(\varphi^{\Sx}\) on a sublevel set that contains all iterates; the existence of such an \(L\) is guarateed by \cref{thm:PP:cost,thm:Lip}, having \(-\infty<\inf\varphi\leq\varphi^{\Sx}\leq\varphi\).
					Therefore,
					\[
						\bigl(
							\varphi^{\Sx}(z^i)
							{}-{}
							\varphi_\star
						\bigr)
						{}-{}
						\bigl(
							\varphi^{\Sx}(z^{i+1})
							{}-{}
							\varphi_\star
						\bigr)
					{}\overrel[\geq]{\eqref{eq:inexactSDMoreau}}{}
						\tfrac14\|z^{i+1}-z^i\|_{\Sx^{-1}}^2
						{}-{}
						c_3\rho^i
					{}\overrel[\geq]{\eqref{eq:inexactEB}}{}
						\tfrac{1}{4c_1}
						\bigl(
							\varphi^{\Sx}(z^i)
							{}-{}
							\varphi_\star
						\bigr)
						{}-{}
						c_4\rho^i.
					\]
					Denoting \(\bar\rho\coloneqq\max\set{\rho,\,1-\nicefrac{1}{4c_1}}<1\), we have
					\begin{align*}
						\bigl(
							\varphi^{\Sx}(z^{i+1})
							{}-{}
							\varphi_\star
						\bigr)
					{}\leq{} &
						\bar\rho
						\bigl(
							\varphi^{\Sx}(z^i)
							{}-{}
							\varphi_\star
						\bigr)
						{}+{}
						c_4\bar\rho^i
					\\
					{}\leq{} &
						\bar\rho^{i+1}
						\bigl(
							\varphi^{\Sx}(z^0)
							{}-{}
							\varphi_\star
						\bigr)
						{}+{}
						c_4\sum_{j=0}^i{
							\bar\rho^{i-j}
							\bar\rho^j
						}
					\\
					{}={} &
						\bigl(
							\bar\rho(\varphi^{\Sx}(z^0)-\varphi_\star)
							{}+{}
							c_4i
						\bigr)
						\bar\rho^i
					{}\leq{}
						c_5(\sqrt{\bar\rho})^{i},
					\end{align*}
					where \(c_5\) is any such that
					\(
						\bigl(
							\bar\rho(\varphi^{\Sx}(z^0)-\varphi_\star)
							{}+{}
							c_4i
						\bigr)
						(\sqrt{\bar\rho})^i
					{}\leq{}
						c_5
					\)
					for every \(i\in\N\).
				}%
				\end{proof}
			\end{thm}

	\section{Subproblem minimization}\label{sec:Subproblem}

		The previous section outlined the overall strategy employed by QPALM, the proximal augmented Lagrangian method.
		This section describes our approach to the inner minimization in \cref{state:PALM:x}, which is clearly the most computationally expensive step of \cref{alg:PALM}.
		QPALM uses an iterative method to solve the convex unconstrained optimization of \eqref{eq:xyz_explicit}, computing a semismooth Newton direction and the optimal stepsize at every iteration.
		Given the convex nature of the inner subproblem, the method we propose here has not changed from \cite{hermans2019qpalm}.

		\subsection{Semismooth Newton method}

			Let \(\varphi_k(x) \) denote the objective function of \eqref{eq:xyz_explicit}.
			This can be written as 
			\begin{align*}
				\varphi_k(x) 
			{}={} &
				f(x) 
				{}+{} 
				g^{\Sy_k^{-1}}(Ax+\Sy_k^{-1} y^{k})
				{}+{}
				\tfrac12\|x-\hat x^k\|_{\Sx^{-1}}^2 
			\\
			{}={} &
				f(x) 
				{}+{} 
				\dist^2_{\Sy_k}\bigl(Ax+\Sy_k^{-1} y^{k}\bigr) 
				{}+{} 
				\tfrac12\|x-\hat x^k\|_{\Sx^{-1}}^2,
			\end{align*}
			and its gradient is given by 
			\begin{align*}
				\nabla\varphi_k(x) 
			{}={}
				\nabla f(x) 
				{}+{} 
				\trans A ( y^{k} + \Sy_k(Ax-Z_k(x))) 
				{}+{}
				\Sx^{-1}(x-\hat x^k),
			\end{align*}
			with  
			\begin{align*}
				Z_k(x)
				&
			{}={}
				\prox_g^{\Sy_k^{-1}}(A x+\Sy_k^{-1} y^{k})
			{}={}
				\proj_C(A x+\Sy_k^{-1} y^{k}) 
				\\
				&
			{}={}	
					Ax+\Sy_k^{-1} y^{k}
					{}+{}
					[\ell-Ax-\Sy_k^{-1} y^{k}]_+
					{}-{}
					[Ax+\Sy_k^{-1} y^{k}-u]_+.
			\end{align*}
			Note that this gradient also appears in \eqref{eq:grad phi_j}, with trial point \(\tilde y^{k+1} =  y^{k} + \Sy_k(Ax-Z_k(x))\).
			Furthermore, because of the projection operator in \(Z_k\), the gradient is not continuously differentiable.
			However, we can use the generalized Jacobian \cite[\S7.1]{facchinei2003finite} of \(\proj_C\) at \(Ax+\Sy_k^{-1} y^{k}\), one element of which is the diagonal matrix \(P_k(x)\) with entries
			\[
				(P_k(x))_{ii}
			{}={}
				\begin{ifcases}
					1 & \ell_i\leq (Ax+\Sy_k^{-1} y^{k})_i\leq u_i\\
					0 \otherwise,
				\end{ifcases}
			\]
			see \eg \cite[\S6.2.d]{themelis2018acceleration}.
			Therefore, one element of the generalized Hessian of \(\varphi_k\) is
			\[
				H_k(x)
			{}={}
				Q+\trans A\Sy_k(\review{\Id}-P_k(x))A
				{}+{}
				\Sx^{-1}.
			\]
			Denoting the set of \emph{active} constraints as 
			\begin{equation}\label{eq:J}
				\mathcal J_k(x)
			{}\coloneqq{}
				\set{i}[
					(Ax+\Sy_k^{-1} y^{k})_i
				{}\notin{}
					{[\ell_i,u_i]}
				],
			\end{equation}
			one has that \((\review{\Id}-P_k(x))_{ii}\) is 1 if \(i\in\J_k(x)\) and 0 otherwise.
			In the remainder of the paper, when \(\J_k(x)\) is used to indicate a submatrix (in subscript), its dependency on \(k\) and \(x\) will be omitted for the sake of brevity of notation.
			\(H_k(x)\) can now be written as
			\begin{equation}\label{eq:H}
				H_k(x)
			{}={}
				Q+\trans A\submatrix{\J\cdot}(\Sy_k)\submatrix{\J\J}A\submatrix{\J\cdot}
				{}+{}
				\Sx^{-1},
			\end{equation}
			
			The semismooth Newton direction \(d\) at \(x\) satisfies 
			\begin{equation}\label{eq:Newton_SCHUR}
			H_k(x)d=-\nabla\varphi_k(x).
			\end{equation}
			 Denoting \(\lambda\coloneqq(\Sy_k)\submatrix{\J\J}A\submatrix{\J\cdot} d\), the computation of \(d\) is equivalent to solving the following extended linear system
			\begin{equation}\label{eq:Newton_KKT}
				\K_k(x) \hspace{0.1cm} \begin{bmatrix}
					d\\
					\lambda
				\end{bmatrix}
			{}={}
				\begin{bmatrix}
					Q+\Sx^{-1} & \trans A\submatrix{\J\cdot} \\
					A\submatrix{\J\cdot} & -(\Sy_k)\submatrix{\J\J}^{-1}
				\end{bmatrix}
				\begin{bmatrix}
					d\\
					\lambda
				\end{bmatrix}
			{}={}
				\begin{bmatrix}
					-\nabla\varphi_k(x)\\
					0
				\end{bmatrix}.
			\end{equation}
			Finding the solution of either linear system \eqref{eq:Newton_SCHUR} or \eqref{eq:Newton_KKT} in an efficient manner is discussed in further detail in \cref{subsec:Solving linear systems}.

		\subsection{Exact linesearch}

			Once a suitable direction \(d\) has been found, a stepsize \(\tau\) needs to be computed.
			This is typically done via a linesearch on a suitable merit function.
			QPALM can compute the optimal stepsize by using the piecewise quadratic function \(\psi(\tau) = \varphi_k(x+\tau d)\) as the merit function.
			Finding the optimal stepsize is therefore equivalent to finding a zero of 
			
			{\ifarxiv\else\mathtight[0.85]\fi
				\begin{align*}
					\psi'(\tau)
				{}={} &
					\innprod{\nabla\varphi_k(x+\tau d)}{d}
				\\
				{}={} &
					\innprod{d}{\nabla f(x+\tau d)+\Sx^{-1}(x+\tau d -\hat x^k)}
					{}+{}
					\innprod{Ad}{
						 y^{k}
						{}+{}
						\Sy_k\bigl(
							A(x+\tau d)-Z_k(x+\tau d)
						\bigr)
					}
				\\
				{}={} &
					\tau\innprod{d}{(Q+\Sx^{-1})d}
					{}+{}
					\innprod{d}{Qx+\Sx^{-1}(x-\hat x^k)+q} 
				\\
				&
					{}+{} 
					\innprod{\Sy_k Ad}{
						\bigl[
							Ax+\Sy_k^{-1} y^{k}-u
							{}+{}
							\tau Ad
						\bigr]_+
					}
					{}-{}
						\innprod{\Sy_k Ad}{
							\bigl[
								\ell-Ax-\Sy_k^{-1} y^{k}
								{}-{}
								\tau Ad
							\bigr]_+
						}
				\\\numberthis\label{eq:psi'}
				{}={} &
					\eta\tau+\beta+\innprod{\delta}{[\delta\tau-\alpha]_+},
				\end{align*}
			}%
				where
				\begin{equation} \label{eq:psi' symbols}
					\begin{cases}[r @{{}\ni{}} rl]
						\R
						&
						\eta
					{}\coloneqq{}&
						\innprod{d}{(Q+\Sx^{-1})d},
					\\
						\R
						&
						\beta
					{}\coloneqq{}&
						\innprod{d}{Qx+\Sx^{-1}(x-\hat x^k)+q},
					\\
						\R^{2m}
						&
						\delta
					{}\coloneqq{}&
						\bigl[-\Sy_k^{\nicefrac12}Ad~~~\Sy_k^{\nicefrac12}Ad\bigr],
					\\
						\R^{2m}
						&
						\alpha
					{}\coloneqq{}&
						\Sy_k^{-\nicefrac12}\bigl[ y^{k}+\Sy_k(Ax-\ell)~~~\Sy_k(u-Ax)\mathrlap{{}- y^{k}\bigr].}
					\end{cases}
				\end{equation}
			
			Note that \(\psi'\) is a monotonically increasing piecewise affine function.
			The zero of this function can be found by sorting all the breakpoints, and starting from \(0\) going through these points \(t_i\) until \(\psi'(t_i) > 0\).
			The optimal stepsize is then in between this and the previous breakpoint and can easily be retrieved by means of interpolation.
			This procedure is outlined in \cref{alg:LS}.
			 
			\begin{algorithm}[t]
				\algcaption{Exact linesearch}
				\label{alg:LS}

		\begin{algorithmic}[1]
		\Require
			\(x,d\in\R^n\),~
			diagonal \(\Sigma\in\symm_{++}(\R^n)\)
		\Provide
			optimal stepsize \(\tau_\star\in\R\)
		\State
			Let \(\func{\psi'}{\R}{\R}\), \(\alpha,\beta\in\R\) and \(\delta,\eta\in\R^{2m}\) be as in \eqref{eq:psi'} and \eqref{eq:psi' symbols}
		\State
			\begin{tabular}[t]{@{}l@{}}
				Define the set of breakpoints of \(\psi'\)
			\\
				\(
					T
				{}={}
					\set{
						\frac{\alpha_i}{\delta_i}
					}[
						i=1,\ldots,2m,~
						\delta_i\neq 0
					]
				\)
			\end{tabular}
		\State
			Sort \(T=\set{t_1,t_2,\ldots}\) such that \(t_i<t_{i+1}\) for all \(i\)
		\State
			Let \(t_i\in T\) be the smallest such that \(\psi'(t_i)\geq0\)
		\State
			\Return
				\(
					\tau_\star
				{}={}
					t_{i-1}
					{}-{}
					\frac{t_i-t_{i-1}}{\psi'(t_i)-\psi'(t_{i-1})}\psi'(t_{i-1})
				\)
				\quad{\footnotesize
					(if \(i=1\) then \(t_{i-1} = 0\))%
				}%
		\end{algorithmic}
			\end{algorithm}

	\section{Linear algebra code}\label{sec:Code}

		\review{The QPALM algorithm is implemented in standalone open-source C code,\footnote{%
			\url{https://github.com/Benny44/QPALM_vLADEL}%
		}
		licensed under the GNU Lesser General Public License version 3 (LGPL v3).
		QPALM also provides interfaces to MATLAB, Python and Julia.}
		
		This section further discusses the relevant linear algebra used in QPALM, which is implemented in a standalone C package LADEL,\footnote{\url{https://github.com/Benny44/LADEL}} and the routine that is used to compute the minimum eigenvalue of a symmetric matrix\review{, which, as mentioned in \cref{sec:PALM}, is used to guarantee convexity of the subproblems.}

		\subsection{Solving linear systems}\label{subsec:Solving linear systems}

			The most computationally expensive step in one iteration of QPALM is solving the semismooth Newton system \eqref{eq:Newton_SCHUR} or \eqref{eq:Newton_KKT}.
			The matrix \(\K_k(x)\) in \eqref{eq:Newton_KKT}, without penalty parameters, is readily recognized as the system of equations that represent the first-order necessary conditions of equality-constrained QPs \cite[\S16.1]{nocedal2006numerical}, and is therefore dubbed the Karush-Kuhn-Tucker (KKT) matrix.
			The matrix \(H_k(x)\) is the Schur complement of \(\K_k(x)\) with respect to the \(-(\Sy_k)\submatrix{\J, \J}^{-1}\) block, and is therefore dubbed the Schur matrix.
			Solving either of the two systems results in a direction along which we can update the primal iterate \(x\).
			The reader is referred to \cite{benzi2005numerical} for a broad overview of solution methods for such systems, including direct and iterative methods.
			In the case of QPALM, the matrix \(\K_k(x)\) or \(H_k(x)\) is decomposed as the product of a unit lower diagonal matrix \(L\), a diagonal matrix \(D\) and the transpose of \(L\).
			This is more commonly known as an \(\LDL\) factorization.
			The factorization and updates are slightly different for \(\K_k(x)\) and \(H_k(x)\), so these cases will be discussed separately.

			\subsubsection{KKT system}
				It is not guaranteed that an \(\LDL\) factorization, with \(D\) diagonal, can be found for every matrix.
				However, because \(\Sy_k\) and \(Q+\Sx^{-1}\) are symmetric positive definite by construction, \(\K_k(x)\) can readily be recognized as a symmetric quasidefinite matrix, which is strongly factorizable \cite[Theorem 2.1]{vanderbei1995symmetric}.
				A symmetric matrix \(K\) is strongly factorizable if for any symmetric \review{permutation \(P\) there exists} a unit lower diagonal matrix \(L\) and a diagonal matrix \(D\) such that \(PK\trans P = \LDL\).
				In other words, we should always be able to find an \(\LDL\) factorization of \(\K_k(x)\) with \(D\) diagonal.
				To find such a factorization, LADEL has implemented a simple uplooking Cholesky method with separation of the diagonal elements, see \cite{davis2005algorithm}.
			
				A crucial step in maintaining sparsity during the factorization is to find an effective permutation.
				Moreover, permutations are sometimes used to prevent \review{ill con\-di\-tion\-ing}.
				However, finding the optimal permutation is an NP-hard problem \cite{yannakakis1981computing}.
				Various heuristics have been developed, an overview of which can be found in \cite[\S7.7]{davis2006direct}.
				LADEL uses the open-source (BSD-3 licensed) implementation\footnote{\url{https://github.com/DrTimothyAldenDavis/SuiteSparse/tree/master/AMD}} of the approximate minimum degree (AMD) ordering algorithm discussed in \cite{amestoy2004algorithm}.
			
				In QPALM, a fill-reducing ordering of the full KKT system, i.e. with \(\J = [1,m]\), is computed using AMD once before the first factorization and is used during the remainder of the solve routine.
				Hence, this permutation minimizes the fill-in of the worst case, that is with all constraints active.
				In fact, when solving the KKT system, we will not consider \(\K_k(x)\) directly, but rather an augmented version \[\widetilde \K_k(x) = \begin{bmatrix}
						Q+\Sx^{-1} &  \trans{(A\supmatrix{\J\cdot})} \\
						A\supmatrix{\J\cdot} & -(\Sy_k)^{-1}
					\end{bmatrix}.\]
				Note that, as mentioned in \cref{subsec:Notation}, \(A\supmatrix{\J}\) is the \(m \times n\) matrix, with \(A\supmatrix{\J}_{j\cdot} = A_{j\cdot}\) if \(j\in \J\) and zero otherwise.
				The size of \(\widetilde \K_k(x)\) is therefore always \((m+n)\times (m+n)\), but due to \eqref{eq:supmatrix} all the inactive constraints give rise to rows and columns that are \(\mathbf{0}\) apart from the diagonal element.
				Combined with \eqref{eq:Newton_KKT}, it immediately follows that \(\lambda_j = 0\) for \(j \not \in \J_k(x)\).

				Before the condition of \cref{state:PALM:x} is satisfied, several Newton steps may be required.
				However, during these iterations \(k\) remains constant, and so does \(\Sy_k\).
				Therefore, the only manner in which \(\widetilde \K_k(x)\) changes is as a result of the change in active constraints when \(x\) is updated.
				Instead of refactorizing the matrix \(\widetilde \K_k(x)\), we can instead use sparse factorization update routines to update the existing factorization matrices \(L\) and \(D\).
				In particular, LADEL has implemented the row addition and row deletion algorithms of \cite{davis2005row}, with minor modifications to allow for negative diagonal elements (indefinite systems), as outlined in \cref{alg:rowadd} and \cref{alg:rowdel}.
			
				\begin{algorithm}[t]
					\algcaption{Row addition (see \cite{davis2005row}, with modifications in \cref{state:rowadd:w} and \cref{state:rowadd:updown})}
					\label{alg:rowadd}

			\begin{algorithmic}[1]
			\Require
				\begin{tabular}[t]{@{}l@{}}
					\(\LDL\) factors \(L\) and \(D\) of a matrix \(C \in \R^{n\times n}\) with \(C_{\beta\cdot} = C_{\cdot \beta} = \mathbf{0}\) except for \(C_{\beta\beta} = \epsilon\)
				\\
					Let \(\alpha = 1:\beta-1\) and \(\gamma = \beta+1:n\), then
				\\
				\(\displaystyle
					\LDL
				{}={}
					\begin{bmatrix}
					L_{\alpha \alpha} & & \\
					0 & 1 & \\
					L_{\gamma \alpha} & 0 & L_{\gamma\gamma} 
					\end{bmatrix}
					\begin{bmatrix}
					D_{\alpha\alpha} & & \\
					& d_{\beta\beta} & \\
					& & D_{\gamma\gamma}
					\end{bmatrix}
					\begin{bmatrix}
					\trans L_{\alpha \alpha} & & \trans L_{\gamma \alpha}\\
					& 1 & 0\\
					& & \trans L_{\gamma\gamma} 
					\end{bmatrix}
				{}={}
					\begin{bmatrix}
					C_{\alpha\alpha} & 0 & \trans C_{\gamma\alpha}\\ 
					0 & \epsilon & 0 \\
					C_{\gamma\alpha} & 0 & C_{\gamma\gamma}
					\end{bmatrix}
				\)
				\end{tabular}
			\Provide
				\begin{tabular}[t]{@{}l@{}}
					Updated \(\LDL\) factors \(\bar L\) and \(\bar D\) of \(\bar C\) which is equal to \(C\) except with the \(\beta\)-th row and column
				\\
					replaced by \(\trans{\bar c}_{\cdot \beta}\) and \(\bar c_{\cdot \beta}\) respectively, i.e.
				\\
				\(\displaystyle
					\bar L \bar D \trans{\bar L}
				{}={}
					\begin{bmatrix}
					L_{\alpha \alpha} & & \\
					\trans{\bar l}_{\alpha \beta} & 1 & \\
					L_{\gamma \alpha} & \bar l_{\gamma \beta} & \bar L_{\gamma\gamma} 
					\end{bmatrix}
					\begin{bmatrix}
					D_{\alpha\alpha} & & \\
					& \bar d_{\beta\beta} & \\
					& & \bar D_{\gamma \gamma}
					\end{bmatrix}
					\begin{bmatrix}
					\trans L_{\alpha \alpha} & \bar l_{\alpha \beta} & \trans L_{\gamma \alpha}\\
					& 1 & \trans{\bar l}_{\gamma \beta}\\
					& & \trans{\bar L}_{\gamma \gamma} 
					\end{bmatrix}
				{}={}
					\begin{bmatrix}
					C_{\alpha\alpha} & \bar c_{\alpha \beta} & \trans C_{\gamma\alpha} \\ 
					\trans{\bar c}_{\alpha \beta} & \bar c_{\beta \beta} & \trans{\bar c}_{\gamma \beta} \\
					C_{\gamma\alpha} & \bar{c}_{\gamma \beta} & C_{\gamma\gamma}
					\end{bmatrix}
				\)
				\end{tabular}
			\State
				Solve the lower triangular system \(L_{\alpha \alpha} D_{\alpha \alpha} \bar l_{\alpha \beta} = \bar c_{\alpha \beta}\) to find \(\bar l_{\alpha \beta}\)
			\State
				\(\bar d_{\beta \beta} = \bar c_{\beta \beta} - \trans{\bar l}_{\alpha \beta} D_{\alpha \alpha} \bar l_{\alpha \beta}\)
			\State
				\(\bar l_{\gamma \beta} = {\bar d}^{-1}_{\beta \beta} (\bar c_{\gamma \beta} - L_{\gamma \alpha} D_{\alpha \alpha} \bar l_{\alpha \beta})\)
			\State\label{state:rowadd:w}%
				\(w = \bar l_{\gamma \beta} \sqrt{|\bar d_{\beta \beta}|}\)
			\State\label{state:rowadd:updown}%
				Perform the rank-1 update or downdate \(\bar L_{\gamma \gamma} \bar D_{\gamma \gamma} \trans{\bar L}_{\gamma \gamma} = L_{\gamma \gamma}D_{\gamma \gamma} \trans L_{\gamma \gamma} - \sign(d_{\beta \beta}) w \trans w \)
			\end{algorithmic}
				\end{algorithm}
			
				\begin{algorithm}[t]
					\algcaption{Row deletion (see \cite{davis2005row}, with modifications in \cref{state:rowdel:w} and \cref{state:rowdel:updown})}
					\label{alg:rowdel}

			\begin{algorithmic}[1]
			\Require
				\begin{tabular}[t]{@{}l@{}}
					\(\LDL\) factors \(L\) and \(D\) of a matrix \(C \in \R^{n\times n}\).
					Let \(\alpha = 1:\beta-1\) and \(\gamma = \beta+1:n\), then 
				\\
					\(\displaystyle
						\LDL
					{}={}
						\begin{bmatrix}
						L_{\alpha \alpha} & & \\
						\trans{l}_{\alpha \beta} & 1 & \\
						L_{\gamma \alpha} & l_{\gamma \beta} & L_{\gamma\gamma} 
						\end{bmatrix}
						\begin{bmatrix}
						D_{\alpha\alpha} & & \\
						& d_{\beta\beta} & \\
						& & D_{\gamma \gamma}
						\end{bmatrix}
						\begin{bmatrix}
						\trans L_{\alpha \alpha} & l_{\alpha \beta} & \trans L_{\gamma \alpha}\\
						& 1 & \trans{l}_{\gamma \beta}\\
						& & \trans{L}_{\gamma \gamma} 
						\end{bmatrix}
					{}={}
						\begin{bmatrix}
						C_{\alpha\alpha} & c_{\alpha \beta} & \trans C_{\gamma\alpha} \\ 
						\trans{c}_{\alpha \beta} & c_{\beta \beta} & \trans{c}_{\gamma \beta} \\
						C_{\gamma\alpha} & c_{\gamma \beta} & C_{\gamma\gamma}
						\end{bmatrix}
					\)
				\end{tabular}
			\Provide
				\begin{tabular}[t]{@{}l@{}}
					Updated \(\LDL\) factors \(\bar L\) and \(\bar D\) of \(\bar C\) which is equal to \(C\) except with the \(\beta\)-th row and column
				\\
					deleted and the diagonal element \(c_{\beta \beta}\) replaced by \(\epsilon\), i.e.
				\\
					\(
						\bar L \bar D \trans{\bar L}
					{}={}
						\begin{bmatrix}
						L_{\alpha \alpha} & & \\
						0 & 1 & \\
						L_{\gamma \alpha} & 0 & \bar L_{\gamma\gamma} 
						\end{bmatrix}
						\begin{bmatrix}
						D_{\alpha\alpha} & & \\
						& \bar d_{\beta\beta} & \\
						& & \bar D_{\gamma\gamma}
						\end{bmatrix}
						\begin{bmatrix}
						\trans L_{\alpha \alpha} & & \trans L_{\gamma \alpha}\\
						& 1 & 0\\
						& & \trans{\bar L}_{\gamma\gamma} 
						\end{bmatrix}
					{}={}
						\begin{bmatrix}
						C_{\alpha\alpha} & 0 & \trans C_{\gamma\alpha} \\ 
						0 & \epsilon & 0 \\
						C_{\gamma\alpha} & 0 & C_{\gamma\gamma}
						\end{bmatrix}
					\)
				\end{tabular}
			\State
				\(\bar l_{\alpha \beta} = 0 \)
			\State
				\(\bar d_{\beta \beta} = \epsilon \)
			\State
				\(\bar l_{\gamma \beta} = 0 \)
			\State\label{state:rowdel:w}%
				\(w = l_{\gamma \beta} \sqrt{|d_{\beta \beta}|}\)
			\State\label{state:rowdel:updown}%
				Perform the rank-1 update or downdate \(\bar L_{\gamma \gamma} \bar D_{\gamma \gamma} \trans{\bar L}_{\gamma \gamma} = L_{\gamma \gamma}D_{\gamma \gamma} \trans L_{\gamma \gamma} + \sign(d_{\beta \beta}) w \trans w \)
			\end{algorithmic}
				\end{algorithm}

			\subsubsection{Schur system}
				The Schur matrix \(H_k(x)\) is symmetric positive definite, since it is the sum of a positive definite matrix (\(Q+\Sx^{-1}\)) and a positive semidefinite matrix.
				Therefore, a Cholesky factorization of \(H_k(x)\) exists.
				Furhtermore, when \(k\) remains constant and \(x\) changes to \(x^+\), the difference between \(H_k(x^+)\) and \(H_k(x)\) is given by
				\[
				H_k(x^+) - H_k(x) = \trans A\submatrix{\J^e\cdot}(\Sy_k)\submatrix{\J^e \J^e}A\submatrix{\J^e\cdot} - \trans A\submatrix{\J^l\cdot}(\Sy_k)\submatrix{\J^l\J^l}A\submatrix{\J^l\cdot},
				\]
				with \(\J^e = \J_k(x^+) \setminus \J_k(x)\) and \(\J^l = \J_k(x) \setminus \J_k(x^+)\) the sets of constraints respectively entering and leaving the active set.
				Therefore, two \review{low-rank} Cholesky factorization updates can be performed \cite{davis1999modifying,davis2001multiple}.
				LADEL has implemented the one-rank update routines in \cite{davis1999modifying}, which are slightly less efficient than the multiple-rank routines outlined in \cite{davis2001multiple}, as implemented in CHOLMOD \cite{chen2008algorithm}.

			\subsubsection{Choosing a system}
				Let \(H\) and \(\K\) denote the ``full'' Schur and KKT matrices, that is with \(\J = [1,m]\).
				In QPALM, we automatically choose which of these systems to factorize, depending on an estimate of the floating point operations required for each.
				The work required to compute an \(\LDL\) factorization is \(\sum |L_{\cdot i}|^2\).
				However, we do not have access to the column counts of the factors before the symbolic factorization.
				Therefore, we try to compute a rough estimate of the column counts of the factor via the column counts of the matrices themselves.
				Moreover, we consider an average column count for each matrix rather than counting the nonzeros in each individual column.
				As such, QPALM uses the following quantity to determine the choice of linear system:
				\begin{align*}
				\frac{\displaystyle\sum_{i = 1}^{n+m} |L^{ \K}_{\cdot i}|^2}{\displaystyle\sum_{i = 1}^{n} |L^{H}_{\cdot i}|^2} 
				{}\approx{}
				\frac{\displaystyle\sum_{i = 1}^{n+m} | \K_{\cdot i}|^2}{\displaystyle\sum_{i = 1}^{n} |H_{\cdot i}|^2} 
				{}\approx{}
				\frac{\displaystyle\sum_{i = 1}^{n+m} \left(\frac{| \K|}{n+m}\right)^2}{\displaystyle\sum_{i = 1}^{n} \left(\frac{|H|}{n}\right)^2} 
				{}={}
				\frac{n}{n+m}\frac{| \K|^2}{|H|^2}
				{}\approx{}
				\frac{n}{n+m}\frac{| \K|^2}{|\widetilde H|^2} 
				,
				\end{align*} 
				with \(L^{\widetilde \K}\) and \(L^{H}\) the lower diagonal factors of the corresponding matrix.
				Computing \(|H|\) exactly requires the same order of work as computing \(H\) itself.
				Depending on the sparsity pattern of \(Q\) and \(A\), \(H\) can be much denser than \(\K\).
				Hence, we do not want to compute \(H\) before choosing between the two systems.
				Instead, we further (over)estimate \(|H|\) by \(|\widetilde H|\) considering separate contributions from \(Q + \Sx^{-1}\) and from \(\trans A \Sy_k A\).
				Note that a row in \(A\) with \(|A_{i\cdot}|\) nonzero elements contributes a block in \(\trans A A\) with \(|A_{i\cdot}|^2\) elements.
				By discounting the diagonal elements, which are present in \(\Sx^{-1}\), this becomes \(|A_{i\cdot}|^2-|A_{i\cdot}|\).
				The overlap of different elements cannot be accounted for (cheaply).
				Therefore, in our estimate we deduct the minimum (possible) amount of overlap of each block with the biggest block.
				Denoting \(\hat A = |A_{\hat i\cdot}| = \max_i(|A_{i\cdot}|)\), this overlap, again discounting diagonal elements, is given as \([\hat A + |A_{i\cdot}| - n]_+^2-[\hat A + |A_{i\cdot}| - n]_+\), and so our estimate for \(|H|\) is 
				\[
				|\widetilde H| {}={} |Q + \Sx^{-1}| {}+{} \hat A^2-\hat A {}+{} \displaystyle\sum_{i\neq \hat i} |A_{i\cdot}|^2-|A_{i\cdot}| - [\hat A + |A_{i\cdot}| - n]_+^2 + [\hat A + |A_{i\cdot}| - n]_+.
				\]
			
\begin{figure}
					\centering
					\includetikz[width=\linewidth]{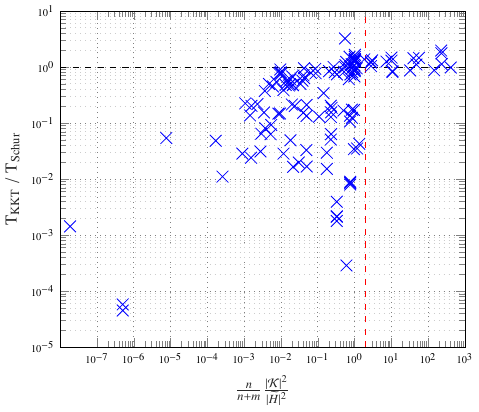}
					\caption{Runtime comparison of KKT and Schur complement methods when applying QPALM to the Maros Meszaros test set.}
					\label{fig:kkt_vs_schur}
				\end{figure}

				Finally, \cref{fig:kkt_vs_schur} compares the runtimes of QPALM solving either the KKT or the Schur system applied to the Maros Meszaros test set.
				Note that the runtime of QPALM using the KKT system can still be much lower than that using the Schur system for an estimated nonzero ratio of 1.
				This is why a heuristic threshold value of 2, indicated by the red dotted line, is chosen for this ratio such that the default option of QPALM switches automatically between the two systems.
				The user also has the option to specify beforehand which system to solve.
				
		\subsection{Computing the minimum eigenvalue}		

			
			Finding the solution of a large symmetric eigenvalue problem has been the topic of a substantial body of research, and many methods exist.
			They are typically divided into two categories: direct methods, which find all eigenvalues, and iterative methods, which find some (or all) eigenvalues.
			The reader is referred to \cite[\S8]{golub2013matrix} for a detailed overview of (the origin of) these methods.
			In our case, since we only need the minimum eigenvalue of \(Q\), iterative methods seem more promising.
			Of these, the Locally Optimal Block Preconditioned Conjugate Gradient (LOBPCG) method, developed by Knyazev \cite{knyazev2001toward}, demonstrated the best performance regarding robustness and speed of convergence in our tests.
			A dedicated implementation of LO(B)PCG to find only the minimum eigenvalue and its corresponding eigenvector was added in QPALM.
			This method iteratively minimizes the Rayleigh quotient \(\frac{\trans x Qx}{\trans x x}\) in a Krylov subspace spanned by three vectors: the current eigenvector estimate \(x^k\), the current residual \(w^k = Qx^k - \lambda^k x^k\), and a conjugate gradient direction \(p^k\).
			The details of the implementation in QPALM can be found in \cref{alg:lobpcg}.
			The computational cost of this algorithm per iteration is essentially a matrix-vector product and solving a 3-by-3 generalized eigenvalue system.
			\review{The latter is performed in our code by finding the roots of a one-dimensional third-order polynomial and performing a Gaussian elimination of a 3-by-3 system.}
			Note that \cref{alg:lobpcg} is very similar to \cite[Algorithm 4.1]{knyazev2001toward}, aside from some scaling.
			
			\begin{algorithm}[t]
				\algcaption{LO(B)PCG}
				\label{alg:lobpcg}

		\begin{algorithmic}[1]
		\Require
			\(x^0 \in \R^n, \varepsilon > 0\) and \(Q\in \symm(\R^n)\)
		\Provide
			Lower bound on \(\lambda^*=\lambda_{\rm min}(Q)\) and estimate of the corresponding eigenvector \(x^*\) of \(Q\)
		\State%
			Initialize \(\lambda^0 = \innprod{x^0}{Ax^0}\), and \(w^0 = Qx^0 - \lambda^0 x^0\), and let \(S = [x^0, w^0]\)
		\State%
			Solve \(\trans S QS y = \mu Sy\), set \(\lambda^1 = \min(\mu)\) and let \(\tilde y\) denote the corresponding eigenvector
		\State%
			\(x^1 = \tilde y_1 x^0 + \tilde y_2 w^0\)
		\State%
			\(p^1 = \tilde y_2 w^0\)
		\For{ \(k=1,2,\ldots\) }
			\State%
				\(w^k = Qx^k - \lambda^k x^k\)
			\If{ \(\|w^k\|_2\leq\varepsilon\) }
				\State
					Return \(\lambda^* = \lambda^k - \|w^k\|_2\) and \(x^* = x^k\)
			\EndIf{}
			\State
				Let \(S = [x^k, w^k, p^k]\)
			\State
				Solve \(\trans S QS y = \mu Sy\), set \(\lambda^k = \min(\mu)\) and let \(\tilde y\) denote the corresponding eigenvector
			\State
				\(x^{k+1} = \tilde y_1 x^k + \tilde y_2 w^k + \tilde y_3 p^k\)
			\State
				\(p^{k+1} = \tilde y_2 w^k + \tilde y_3 p^k\)
		\EndFor{}
		\end{algorithmic}
			\end{algorithm}

	\section{Parameter selection}\label{sec:Parameter}

			The little details can make or break the practical performance of an algorithm.
			This section discusses aspects that make QPALM more robust, such as preconditioning of the problem data, and the most important parameter settings.
			Some of these parameters and parameter update criteria have been tuned manually and some are computed automatically based on the problem data or current iterates.
			The last subsection also lays out in detail the termination criteria employed by QPALM.

		\subsection{Factorization updates}

			As mentioned in \cref{subsec:Solving linear systems}, in between Newton iterations we can update the factorization instead of having to refactorize from scratch.
			However, in practice, the factorization update routines will only be more efficient if the number of constraints entering and leaving the active set is relatively low.
			Hence, when the active set changes significantly, we want to recompute the factorization instead.
			After some experimental tuning, we decided on the following criterion to do an update:
			\[
			|\J^e| + |\J^l| \leq \min(\texttt{max\_rank\_update}, \texttt{ max\_rank\_update\_fraction} \cdot (n+m)),
			\] 
			with \(\texttt{max\_rank\_update} = 160\) an absolute limit and \texttt{max\_rank\_update\_fraction} \(= 0.1\) a relative limit on the number of changing active constraints.
			Both of these parameters can also be set by the user.

		\subsection{Preconditioning}

			Most optimization solvers will perform a scaling of the problem in order to prevent too much \review{ill conditioning}.
			A standard idea in nonlinear optimization is to evaluate the objective and constraints and/or the norm of their gradients at a representative point and scale them respectively with their inverse, see for example \cite[\S12.5]{birgin2014practical}.
			However, the quality of this scaling depends of course on the degree to which the initial point \review{is representative} for the iterates, and by extension the solution.
			Furthermore, since we are dealing with a QP, the constraints and objective are all determined by matrices.
			Therefore, it makes sense to equilibrate these matrices directly, as is done in OSQP for example \cite[\S5.1]{stellato2020osqp}.
			OSQP applies a modified Ruiz equilibration \cite{ruiz2001scaling} to the KKT matrix.
			This equilibration routine iteratively scales the rows and columns in order to make their infinity norms go to 1.
			OSQP adds an additional step that scales the objective to take into account also the linear part \(q\).
			We have found in our benchmarks, however, that instead of this scaling it is better to apply Ruiz equilibration to the constraints only, and to scale the objective by a single constant.
			Why exactly this is \review{a better strategy} is unknown to us, but we suspect that the constraints are more sensitive to the scaling, so it might be better to deal with them separately.
			
			\begin{algorithm}[t]
				\algcaption{Ruiz equilibration \cite{ruiz2001scaling}}
				\label{alg:ruiz}

		\begin{algorithmic}[1]
		\Require
			\(A \in \R^{m\times n}\)
		\Provide
			\(D\in\R^n, E\in\R^m\) and \(\bar A = EAD\)
		\State%
			Initialize \(\bar A = A\),~ \(D = \Id_n\),~ \(E = \Id_m\)
		\For{ \(k=1,\ldots,\texttt{scaling}\) }
			\For{ \(i=1,\ldots,m\) }
				 \State%
				 	\(\bar E_{ii} = \sqrt{\|\hat A_{i\cdot}\|_\infty}\)
			\EndFor{}
			\For{ \(j=1,\ldots,n\) }
				 \State%
				 	\(\bar D_{jj} = \sqrt{\|\bar A_{\cdot j}\|_\infty}\)
			\EndFor{}
			\State%
				\(\bar A = {\bar E}^{-1} \bar A {\bar D}^{-1}\)
			\State%
				\(D = D {\bar D}^{-1} \)
			\State%
				\(E = E {\bar E}^{-1} \)
		\EndFor{}
		\end{algorithmic}
			\end{algorithm}
			
			In QPALM we apply Ruiz equilibration, outlined in \cref{alg:ruiz}, to the constraint matrix \(A\), yielding \(\bar A = EAD\).
			The setting \texttt{scaling} denotes the number of scaling iterations which can be set by the user and defaults to 10.
			The objective is furthermore scaled by \(c = \max(1.0, \|D(Q x^0 + q)\|_\infty)^{-1}\).
			In conclusion, we arrive at a scaled version of \eqref{eq:P}
			\[
				\minimize_{\bar x\in\R^n}\tfrac12\trans{\bar x} \bar Q \bar x+\trans{\bar q} \bar x
			\quad\stt{}
				\bar A \bar x\in \bar C,
			\]
			with \(\bar x = D^{-1}x\), \(\bar Q = cDQD\), \(\bar q = cDq\), \(\bar A = EAD\), \(\bar C=\set{z\in\R^m}[\bar \ell\leq z \leq \bar u]\), \(\bar \ell = E\ell\) and \(\bar u = Eu\).
			The dual variables in this problem are \(\bar y = cE^{-1}y\).
			This is the problem that is actually solved in QPALM, although the termination criteria are unscaled, that is they apply to the original problem \eqref{eq:P}, see \cref{subsec:termination}.

		\subsection{Penalty parameters}

			The choice of the penalty parameters, and the rules used to update them have been found to be a decisive factor in the performance of QPALM.
			In this section we discuss both the traditional penalty parameters arising in the augmented Lagrangian formulation \(\Sy\), and the proximal penalty parameters \(\Sx\).
			
			\subsubsection{Dual penalty parameters}
				The dual penalty parameters \(\Sy\) play an integral role in slowly but surely enforcing feasibility over the iterations.
				Because the inner subproblems solved by QPALM are strongly convex, there is no theoretical requirement on the penalty parameters, other than the obvious one of them being positive.
				However, experience with augmented Lagrangian methods suggests that high values can inhibit convergence initially, as they then introduce \review{ill conditioning} in the problem, whereas high values near a solution are very useful to enforce feasibility.
				As such, these penalty parameters are typically increased during the solve routine depending on the constraint violations of the current iterate.
				A standard rule is to (only) increase the parameters when the respective constraint violations have not decreased sufficiently, see \cite[\S17.4]{nocedal2006numerical}.
				Furthermore, an added rule in QPALM that is observed to work well is to increase the penalties proportional to their corresponding constraint violation.
				Hence, we employ the following strategy to find \(\Sy_{k+1}\) in \cref{state:y,state:PALM:inner} \review{(based on scaled quantities):}
				\begin{equation} \label{eq:update sigma} \mathtight
				\frac{(\Sy_{k+1})_{ii}}{(\Sy_{k})_{ii}} = \begin{ifcases}
						1.0 & |(\bar A\bar x^{k+1}-\bar z^{k+1})_i| < \theta |(\bar A\bar x^{k}-\bar z^{k})_i| \\
						\min\left[\frac{\sigma_{\max}}{(\Sy_{k})_{ii}}, \max\left(\mathbf{\Delta}\frac{|(\bar A\bar x^{k+1}-\bar z^{k+1})_i|}{\|\bar A\bar x^{k+1}-\bar z^{k+1}\|_\infty}, 1.0\right)\right] \otherwise.
					\end{ifcases}
				\end{equation}
				The default parameters here are \(\theta = 0.25\), \(\mathbf{\Delta} = 100\), and \(\sigma_{\max} = 10^9\) and can all be set by the user.
				The usage of this rule, in particular letting the factor depend on the constraint violation itself, has been a crucial step in making the performance of QPALM more robust.
				\review{In fact, it seems to us to be the key difference when compared to the corresponding penalty update in OSQP \cite[\S5.2]{stellato2020osqp}, which uses the same factor for all the constraints and updates it only seldom based on the setup time and current runtime, in contrast to QPALM where each penalty is updated whenever the corresponding constraint violation has not sufficiently decreased (and is not already relatively small).}
				
				Note that in case only a few penalties are modified, the factorization of either \(\widetilde \K\) or \(H\) may be updated using low-rank update routines.
				In practice, we set the limit on the amount of changing penalties a bit lower as we expect an additional update to be required for the change in active constraints.
			
				As with regards to an initial choice of penalty parameters, the formula proposed in \cite[\S12.4]{birgin2014practical} was found to be effective after some tweaking of the parameters inside.
				As such we use the following rule to determine initial values of the penalties
				\begin{equation} \label{eq:sigma init}
				(\Sy_0)_{ii} = \max\left[10^{-4}, \min\left(\sigma_{\textrm{init}} \frac{\max (1.0, |\tfrac12 \trans{({\bar x}^0)} \bar Q \bar x^0 + \trans{\bar q} \bar x^0|)}{\max (1.0, \tfrac12\|\bar A {\bar x}^0-\proj_{\bar C}(\bar A {\bar x}^0)\|^2)}, 10^{4}\right) \right], 
				\end{equation}
				with \(\sigma_{\textrm{init}}\) a parameter with a default value of 20 and which can also be set by the user.
				Setting the initial penalty parameters to a high value can be very beneficial when provided with a good (feasible) initial guess, as therefore feasibility will not be lost.
				An investigation into this and warm-starting QPALM in general is a topic for future work.

			\subsubsection{Primal penalty parameters}
				The primal, or proximal, penalty parameters \(\Sx\) serve to regularize the QP around the ``current'' point \({\hat x}^k\).
				An appropriate choice makes it so \review{that} the subproblems are strongly convex, as discussed before.
				In many problems, the user knows whether \review{the QP at hand is convex or not}.
				Therefore, QPALM allows the user to indicate which case is dealt with.
				If the user indicates the problem is (or might be) nonconvex, i.e. that \(Q\) is not necessarily positive semidefinite, QPALM uses \cref{alg:lobpcg} to obtain a tight lower bound \(\lambda^*\) on the minimum eigenvalue.
				If this value is negative, we set \(\forall i: \Sx_{ii} = \frac{1}{|\lambda^*-10^{-6}|}\).
				Otherwise, or \review{in case} the user indicates the problem is convex, the default value \review{is \(\Sx_{ii}^{-1} = 10^{-7}\)}, a reasonably low value to not interfere with the convergence speed while guaranteeing that \(H_k(x)\) or \(\widetilde \K_k(x)\) is positive definite or quasidefinite respectively.
				\review{This penalty has then a very similar effect as the Hessian regularization in OSQP \cite{stellato2020osqp}, where the default is \(10^{-6}\Id\).}
				Furthermore, in the convex case, if the convergence is slow but the primal termination criterion \eqref{eq:primal termination} is already satisfied, \review{\(\Sx_{ii}^{-1}\) may be further decreased to \(10^{-12}\)} depending on an estimate of the (machine accuracy) errors that would be accumulated in \(H_k(x)\).
				Finally, QPALM also allows the selection of an initial \(\Sx_{0,ii} = \gamma_{\rm init}\), and an update rule \(\Sx_{k+1,ii} = \min(\gamma_{\rm upd} \Sx_{k+1,ii}, \gamma_{\rm max})\), but this is not beneficial in practice.
				Not only does it not seem to speed up convergence on average, but every change in \(\Sx\) also forces QPALM to refactorize the system.
			
			

		\subsection{Termination}\label{subsec:termination}

			This section discusses the termination criteria used in QPALM.
			Additionally to the criteria to determine a stationary point, we also discuss how to determine whether the problem is primal or dual infeasible.
			
			\subsubsection{Stationarity}
				Termination is based on the unscaled residuals, that is the residuals pertaining to the original problem \eqref{eq:P}.
				In QPALM, we allow for an absolute and a relative tolerance for both the primal and dual residual.
				As such, we terminate on an approximate stationary primal-dual pair \((\bar x, \bar y)\), with associated \({\bar z}^{k} {}={} \proj_C(\bar A {\bar x} + \Sy_k^{-1} {\bar y})\), if
				\begin{subequations}\label{eq:termination}
					\begin{align}
						\tfrac{1}{c}\|D^\smallinv(\bar Q {\bar x} {+} \bar q {+} \trans {\bar A} {\bar y})\|_\smallinf
					{}\leq{} &
						\varepsilon_{\textrm{a}} {+} \tfrac{\varepsilon_\textrm{r}}{c}\max(\|D^\smallinv\bar Q {\bar x}\|_\smallinf, \|D^\smallinv \bar q\|_\smallinf, \|D^\smallinv\trans{\bar A} {\bar y}\|_\smallinf)
					\label{eq:dual termination}
					\\
						\|E^\smallinv(\bar A {\bar x} {-} {\bar z}^k)\|_\smallinf &{}\leq{} \varepsilon_{\textrm{a}} {+} \varepsilon_\textrm{r}\max(\|E^\smallinv\bar A {\bar x}\|_\smallinf, \|E^\smallinv {\bar z}^k\|_\smallinf).
					\label{eq:primal termination}
					\end{align}
				\end{subequations}
				Here, the tolerances \(\varepsilon_{\textrm{a}}\) and \(\varepsilon_{\textrm{r}}\) are by default \(10^{-4}\) and can be chosen by the user.
				In the simulations of \cref{sec:Simulations}, these tolerances were always set to \(10^{-6}\).
			
				To determine termination of the subproblem in \cref{state:PALM:x}, following \eqref{eq:grad phi_j}, the termination criterion
				\begin{equation}\label{eq:inner termination}\mathtight\textstyle
					\frac{1}{c}\|D^\smallinv (\bar Q {\bar x} {+} \bar q {+} \Sx^\smallinv({\bar x} - {\hat{\bar x}}^k) {+} \trans {\bar A} {\bar y})\|_\smallinf 
				{}\leq{}  
					\delta_{\textrm{a},k} {+} \frac{\delta_{\textrm{r},k}}{c}\max(\|D^\smallinv\bar Q {\bar x}\|_\smallinf, \|D^\smallinv \bar q\|_\smallinf, \|D^\smallinv\trans{\bar A} {\bar y}\|_\smallinf)
				\end{equation}
				is used.
				Here, the absolute and relative intermediate tolerances \(\delta_{\textrm{a},k}\) and \(\delta_{\textrm{r},k}\) start out from \(\delta_{\textrm{a},0}\) and \(\delta_{\textrm{r},0}\), which can be set by the user and default to \(10^0\).
				In \cref{state:PALM:delta} they are updated using the following rule
				\begin{align*}
					\delta_{\textrm{a},k+1} &{}={} \max(\rho \delta_{\textrm{a},k+1}, \varepsilon_{\textrm{a}}), \\
					\delta_{\textrm{r},k+1} &{}={} \max(\rho \delta_{\textrm{r},k+1}, \varepsilon_{\textrm{r}}),
				\end{align*}
				with \(\rho\) being the tolerance update factor, which can be set by the user and which defaults to \(10^{-1}\).
				Note that, in theory, these intermediate tolerances should not be lower bounded but instead go to zero.
				In practice, this is however not possible due to machine accuracy errors.
				Furthermore, we have not perceived any inhibition on the convergence as a result of this lower bound.
				This makes sense as the inner subproblems are solved up to machine accuracy by the semismooth Newton method as soon as the correct active set is identified.

			\subsubsection{Infeasibility detection}
				Detecting infeasibility of a (convex) QP from the primal and dual iterates has been discussed in the literature \cite{banjac2019infeasibility}.
				The relevant criteria have also been implemented in QPALM, with a minor modification of the dual infeasibility criterion for a nonconvex QP.
				As such, we determine that the problem is primal infeasible if for a \(\delta {\bar y} \neq 0\) the following two conditions hold
				\begin{subequations}\label{eq:primal infeasibility}
					\begin{align}
						\|D^\smallinv \trans {\bar A} \delta {\bar y}\|_\smallinf &{}\leq{} \varepsilon_{\rm pinf}\|E \delta {\bar y}\|_\smallinf,
					\\
						\trans {\bar u} [\delta {\bar y}]_+ - \trans {\bar l} [-\delta {\bar y}]_+ &{}\leq -\varepsilon_{\rm pinf} \|E \delta {\bar y}\|_\smallinf,
						\label{eq:primal infeas 2}
					\end{align}
				\end{subequations} 
				with the certificate of primal infeasibility being \(\frac{1}{c}E\delta {\bar y}\).
			
				The problem is determined to be dual infeasible if for a \(\delta {\bar x} \neq 0\)
				\begin{subequations} \label{eq:dual infeasibility}
					\begin{equation}
						(E^\smallinv \bar A \delta {\bar x})_i
						\begin{ifcases}
							\in [-\varepsilon_{\rm dinf}, \varepsilon_{\rm dinf}] \|D \delta {\bar x}\|_\smallinf & \bar u_i, \bar \ell_i \in \R
						\\
							\geq -\varepsilon_{\rm dinf} \|D \delta {\bar x}\|_\smallinf & \bar{u}_i = +\infty
						\\
							\leq \varepsilon_{\rm dinf} \|D \delta {\bar x}\|_\smallinf & \bar{\ell}_i = -\infty
						\end{ifcases}
					\end{equation}
					holds for all \(i \in [1,m]\), and either 
					\begin{align} \label{eq:convex dinf}
						\|D^\smallinv \bar Q \delta {\bar x}\|_\smallinf &{}\leq{} c\varepsilon_{\rm dinf}\|D\delta {\bar x}\|_\smallinf,
					\\
						\trans {\bar q} \delta {\bar x} &{}\leq -c\varepsilon_{\rm dinf} \|D \delta {\bar x}\|_\smallinf
						\label{eq:convex dinf 2}
					\end{align}
					or
					\begin{equation}\label{eq:nonconvex dinf}
						\trans{(\delta {\bar x})} \bar Q \delta {\bar x} \leq -c\varepsilon_{\rm dinf}^2\|\delta {\bar x}\|^2
					\end{equation}
					hold.
					Equations \eqref{eq:convex dinf} and \eqref{eq:convex dinf 2} express the original dual infeasibility for convex QPs, that is the conditions that \(\delta {\bar x}\) is a direction of zero curvature and negative slope, whereas \eqref{eq:nonconvex dinf} is added in the nonconvex case to determine whether the obtained \(\delta {\bar x}\) is a direction of negative curvature.
					In the second case, the objective would go to \(-\infty\) quadratically along \(\delta {\bar x}\), and in the first case only linearly.
					Therefore, the square of the tolerance, assumed to be smaller than one, is used in \eqref{eq:nonconvex dinf}, so as to allow for earlier detection of this case.
					Note that we added minus signs in equations \eqref{eq:primal infeas 2} and \eqref{eq:convex dinf 2} in comparison to \cite{banjac2019infeasibility}.
					The reason for this is that the interpretation of our tolerance is different.
					In essence, \cite{banjac2019infeasibility} may declare some problems infeasible even though they are feasible.
					Our version prevents such false positives at the cost of requiring sufficient infeasibility and possibly a slower detection.
					We prefer this, however, over incorrectly terminating a problem, as many interesting problems in practice may be close to infeasible.
					When the tolerances are very close to zero, of course both versions converge to the same criterion.
				\end{subequations}
				The tolerances \(\varepsilon_{\rm pinf}\) and \(\varepsilon_{\rm dinf}\) can be set by the user and have a default value of \(10^{-5}\).

	\section{The full QPALM algorithm}\label{sec:FullQPALM}

		\Cref{alg:FullQPALM} synopsizes all steps and details that make up the QPALM algorithm.
		Herein we set \(\varepsilon_{{\rm a},0} = \delta_{{\rm a}, 0}\), \(\varepsilon_{{\rm r},0} = \delta_{{\rm r}, 0}\).
		For brevity, the details on factorizations and updates necessary for \cref{state:QPALM:semismooth Newton}, which have been discussed prior in \cref{subsec:Solving linear systems}, have been omitted here.
		
		\begin{algorithm}[p]
			\algcaption{QPALM for the nonconvex problem \eqref{eq:P}}
			\label{alg:FullQPALM}

		\begin{algorithmic}[1]\linespread{1.49}\selectfont
		\Require
			\begin{tabular}[t]{@{}l@{}}
				Problem data: \(Q \in \symm(\R^n)\);~~
				\(q \in \R^n\);~~
				\(A \in \R^{m \times n}\);~~
				\(\ell, u \in R^m\) with \(\ell \leq u\);~~ 
				\\		
				\((x^0,y^0)\in\R^n\times\R^m\);~~
				\(\varepsilon_{\rm a},\varepsilon_{\rm r}, \delta_{\rm a,0}, \delta_{\rm r,0}, \varepsilon_{\rm pinf}, \varepsilon_{\rm dinf}, { \varepsilon_{{\rm a},0}, \varepsilon_{{\rm r},0},} \sigma_{\rm init}, \sigma_{\rm max}, \gamma >0\);~~
				\(\rho,\theta\in(0,1)\);~~
					\\
				\(\mathbf{\Delta} > 1\);~~ 
				\(\texttt{scaling} \in \N\)
			\end{tabular}
			\State\label{state:QPALM:scaling}%
				\begin{tabular}[t]{@{}l@{}}
					Use \cref{alg:ruiz} to find \(D\) and \(E\), and let \(c = \max(1.0, \|D(Q x^0 + q)\|_\infty)^{-1}\)
					Convert the data using
				\\
					the scaling factors: \({\bar x}^0 = D^{-1}x^0\), \({\bar y}^0 = cE^{-1}y^0\), \(\bar Q = cDQD\), \(\bar q = cDq\), \(\bar A = EAD\), \(\bar \ell = E\ell\) and \(\bar u = Eu\)%
				\end{tabular}
			\State\label{state:QPALM:initialization}%
				Initialize \({\hat {\bar x}}^0 = {\bar x}^0\), \(\Sy_0\) from \eqref{eq:sigma init} and \(\delta \bar x = 0\)
			\State*\label{state:QPALM:init Sx}%
				Compute \(\lambda^*\) using \cref{alg:lobpcg}
			\If*{ \(\lambda^* < 0\) }
				\State*\label{state:QPALM:nonconvex Sx}%
				\(\Sx_{ii} = \frac{1}{|\lambda^*-10^{-6}|}, \quad i = 1, \ldots, n\)
			\Else*{}
				\State\label{state:QPALM:convex Sx}%
				\(\Sx_{ii} = \gamma, \quad i = 1, \ldots, n\)
			\EndIf{}
		\For{ \(k=0,1,\ldots\) }%
			\State
				Set \({\bar x}^{k, 0} = {\bar x}^{k}\)
			\For{ \(\nu=0,1,\ldots\) }%
				\State
					\(
						{\bar z}^{k,\nu} 
					{}={} 
						\proj_{\bar C}(\bar A {\bar x}^{k, \nu} + \Sy_k^{-1} {\bar y}^k)
					\)
				\State
					\( \delta {\bar y} = \Sy_k(\bar A {\bar x}^{k, \nu}-{\bar z}^{k,\nu})  \)
				\If{ \eqref{eq:termination} is satisfied at \(({\bar x}^{k, \nu}, {\bar y}^k \review{{}+ \delta\bar y})\) }
					\State
						\Return \(({\bar x}^{k, \nu}, {\bar y}^k + \delta {\bar y})\)
				\ElsIf{ \eqref{eq:primal infeasibility} is satisfied at \(\delta {\bar y}\) }
					\State
						\Return \(c^\smallinv E\delta {\bar y}\) as the certificate of primal infeasibility
				\ElsIf{ \eqref{eq:dual infeasibility} is satisfied at \(\delta {\bar x}\) }
					\State
						\Return \(D\delta {\bar x}\) as the certificate of dual infeasibility
				\ElsIf{ \eqref{eq:inner termination} is satisfied at \(({\bar x}^{k, \nu}, {\bar y}^k + \delta {\bar y})\) }
					\State\textbf{break}
				\Else{}
					\State\label{state:QPALM:semismooth Newton}%
						Find \(d\) by solving either \eqref{eq:Newton_SCHUR} or \eqref{eq:Newton_KKT}
					\State\label{state:QPALM:exact linesearch}%
						Find \(\tau\) using \cref{alg:LS}
					\State\label{state:QPALM:update delta x}%
						\(\delta {\bar x} = \tau d\)
					\State\label{state:QPALM:update x inner}%
						\({\bar x}^{k, \nu+1} = {\bar x}^{k, \nu} + \delta {\bar x}\)
				\EndIf{}
			\EndFor{}
			\State\label{state:QPALM:update x outer}%
				Set \({\bar x}^{k+1} = {\bar x}^{k, \nu}\), \({\bar z}^{k+1} = {\bar z}^{k, \nu}\) and \({\bar y}^{k+1} 
				{}={}
					{\bar y}^{k}
					{}+{}
					\delta \bar y\)
			\If*{ { \(\|E^\smallinv (\bar A {\bar x}^{k+1}-{\bar z}^{k+1})\|_\smallinf\leq \varepsilon_{{\rm a}, k} +  \varepsilon_{{\rm r},k}\max(\|E^\smallinv\bar A {\bar x}^{k+1}\|_\smallinf, \|E^\smallinv {\bar z}^{k+1}\|_\smallinf)\)} }\label{step:QPALM:innerTermination}%
				\State\label{state:QPALM:inner success}%
					Update \({\hat {\bar x}}^{k+1} ={\bar x}^{k+1}\)
				\State*
					\(\varepsilon_{{\rm a},{k+1}} = \max(\rho\varepsilon_{{\rm a},{k}}, \varepsilon_{\rm a})\)
					and \(\varepsilon_{{\rm r},{k+1}} = \max(\rho\varepsilon_{{\rm r},{k}}, \varepsilon_{\rm r})\)
			\Else*{}
				\State*\label{state:QPALM:inner fail}%
					{ Set~
					\({\hat {\bar x}}^{k+1}={\hat {\bar x}}^{k}\),
					\(\varepsilon_{{\rm a},{k+1}} = \varepsilon_{{\rm a},{k}}\) 
					and \(\varepsilon_{{\rm r},{k+1}} = \varepsilon_{{\rm r},{k}}\)}
			\EndIf{}
			\State\label{state:QPALM:update Sy}%
				Update \(\Sy_{k+1}\) according to \eqref{eq:update sigma}
			\State\label{state:QPALM:delta}%
				\(\delta_{{\rm a},{k+1}} = \max(\rho\delta_{{\rm a},{k}}, \varepsilon_{\rm a})\) 
					and \(\delta_{{\rm r},{k+1}} = \max(\rho\delta_{{\rm r},{k}}, \varepsilon_{\rm r})\)
		\EndFor{}
		\end{algorithmic}
		\end{algorithm}
		
		It is interesting to note that QPALM algorithm presented here differs from its antecedent convex counterpart \cite{hermans2019qpalm} only by the addition of the lines marked with a star ``\(\star\)'', namely for the setting of \(\Sx\) and the inner termination criteria.
		In the convex case, the starred lines are ignored and \cref{state:QPALM:convex Sx} and \cref{state:QPALM:inner success} will always activate.
		It is clear that the routines in QPALM require minimal changes when extended to nonconvex QPs.
		Furthermore, in numerical experience with nonconvex QPs the criterion of \cref{step:QPALM:innerTermination} seemed to be satisfied most of the time.
		Therefore, aside from the computation of a lower bound of the minimum eigenvalue of \(Q\), QPALM behaves in a very similar manner for convex and for nonconvex QPs.
		Nevertheless, in practice convergence can be quite a bit slower due to the (necessary) heavy regularization induced by \(\Sx\) if Q has a negative eigenvalue with a relatively large magnitude.
		
		\Cref{tab:default parameters} lists the main \review{user-settable} parameters used in QPALM alongside their default values.
		
		\begin{table}
			\centering
			\begin{tabular}{l|cl} 
				Name & Default value & \multicolumn{1}{c}{Description}\\ \hline

			\(\varepsilon_{\rm a}\) & \(10^{-4}\) & Absolute tolerance on termination criteria \\
			\(\varepsilon_{\rm r}\) & \(10^{-4}\) & Relative tolerance on termination criteria \\
			\(\delta_{{\rm a},0}\)  & \(10^0\)    & Starting value of the absolute intermediate tolerance \\
			\(\delta_{{\rm r},0}\)  & \(10^0\)    & Starting value of the relative intermediate tolerance \\
			\(\rho\) 				& \(10^{-1}\) & Update factor for the intermediate tolerance \\
			\(\sigma_{\rm init}\)   & \(20\)	  & Used in the determination of the starting penalty parameters (cf. \eqref{eq:sigma init}) \\
			\(\sigma_{\max}\)			& \(10^9\)    & Cap on the penalty parameters \\
			\(\mathbf{\Delta}\)		& \(100\)	  & Factor used in updating the penalty parameters (cf. \eqref{eq:update sigma}) \\
			\(\theta\)				& \(0.25\)	  & Used in determining which penalties to update (cf. \eqref{eq:update sigma}) \\
			\(\gamma_{\rm init}\)   & \(10^7\)    & Initial value of the proximal penalty parameter (convex case) \\
			\(\gamma_{\rm upd}\)    & \(10\)	  & Update factor for the proximal penalty parameter (convex case) \\
			\(\gamma_{\rm max}\)	& \(10^7\)	  & Cap on the proximal penalty parameter (convex case) \\
			\texttt{scaling}		& \(10\)	  & Number of Ruiz scaling iterations applied to \(A\)
			\end{tabular}
			\caption{Main parameters used in QPALM and their default values.}
			\label{tab:default parameters}%
		\end{table}

	\section{Numerical Results}\label{sec:Simulations}

		The performance of QPALM is benchmarked against other state-of-the-art solvers.
		For convex QPs, we chose the interior-point solver Gurobi \review{(version 9.1.2)} \cite{gurobi2018gurobi}, the \review{(parametric)} active-set solver qpOASES \review{(version 3.2.1)} \cite{ferreau2014qpoases}, and the operator splitting based solver OSQP \review{(version 0.6.2)} \cite{stellato2020osqp}.
		\review{In addition, for the simulations on optimal control we added the tailored interior-point solver HPIPM (version 0.1.4) \cite{frison2020hpipm}.} There are many other solvers available, some of which are tailored to certain problem classes, but the aforementioned ones provide a good sample of the main methods used for general convex QPs.
		For nonconvex QPs, however, no state-of-the-art open-source (local) optimization solver exists to our knowledge.
		Some \review{open-source} indefinite QP algorithms have been proposed, such as in \cite{absil2007newton}.
		However, their solver was found to run into numerical issues very often.
		\review{The active-set solvers SQIC \cite{gill2015methods} and qpOASES \cite{ferreau2014qpoases} also work on indefinite QPs, although the former is not publicly available and the latter fails on most large-scale sparse problems in our benchmarks.} Hence, we did not compare against a QP solver specifically, but rather against a state-of-the-art nonlinear optimization solver, IPOPT \cite{wachter2006implementation}, when dealing with nonconvex QPs.
		All simulations were performed on a notebook with Intel(R) Core(TM) i7-7600U CPU @ 2.80GHz x 2 processor and 16 GB of memory.
		\review{The convex problems are solved to a low and a medium-high accuracy value, with the termination tolerances \(\varepsilon_{\rm a}, \varepsilon_{\rm r}\) both set to \(10^{-3}\) or \(10^{-6}\) respectively for QPALM.
		In other solvers, the corresponding termination tolerances were similarly set to \(10^{-3}\) or \(10^{-6}\).
		This inevitably causes a bit of bloating in the presentation of the results, but it is the fairest way to compare the different algorithms, since the performance of certain solvers varies widely for different accuracies.
		For example, OSQP, being a first-order method, tends to find solutions at low accuracy quickly, but often fails to find solutions at high accuracy in a competitive time, as will be shown.
		The nonconvex problems were solved to an accuracy of \(10^{-6}\) since we only compare against an interior-point method.} Furthermore, for all solvers and all problems, the maximum number of iterations was set to infinity, and a time limit of 3600 seconds was specified.

		\subsection{Comparing runtimes}

			Comparing the performance of solvers on a benchmark test set is not straightforward, and the exact statistics used may influence the resulting conclusions greatly.
			In this paper, we will compare runtimes of different solvers on a set of QPs using two measures, the shifted geometric means (sgm) and the performance profiles.
			When dealing with specific problem classes, such as in \cref{subsec:Portfolio} and \cref{subsec:MPC}, we will not use these statistics but instead make a simple plot of the runtime of the various solvers as a function of the problem dimension.

			\subsubsection{Shifted geometric means}
				Let \(t_\texttt{s,p}\) denote the time required for solver \texttt{s} to solve problem \texttt{p}.
				Then the shifted geometric means \(\bar{t}_\texttt{s}\) of the runtimes for solver \texttt{s} on problem set \texttt{P} is defined as
				\[
					\bar{t}_\texttt{s} = \sqrt[\texttt{\small |P|}]{\prod_{\texttt{p} \in \texttt{P}} (t_\texttt{s,p} + \zeta)}  - \zeta = e^{\frac{1}{\texttt{|P|}}\sum_{\texttt{p} \in \texttt{P}}\ln{(t_\texttt{s,p} + \zeta)}} - \zeta,
				\]
				where the second formulation is used in practice to prevent overflow when computing the product.
				In this paper, runtimes are expressed in seconds, and a shift of \(\zeta = 1\) is used.
				Also note that we employ the convention that when a solver \texttt{s} fails to solve a problem \texttt{p} (within the time limit), the corresponding \(t_\texttt{s,p}\) is set to the time limit for the computation of the sgm.

			\subsubsection{Performance profile}
				To compare the runtime performance in more detail, also performance profiles \cite{dolan2002benchmarking} are used.
				Such a performance profile plots the fraction of problems solved within a runtime of \(f\) times the runtime of the fastest solver for that problem.
				Let \texttt{S} be the set of solvers tested, then
				\[
					r_\texttt{s,p} = \frac{t_\texttt{s,p}}{\min_{\texttt{s} \in \texttt{S}}t_\texttt{s,p}},
				\]   
				denotes the performance ratio of solver \texttt{s} with respect to problem \texttt{p}.
				Note that by convention \(r_\texttt{s,p}\) is set to \(\infty\) when \texttt{s} fails to solve \texttt{p} (within the time limit).
				The fraction of problems \(q_\texttt{s}(f)\) solved by \texttt{s} to within a multiple \(f\) of the best runtime, is then given as
				\[
					q_\texttt{s}(f) = \frac{1}{\texttt{|P|}} \sum_{\texttt{P} \ni \texttt{p} : r_\texttt{s,p} \leq f} 1.
				\] 
				Performance profiles have been found to misrepresent the performance when more than two solvers were compared at the same time \cite{gould2016note}.
				As such, we will construct only the performance profile of each other solver and QPALM, and abstain from comparing the other solvers amongst each other.

		\subsection{Nonconvex QPs}	

			Nonconvex QPs arise in several application domains, such as in a reformulation of mixed integer quadratic programs and in the solution of partial differential equations.
			Furthermore, an indefinite QP has to be solved at every iteration of a sequential quadratic programming method applied to a nonconvex optimization problem.
			To have a broad range of sample QPs, we consider in this paper the set of nonconvex QPs included in the Cutest test framework \cite{gould2015cutest}.
			\Cref{tab:cutest_long} lists for each of those QPs the number of primal variables \(n\) and the number of constraints \(m\), excluding bound constraints.
			In addition, it lists a comparison of the runtime and final objective value for both QPALM and IPOPT.
			Given that both solvers only produce an (approximate) stationary point, and not necessarily the same, these results have been further analyzed to produce \cref{tab:cutest_summary}.
			Here, the problems have been divided according to whether both solvers converged to the same point or not, the criterion of which was set to a relative error on the primal solutions of \(10^{-6}\).
			
\begin{figure}[ht]
				\centering
			  	\includetikz[width=0.7\textwidth]{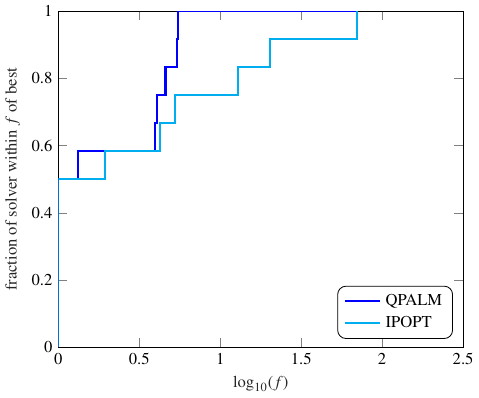}
				\caption{%
					Performance profile for QPALM and IPOPT on the nonconvex QPs of the Cutest test set where both converged to the same approximate stationary point.%
				}%
			    \label{fig:PPcutest}%
			\end{figure}

			On the one hand, the runtimes of the problems where the same solution was found have been listed as shifted geometric means.
			It is clear that \review{on average QPALM is competitive against IPOPT} for these problems.
			These runtimes were further compared in the performance profile of \cref{fig:PPcutest}.
			This shows \review{again that QPALM was competitive against IPOPT in runtimes}.
			On the other hand, for the problems with different solutions, the objective value of the solution was compared and the number of times either QPALM or IPOPT had the lowest objective was counted.
			The resulting tally of \review{46 against 38 in favor} of QPALM suggests there is no clear winner in this case.
			This was to be expected as both solvers report on the first stationary point obtained, and neither uses globalization or restarting procedures to obtain a better one.
			
			\review{The term dead points here refers to first-order stationary points which do not satisfy the second-order necessary condition that the reduced Hessian \(\trans ZQZ\) be positive semidefinite, see \cite[Result 2.2]{gill2015methods}.
			Here, \(Z\) is the nullspace of \(A_{\J\cdot}\), with \(\J\) the set of active constraints, that is the constraints which hold as equality at the given point.
			It is clear from the table that both solvers find a few amount of dead points, although they do so on different problems.}
			
			Finally, also the failure rate was reported.
			It is clear that QPALM outperforms IPOPT by a small margin.
			Furthermore, for the six problems that QPALM failed to solve within the time limit, that is NCVXQP\{1-3,7-9\}, IPOPT also failed to solve in time.
			IPOPT reported two of the problems, A2NNDNIL and A5NNDNIL, as primal infeasible, whereas for these problems QPALM found a point satisfying the approximate stationary conditions.
			In fact, the problems are primal infeasible, and QPALM also reports this once slightly stricter termination tolerances are enforced.
			Hence, we consider both solvers to have succeeded for these two cases.
			
			{\small
				\begin{longtable}[ht]{lcc|cc|cc}
					& & & \multicolumn{2}{c}{Runtime} & \multicolumn{2}{c}{Objective} \\
					\multicolumn{1}{l}{Problem} & n & m & QPALM & IPOPT & QPALM & IPOPT \\
					\hline

			A0ENDNDL & 45006 & 15002 & 1.06e+01 & 1.92e+00 & \phantom{-}1.71e-05 & \phantom{-}1.84e-04\\ 
			A0ENINDL & 45006 & 15002 & 9.85e+00 & 1.83e+00 & \phantom{-}4.61e-05 & \phantom{-}1.84e-04\\ 
			A0ENSNDL & 45006 & 15002 & 4.68e+00 & 3.55e+01 & -9.40e-06 & \phantom{-}1.48e-04\\ 
			A0ESDNDL & 45006 & 15002 & 1.06e+01 & 2.66e+00 & \phantom{-}5.47e-06 & \phantom{-}1.84e-04\\ 
			A0ESINDL & 45006 & 15002 & 9.34e+00 & 2.03e+00 & \phantom{-}2.21e-06 & \phantom{-}1.84e-04\\ 
			A0ESSNDL & 45006 & 15002 & 4.42e+00 & 2.52e+01 & -6.66e-06 & \phantom{-}1.48e-04\\ 
			A0NNDNDL & 60012 & 20004 & 1.76e+02 & 7.88e+00 & \phantom{-}3.22e-04 & \phantom{-}1.84e-04\\ 
			A0NNDNIL & 60012 & 20004 & 3.55e+03 & 3.03e+01 & \phantom{-}2.25e+00 & \phantom{-}1.98e-04\\ 
			A0NNDNSL & 60012 & 20004 & 6.26e+01 & 2.56e+01 & -4.81e-04 & \phantom{-}2.15e-04\\ 
			A0NNSNSL & 60012 & 20004 & 1.87e+01 & 3.26e+01 & -1.86e-04 & \phantom{-}1.54e-04\\ 
			A0NSDSDL & 60012 & 20004 & 3.00e+01 & 5.83e+00 & -2.29e-04 & \phantom{-}1.84e-04\\ 
			A0NSDSDS &  6012 &  2004 & 1.49e+00 & 9.31e-01 & -6.90e-06 & \phantom{-}2.33e-04\\ 
			A0NSDSIL & 60012 & 20004 & 6.54e+02 & 3.00e+01 & -1.77e-04 & \phantom{-}1.97e-04\\ 
			A0NSDSSL & 60012 & 20004 & 3.11e+01 & 1.29e+01 & -3.86e-06 & \phantom{-}1.67e-04\\ 
			A0NSSSSL & 60012 & 20004 & 1.77e+01 & 3.02e+01 & \phantom{-}5.65e-05 & \phantom{-}1.49e-04\\ 
			A2ENDNDL & 45006 & 15002 & 2.63e+01 & 2.76e+00 & \phantom{-}1.04e-06 & \phantom{-}9.88e-04\\ 
			A2ENINDL & 45006 & 15002 & 2.53e+01 & 2.77e+00 & \phantom{-}7.94e-07 & \phantom{-}9.73e-04\\ 
			A2ENSNDL & 45006 & 15002 & 4.68e+00 & 4.11e+01 & \phantom{-}9.99e-07 & \phantom{-}2.10e-02\\ 
			A2ESDNDL & 45006 & 15002 & 2.75e+01 & 2.77e+00 & \phantom{-}1.22e-06 & \phantom{-}9.88e-04\\ 
			A2ESINDL & 45006 & 15002 & 2.52e+01 & 2.81e+00 & \phantom{-}4.34e-07 & \phantom{-}9.73e-04\\ 
			A2ESSNDL & 45006 & 15002 & 4.83e+00 & 3.91e+01 & \phantom{-}2.03e-06 & \phantom{-}2.22e-02\\ 
			A2NNDNDL & 60012 & 20004 & 1.43e+03 & 7.48e+00 & \phantom{-}3.44e-04 & \phantom{-}3.03e-04\\ 
			A2NNDNIL & 60012 & 20004 & 3.14e+01 & PI & \phantom{-}5.11e+01 & /\\ 
			A2NNDNSL & 60012 & 20004 & 3.48e+02 & 8.85e+01 & -9.17e-07 & \phantom{-}2.11e-04\\ 
			A2NNSNSL & 60012 & 20004 & 2.89e+01 & 3.17e+01 & -3.59e-05 & \phantom{-}1.46e-03\\ 
			A2NSDSDL & 60012 & 20004 & 5.15e+01 & 6.14e+00 & \phantom{-}3.13e-06 & \phantom{-}7.76e-04\\ 
			A2NSDSIL & 60012 & 20004 & 4.17e+01 & 4.20e+01 & \phantom{-}5.95e-01 & \phantom{-}1.57e+00\\ 
			A2NSDSSL & 60012 & 20004 & 4.08e+01 & 2.02e+01 & -4.01e-06 & \phantom{-}2.47e-02\\ 
			A2NSSSSL & 60012 & 20004 & 2.55e+01 & 2.94e+01 & -8.81e-06 & \phantom{-}4.69e-04\\ 
			A5ENDNDL & 45006 & 15002 & 4.97e+01 & 2.73e+00 & \phantom{-}1.01e-06 & \phantom{-}2.19e-03\\ 
			A5ENINDL & 45006 & 15002 & 5.38e+01 & 2.75e+00 & \phantom{-}4.48e-07 & \phantom{-}2.25e-03\\ 
			A5ENSNDL & 45006 & 15002 & 5.06e+00 & 3.72e+01 & \phantom{-}4.15e-06 & \phantom{-}5.12e-02\\ 
			A5ESDNDL & 45006 & 15002 & 5.25e+01 & 2.76e+00 & \phantom{-}2.05e-06 & \phantom{-}2.19e-03\\ 
			A5ESINDL & 45006 & 15002 & 5.14e+01 & 2.74e+00 & \phantom{-}6.23e-07 & \phantom{-}2.25e-03\\ 
			A5ESSNDL & 45006 & 15002 & 4.94e+00 & 3.68e+01 & \phantom{-}1.86e-05 & \phantom{-}5.12e-02\\ 
			A5NNDNDL & 60012 & 20004 & 2.02e+03 & 9.30e+00 & \phantom{-}5.05e-04 & \phantom{-}1.86e-03\\ 
			A5NNDNIL & 60012 & 20004 & 3.39e+01 & PI & \phantom{-}1.02e+02 & /\\ 
			A5NNDNSL & 60012 & 20004 & 2.27e+02 & 8.88e+01 & \phantom{-}3.00e-06 & \phantom{-}2.10e-04\\ 
			A5NNSNSL & 60012 & 20004 & 3.29e+01 & 4.22e+01 & \phantom{-}3.63e-06 & \phantom{-}1.38e-03\\ 
			A5NSDSDL & 60012 & 20004 & 9.73e+01 & 5.94e+00 & \phantom{-}3.92e-06 & \phantom{-}1.86e-03\\ 
			A5NSDSDM &  6012 &  2004 & 2.13e+00 & 8.69e-01 & \phantom{-}1.14e-06 & \phantom{-}2.33e-04\\ 
			A5NSDSIL & 60012 & 20004 & 3.76e+01 & 6.47e+01 & \phantom{-}6.68e+00 & \phantom{-}1.29e+00\\ 
			A5NSDSSL & 60012 & 20004 & 5.04e+01 & 1.97e+01 & -1.11e-05 & \phantom{-}1.00e-02\\ 
			A5NSSNSM &  6012 &  2004 & 1.52e+00 & 9.42e-01 & \phantom{-}1.90e-06 & \phantom{-}2.33e-04\\ 
			A5NSSSSL & 60012 & 20004 & 3.13e+01 & 7.58e+01 & -4.78e-06 & \phantom{-}2.45e-04\\ 
			BIGGSC4 &     4 &     7 & 1.48e-04 & 5.41e-02 & -2.45e+01 & -2.45e+01\\ 
			BLOCKQP1 & 10010 &  5001 & 1.64e-01 & 1.39e+01 & -4.99e+03 & -4.99e+03\\ 
			BLOCKQP2 & 10010 &  5001 & 1.41e-01 & 1.83e+00 & -4.99e+03 & -4.99e+03\\ 
			BLOCKQP3 & 10010 &  5001 & 2.69e+02 & 4.35e+02 & -2.49e+03 & -2.49e+03\\ 
			BLOCKQP4 & 10010 &  5001 & 3.74e-01 & 1.96e+00 & -2.50e+03 & -2.50e+03\\ 
			BLOCKQP5 & 10010 &  5001 & 2.56e+02 & 4.31e+02 & -2.49e+03 & -2.49e+03\\ 
			BLOWEYA &  4002 &  2002 & 1.86e+00 & 3.46e+03 & -7.09e-06 & -2.28e-02\\ 
			BLOWEYB &  4002 &  2002 & 3.76e-02 & 3.00e+03 & -3.49e-05 & -1.52e-02\\ 
			BLOWEYC &  4002 &  2002 & 6.64e-01 & F & -2.92e-03 & /\\ 
			CLEUVEN3 &  1200 &  2973 & 6.85e+00 & 2.47e+01 & \phantom{-}3.72e+05 & \phantom{-}2.86e+05\\ 
			CLEUVEN4 &  1200 &  2973 & 5.87e+02 & 5.74e+01 & \phantom{-}2.84e+06 & \phantom{-}2.86e+05\\ 
			CLEUVEN5 &  1200 &  2973 & 6.42e+00 & 2.46e+01 & \phantom{-}3.72e+05 & \phantom{-}2.86e+05\\ 
			CLEUVEN6 &  1200 &  3091 & 5.67e+00 & 2.25e+01 & \phantom{-}2.21e+07 & \phantom{-}2.21e+07\\ 
			FERRISDC &  2200 &   210 & 2.14e+00 & 2.58e+00 & -1.02e-10 & -2.13e-04\\ 
			GOULDQP1 &    32 &    17 & 3.88e-03 & 7.84e-02 & -3.49e+03 & -3.49e+03\\ 
			HATFLDH &     4 &     7 & 1.25e-04 & 5.60e-02 & -2.45e+01 & -2.45e+01\\ 
			HS44 &     4 &     6 & 6.92e-05 & 3.54e-02 & -1.50e+01 & -1.30e+01\\ 
			HS44NEW &     4 &     6 & 6.84e-05 & 4.78e-02 & -1.50e+01 & -1.30e+01\\ 
			LEUVEN2 &  1530 &  2329 & 2.67e+00 & 1.96e+00 & -1.41e+07 & -1.41e+07\\ 
			LEUVEN3 &  1200 &  2973 & 1.05e+03 & 2.45e+02 & -1.38e+09 & -1.99e+09\\ 
			LEUVEN4 &  1200 &  2973 & 1.50e+01 & 4.84e+02 & -4.78e+08 & -1.83e+09\\ 
			LEUVEN5 &  1200 &  2973 & 1.05e+03 & 2.62e+02 & -1.38e+09 & -1.99e+09\\ 
			LEUVEN6 &  1200 &  3091 & 2.09e+03 & 1.05e+02 & -1.17e+09 & -1.19e+09\\ 
			LEUVEN7 &   360 &   946 & 9.57e-02 & 4.14e-01 & \phantom{-}6.95e+02 & \phantom{-}6.95e+02\\ 
			LINCONT &  1257 &   419 & PI & PI & / & /\\ 
			MPC1 &  2550 &  3833 & 3.47e+00 & 6.73e+00 & -2.33e+07 & -2.33e+07\\ 
			MPC10 &  1530 &  2351 & 4.82e+00 & 7.33e-01 & -1.50e+07 & -1.50e+07\\ 
			MPC11 &  1530 &  2351 & 4.10e+00 & 8.09e-01 & -1.50e+07 & -1.50e+07\\ 
			MPC12 &  1530 &  2351 & 4.32e+00 & 7.31e-01 & -1.50e+07 & -1.50e+07\\ 
			MPC13 &  1530 &  2351 & 3.65e+00 & 7.39e-01 & -1.50e+07 & -1.50e+07\\ 
			MPC14 &  1530 &  2351 & 3.54e+00 & 7.69e-01 & -1.50e+07 & -1.50e+07\\ 
			MPC15 &  1530 &  2351 & 3.43e+00 & 7.72e-01 & -1.50e+07 & -1.50e+07\\ 
			MPC16 &  1530 &  2351 & 2.92e+00 & 7.53e-01 & -1.50e+07 & -1.50e+07\\ 
			MPC2 &  1530 &  2351 & 3.74e+00 & 6.67e-01 & -1.50e+07 & -1.50e+07\\ 
			MPC3 &  1530 &  2351 & 2.99e+00 & 7.94e-01 & -1.50e+07 & -1.50e+07\\ 
			MPC4 &  1530 &  2351 & 4.41e+00 & 6.74e-01 & -1.50e+07 & -1.50e+07\\ 
			MPC5 &  1530 &  2351 & 4.12e+00 & 7.73e-01 & -1.50e+07 & -1.50e+07\\ 
			MPC6 &  1530 &  2351 & 4.22e+00 & 7.71e-01 & -1.50e+07 & -1.50e+07\\ 
			MPC7 &  1530 &  2351 & 4.17e+00 & 8.05e-01 & -1.50e+07 & -1.50e+07\\ 
			MPC8 &  1530 &  2351 & 3.98e+00 & 7.60e-01 & -1.50e+07 & -1.50e+07\\ 
			MPC9 &  1530 &  2351 & 4.47e+00 & 7.03e-01 & -1.50e+07 & -1.50e+07\\ 
			NASH &    72 &    24 & PI & PI & / & /\\ 
			NCVXQP1 & 10000 &  5000 & F & F & / & /\\ 
			NCVXQP2 & 10000 &  5000 & F & F & / & /\\ 
			NCVXQP3 & 10000 &  5000 & F & F & / & /\\ 
			NCVXQP4 & 10000 &  2500 & 1.06e+03 & 1.90e+03 & -9.38e+09 & -9.38e+09\\ 
			NCVXQP5 & 10000 &  2500 & 1.52e+03 & 2.76e+03 & -6.63e+09 & -6.63e+09\\ 
			NCVXQP6 & 10000 &  2500 & 2.49e+03 & F & -3.40e+09 & /\\ 
			NCVXQP7 & 10000 &  7500 & F & F & / & /\\ 
			NCVXQP8 & 10000 &  7500 & F & F & / & /\\ 
			NCVXQP9 & 10000 &  7500 & F & F & / & /\\ 
			PORTSNQP & 100000 &     2 & 3.80e+02 & 1.31e+00 & -1.56e+00 & -1.00e+00\\ 
			QPNBAND & 50000 & 25000 & 4.48e+00 & 3.41e+00 & -2.50e+05 & -2.50e+05\\ 
			QPNBLEND &    83 &    74 & 4.46e-03 & 5.49e-02 & -9.14e-03 & -9.13e-03\\ 
			QPNBOEI1 &   384 &   351 & 2.09e-01 & 1.57e+00 & \phantom{-}6.78e+06 & \phantom{-}6.75e+06\\ 
			QPNBOEI2 &   143 &   166 & 2.82e-02 & 8.54e-01 & \phantom{-}1.37e+06 & \phantom{-}1.37e+06\\ 
			QPNSTAIR &   467 &   356 & 6.88e-02 & 9.34e-01 & \phantom{-}5.15e+06 & \phantom{-}5.15e+06\\ 
			SOSQP1 &  5000 &  2501 & 4.07e-02 & 1.74e-01 & \phantom{-}4.24e-07 & -1.03e-10\\ 
			SOSQP2 &  5000 &  2501 & 9.44e-02 & 1.80e-01 & -1.25e+03 & -1.25e+03\\ 
			STATIC3 &   434 &    96 & DI & DI & / & /\\ 
			STNQP1 &  8193 &  4095 & 1.51e+01 & 1.05e+03 & -3.12e+05 & -3.12e+05\\ 
			STNQP2 &  8193 &  4095 & 2.99e+01 & 7.38e+00 & -5.75e+05 & -5.75e+05\\ 
					\caption{%
						Runtime and final objective value comparison for QPALM and IPOPT applied to the nonconvex QPs of the Cutest test set.
						Failure codes: PI = primal infeasible, DI = dual infeasible and F = time limit exceeded (or numerical issues).%
					}%
					\label{tab:cutest_long}
				\end{longtable}
			}%
			
			\begin{table}[ht]
				\centering
				\begin{tabular}{l|cc}
					& QPALM & IPOPT
				\\\hline

			Runtime (sgm) & 4.0861 & 4.0761\\
			Optimal &   46 &   38\\
			Dead points &    7 &    4 \\
			Failure rate [\%] & 5.6075 & 7.4766
				\end{tabular}
				\caption{%
					Statistics of QPALM and IPOPT applied to the nonconvex QPs of the Cutest test set.
					The runtime reported is the mean over the 11 problems which converged to the same stationary point, whereas \textit{optimal} denotes the number of times the solver found the stationary point with the lowest objective in problems where different stationary points were found.
					\review{Dead points indicate points at which second-order necessary conditions are violated.}%
				}%
				\label{tab:cutest_summary}%
			\end{table}
			
		\subsection{Convex QPs}

			Convex QPs arise in multiple well-known application domains, such as portfolio optimization and linear MPC.
			Solving such QPs has therefore been the subject of substantial research, and many methods exist.
			We compare QPALM against the interior-point solver Gurobi \cite{gurobi2018gurobi}, the active-set solver qpOASES \cite{ferreau2014qpoases}, and the operator splitting based solver OSQP \cite{stellato2020osqp}.
			First we compare all solvers on the Maros-Meszaros benchmark test set \cite{maros1999repository}.
			However, qpOASES is excluded in this comparison as it tends to fail on larger problems which are ubiquitous in this set.
			Then, the performance of all solvers is also compared for quadratic problems arising from the two aforementioned application domains, portfolio optimization and MPC.
			\review{In the last case, we add the tailored interior-point solver HPIPM \cite{frison2020hpipm} to the benchmark.}

		\subsubsection{Maros Meszaros}

\begin{figure}
			\centering
				\begin{minipage}{.495\linewidth}
					\centering
					\includetikz[width=\linewidth]{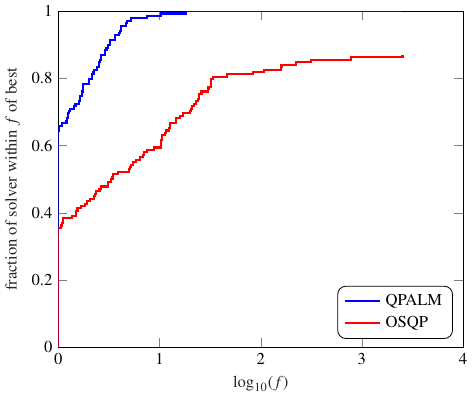}
				\end{minipage}
				\hfill
				\begin{minipage}{.495\linewidth}
					\centering
					\includetikz[width=\linewidth]{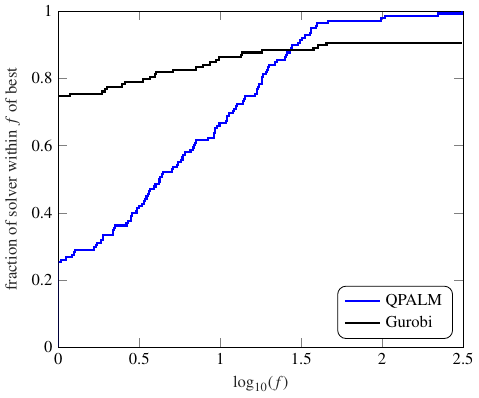}
				\end{minipage}
				\caption{%
					Performance profiles comparing QPALM with OSQP and Gurobi respectively on the Maros Meszaros problem set \review{for a tolerance of \(10^{-6}\)}.%
				}%
				\label{fig:PP_MM_1e6}%
			\end{figure}

			The Maros-Meszaros test set contains 138 convex quadratic programs, and is often used to benchmark convex QP solvers.
			\review{\Cref{tab:MM_small_1e3,tab:MM_small_1e6} list} the shifted geometric mean of the runtime and failure rate of QPALM, OSQP and Gurobi applied to this set \review{with termination tolerances respectively set to \(10^{-3}\) and \(10^{-6}\)}.
			A key aspect of QPALM that is demonstrated here is its robustness.
			The Maros-Meszaros set includes many large-scale and ill-conditioned QPs, and the fact that QPALM succeeds in solving all of them \review{up to an accuracy of \(10^{-6}\)} within one hour is a clear indication that it is very robust with respect to the problem data.
			In runtime it is also faster on average than the other solvers.
			However, Gurobi is faster more often, as is shown in the performance profiles \review{in \cref{fig:PP_MM_1e3,fig:PP_MM_1e6}}.
			The high shifted geometric mean runtime of Gurobi is mostly due to its relatively high failure rate.
			\review{For a tolerance of \(10^{-6}\),} OSQP also has a high failure rate, and is slower than QPALM, both on average and in frequency.
			As a first-order method, in spite of employing a similar preconditioning routine to ours, it seems to still exhibit a lack of robustness with respect to ill conditioning and to somewhat stricter tolerance requirements.
			\review{However, it clearly performs well on the set for a low tolerance of \(10^{-3}\), as demonstrated in \cref{tab:MM_small_1e3,fig:PP_MM_1e3}.
			It only fails once and is faster than QPALM in almost 80\% of the problems.}
			
			\begin{table}
				\centering
			\review{%
				\begin{tabular}{l|ccc}
					& QPALM & OSQP & Gurobi \\ \hline

			Runtime (sgm) & 0.5665 & 0.5943 & 1.4254\\
			Failure rate [\%] & 0.0000 & 0.7246 & 9.4203
				\end{tabular}
				\caption{%
					\review{Shifted geometric mean runtime and failure rate for QPALM, OSQP and Gurobi on the Maros Meszaros problem set with tolerance \(10^{-3}\).}%
				}%
				\label{tab:MM_small_1e3}%
			}%
			\end{table}
			
			\begin{table}
				\centering
				\begin{tabular}{l|ccc}
					& QPALM & OSQP & Gurobi \\ \hline

			Runtime (sgm) & 0.8286 & 7.8153 & 1.4101\\
			Failure rate [\%] & 0.0000 & 13.0435 & 9.4203
				\end{tabular}
				\caption{%
					Shifted geometric mean runtime and failure rate for QPALM, OSQP and Gurobi on the Maros Meszaros problem set \review{with tolerance \(10^{-6}\)}.%
				}%
				\label{tab:MM_small_1e6}%
			\end{table}
			
			\review{%
			\begin{figure}
			\centering
			\review{%
				\begin{minipage}{.495\linewidth}
					\centering
					\includetikz[width=\linewidth]{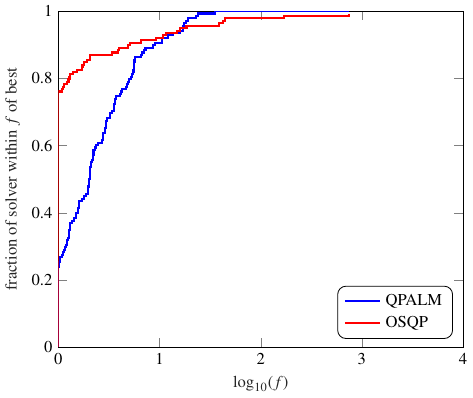}
				\end{minipage}
				\hfill
				\begin{minipage}{.495\linewidth}
					\centering
					\includetikz[width=\linewidth]{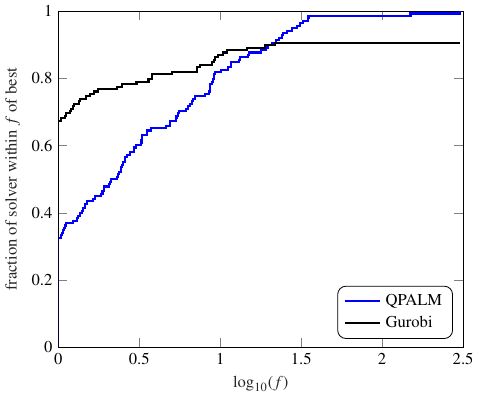}
				\end{minipage}
				\caption{%
					\review{Performance profiles comparing QPALM with OSQP and Gurobi respectively on the Maros Meszaros problem set for a tolerance of \(10^{-3}\).}%
				}%
				\label{fig:PP_MM_1e3}%
			}%
			\end{figure}}%
			
		\subsubsection{Portfolio}\label{subsec:Portfolio}

\begin{figure}
			\centering
				\begin{minipage}{.495\linewidth}
					\centering
					\includetikz[width=\linewidth]{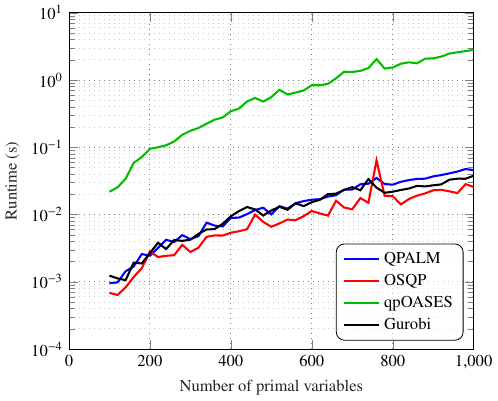}
				\end{minipage}
				\hfill
				\begin{minipage}{.495\linewidth}
					\centering
					\includetikz[width=\linewidth]{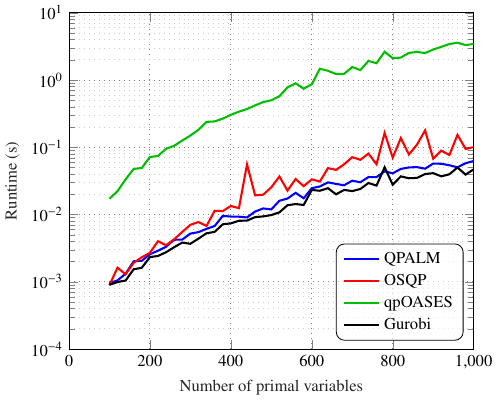}
				\end{minipage}
				\caption{%
					Runtimes of QPALM, OSQP, qpOASES and Gurobi when solving portfolio optimization problems of varying sizes.
					\review{The tolerances are \(10^{-3}\) and \(10^{-6}\) left and right respectively.}%
				}%
				\label{fig:Portfolio}%
			\end{figure}

			In portfolio optimization, the goal is to select a portfolio of assets to invest in to maximize profit taking into account risk levels.
			Given a vector \(x\) denoting the (relative) investment in each asset, the resulting quadratic program is the following
			\begin{align*}
				\minimize_{x\in\R^n} \hspace{0.2cm} & \beta\trans x\Sigma x-\trans \mu x \\
			\stt{} \hspace{0.2cm}
				& x \geq 0, \\
				& \sum_{i=1}^n x_i = 1,
			\end{align*}
			with \(\mu \in \R^{n}\) a vector of expected returns, \(\Sigma \in \symm(\R^n)\) a covariance matrix representing the risk and \(\beta > 0\) a parameter to adjust the aversion to risk.
			Typically, \(\Sigma = F\trans F + D\), with \(F\in \R^{n \times r}\) a low rank matrix and \(D\in \R^{n \times n}\) a diagonal matrix.
			In order not to form the matrix \(\Sigma\), the following reformulated problem can be solved instead in \((x, y)\)
			\begin{align*}
				\minimize_{x\in\R^n} \hspace{0.2cm} & 
			\begin{bmatrix}
				x \\
				y 
			\end{bmatrix}^\top 
			\begin{bmatrix}
				D & \\
				& \Id_r
			\end{bmatrix}		
			\begin{bmatrix}
				x \\
				y 
			\end{bmatrix} - \beta^\smallinv \trans \mu x \\
			\stt{} \hspace{0.2cm}
				& x \geq 0, \\
				& \sum_{i=1}^n x_i = 1, \\
				& y = \trans F x.
			\end{align*}
			We solved this problem for values of \(n\) ranging from 100 to 1000, with \(r = \ceil{\frac{n}{10}}\).
			We choose the elements of \(\mu\) uniformly on \([0,1]\), the diagonal elements \(D_{ii}\)  uniformly on the interval \([0, \sqrt{r}]\), and the matrix \(F\) has \(50\%\) nonzeros drawn from \(\mathcal{N}(0,1)\).
			For each value of \(n\), we solve the problem \review{for ten values of \(\beta\) on a logarithmic scale between \(10^{-2}\) and \(10^2\)} and compute the arithmetic mean of the runtimes.
			\review{The runtimes of QPALM, OSQP, Gurobi and qpOASES solving these problems as such with tolerances \(10^{-3}\) and \(10^{-6}\) for different values of \(n\) are shown in \cref{fig:Portfolio}.
			When warm-starting the problems from the solution of the problem with the previous \(\beta\) value, \cref{fig:Portfolio_sequential} is obtained.} The structure of the portfolio optimization problem is quite specific: the Hessian of the objective is diagonal, and the only inequality constraints are bound constraints.
			It is clear from the figure that Gurobi exhibits the lowest runtimes for this type of problem, followed closely by QPALM and OSQP.
			The latter performs well especially for the small problems and has some robustness issues for larger ones.
			It seems qpOASES exhibits quite a high runtime when compared to the others.
			\review{It is, however, the solver which benefits most from warm-starting in this scenario, as the runtime for warm-started problems is similar to the other solvers.
			Therefore, if very many values for \(\beta\) were tested, the qpOASES curve would coincide with the others in \cref{fig:Portfolio_sequential}.
			QPALM and OSQP already exhibit low runtimes for this problem and do not benefit much from warm-starting here.}
			
			\review{%
			\begin{figure}
			\centering
				\begin{minipage}{.495\linewidth}
					\centering
					\includetikz[width=\linewidth]{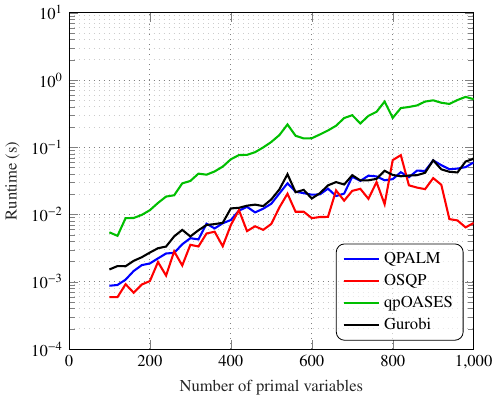}
				\end{minipage}
				\hfill
				\begin{minipage}{.495\linewidth}
					\centering
					\includetikz[width=\linewidth]{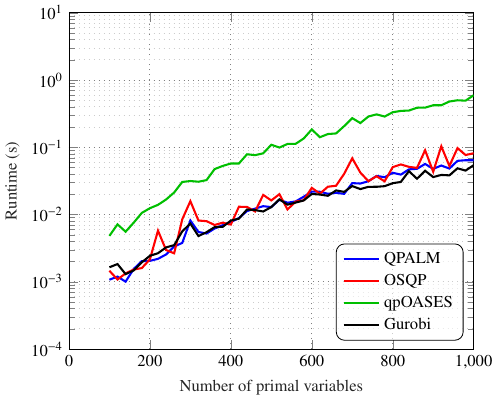}
				\end{minipage}
				\caption{%
					\review{Runtimes of QPALM, OSQP, qpOASES and Gurobi when solving portfolio optimization problems of varying sizes and warm-starting in between different values of \(\beta\).
					The tolerances are \(10^{-3}\) and \(10^{-6}\) left and right respectively.}%
				}%
				\label{fig:Portfolio_sequential}%
			\end{figure}
			}%

		\subsubsection{MPC}\label{subsec:MPC}

			In a \review{model predictive control} strategy, one solves an optimal control problem (OCP) at every sample time to determine the optimal control inputs that need to be applied to a system.
			The OCP considers a control horizon \(N\), that is, it computes a series of \(N\) inputs, of which only the first is applied to the system.
			Given a discrete linear system with \(n_x\) states \(x\) and \(n_u\) inputs \(u\), and its corresponding system dynamics in state-space form, \(x_{k+1} = Ax_k + Bu_k\), the OCP we consider in this paper is one where we control the system from an initial state \(\tilde x\) to the reference state at the origin, which can be formulated as
			\begin{align*}
				\minimize_{z\in\R^{(N+1)n_x + Nn_u}} \hspace{0.2cm} & x_N^\top  Q_N x_N + \sum_{k = 0}^{N-1} x_k^\top Q x_k + u_k^\top R u_k \\
			\stt{} \hspace{0.2cm}
				& x_0 = \tilde x, \\
				& x_{k+1} = Ax_k + Bu_k, \quad k = 0,\ldots, N-1, \\
				& x_k \in \X, \quad k = 0,\ldots, N-1, \\
				& x_N \in \X_N, \\
				& u_k \in \U, \quad k = 0,\ldots, N-1.
			\end{align*}
			
			Here, the decision variable is the collection of \(N+1\) state samples and \(N\) input samples, \(z = (x_0, u_0, x_1, \ldots, u_{N-1}, x_N)\).
			The stage and terminal state cost matrices are positive definite matrices, \(Q, Q_N \in \symm_{++}(\R^{n_x})\) and \(R \in \symm_{++}(\R^{n_u})\).
			\(\X\), \(\X_N\) and \(\U\) represent polyhedral constraints on the states, terminal state and inputs respectively.
			In our example, we consider box constraints \(\X = [-x_b, x_b]\) and \(\U = [-u_b, u_b]\) and determine the terminal constraint as the maximum control invariant set of the system.
			Furthermore, the terminal cost is computed from the discrete-time algebraic Riccati equations.
			
\begin{figure}
			\centering
				\begin{minipage}{.495\linewidth}
					\centering
					\includetikz[width=\linewidth]{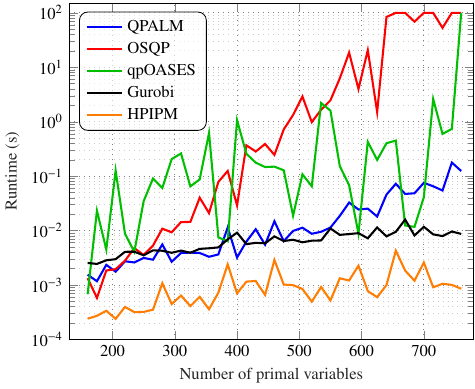}
				\end{minipage}
				\hfill
				\begin{minipage}{.495\linewidth}
					\centering
					\includetikz[width=\linewidth]{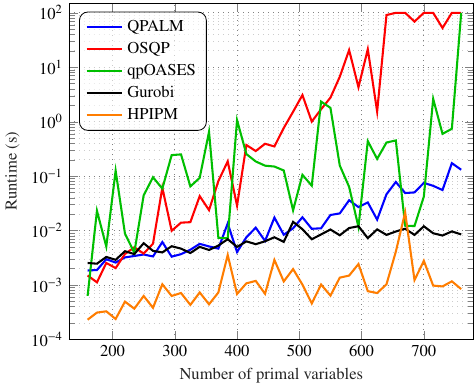}
				\end{minipage}
				\caption{%
					Runtimes of QPALM, OSQP, qpOASES, Gurobi and HPIPM when solving OCPs for varying time horizons.
					\review{The tolerances are \(10^{-3}\) and \(10^{-6}\) left and right respectively.}%
				}%
				\label{fig:randomMPC}%
			\end{figure}

			We solved this problem for a system with 10 states and 5 inputs for different values of the time horizon.
			The state cost matrix is set as \(Q = M\trans M\), with \(M\in \R^{n_x \times n_x}\) consisting of \(50\%\) nonzeros drawn from the normal distribution \(\mathcal N(0,5)\).
			The input cost matrix is chosen to be a small diagonal matrix with \(R_{ii} = 0.01\).
			The system considered is slightly unstable, with the elements of \(A\) drawn from \(\mathcal N(0,2)\) and those of \(B\) from \(\mathcal N(0,1)\).
			The state and input limits \(x_b\) and \(u_b\) are drawn from \(\mathcal N(10,2)\).
			Finally, the initial state is chosen such that it is possible but difficult to satisfy all the constraints, in order to represent a challenging MPC problem.
			The resulting runtimes of solving one such OCP for varying time horizons are shown in \cref{fig:randomMPC} \review{for tolerances \(10^{-3}\) and \(10^{-6}\).
			HPIPM performs best, as expected from a tailored solver, followed by Gurobi and QPALM.}
			OSQP and qpOASES both have issues with robustness given the challenging nature of the problem, although \review{the former performs well on small problems and the latter also exhibits} fast convergence in some cases.
			
			
\begin{figure}
			\centering
				\begin{minipage}{.495\linewidth}
					\centering
					\includetikz[width=\linewidth]{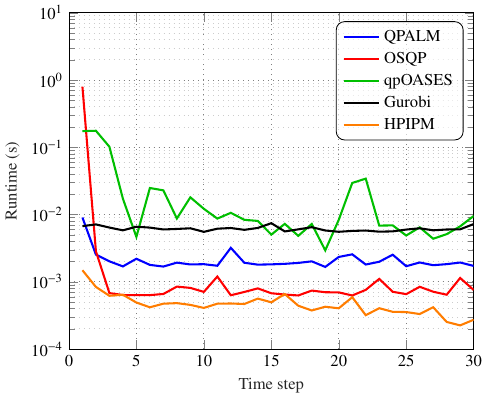}
				\end{minipage}
				\hfill
				\begin{minipage}{.495\linewidth}
					\centering
					\includetikz[width=\linewidth]{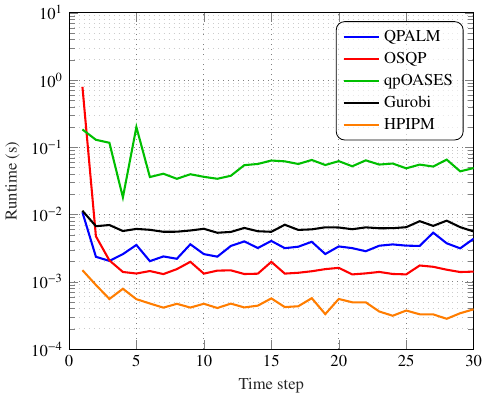}
				\end{minipage}
				\caption{%
					Runtimes of QPALM, OSQP, qpOASES, Gurobi and HPIPM when solving sequential OCPs in an MPC setting, with \(N=30\).
					\review{The tolerances are \(10^{-3}\) and \(10^{-6}\) left and right respectively.}%
				}
				\label{fig:randomMPCsequential}
			\end{figure}

			An important aspect to consider when choosing a QP solver for MPC is the degree to which it can work with an initial guess.
			This is of great import due to the fact that subsequent OCPs are very similar.
			The solution of the previous OCP can therefore be shifted by one sample time and supplied as an initial guess.
			This procedure is also called warm-starting.
			\Cref{fig:randomMPCsequential} shows the result of warm-starting subsequent OCP in this manner.
			Here, we solved 30 subsequent OCPs for a fixed time horizon of 30, corresponding to 460 primal variables.
			Furthermore, when computing the next initial state, we add a small disturbance drawn from the normal distribution \(\mathcal N(0, 0.01)\).
			It is clear that qpOASES, QPALM and OSQP all benefit greatly from this warm-starting.
			\review{However, Gurobi, as is typical of an interior-point method, does not have this advantage.
			For this reason, interior-point methods are typically not considered as solvers for MPC problems.
			HPIPM, however, has incredibly low runtimes regardless, and so may be an excellent choice for optimal control problems.}

	\section{Conclusion}\label{sec:Conclusion}

		This paper presented QPALM, a proximal augmented Lagrangian method for convex and nonconvex quadratic programming.
		On a theoretical level, it is shown that the sequence of inexact solutions of the proximal augmented Lagrangian, shown to be equivalent to inexact proximal point iterations, converges globally at an \(R\)-linear rate to a stationary point for the original problem when the proximal penalty ensures strong convexity of the inner subproblems.
		On a practical level, the implementation of QPALM is considered in great detail.
		The inner subproblems are solved using a direction obtained from a semismooth Newton method which relies on dedicated \(\LDL\)-factorization and factorization update routines, and on the optimal stepsize which can be efficiently computed as the zero of a monotone, piecewise affine function.
		
		The QPALM algorithm is implemented in open-source C code, and parameter selection and update routines have all been worked out carefully.
		The resulting code is shown to strike a unique balance between robustness when faced with hard problems and efficiency when faced with easy problems.
		Given a time limit of one hour, QPALM can find an approximate stationary point or correctly identify infeasibility for 94.39\% of the nonconvex QPs in the Cutest test set, whereas IPOPT does this only for \review{92.52\%}.
		Moreover, QPALM was able to solve all of the convex QPs in the Maros-Meszaros set \review{up to a tolerance of \(10^{-6}\), while} Gurobi and OSQP exhibited a fail rate of \review{9.42\% and 13.04\%,} respectively.
		These results are significant since the Cutest and Maros-Meszaros test-set contain some very large-scale and ill-conditioned QPs.
		Furthermore, QPALM benefits from warm-starting unlike interior-point methods.

\ifarxiv
	\clearpage
\fi
	\begin{appendix}
		\proofsection{thm:PP}

		\begin{appendixproof}{thm:PP}
			The proximal inequality
			\[
				\varphi(x^{k+1})
				{}+{}
				\tfrac12\|x^{k+1}-x^k-e^k\|_{\Sx^{-1}}^2
			{}\leq{}
				\varphi(x^k)
				{}+{}
				\tfrac12\|e^k\|_{\Sx^{-1}}^2,
			\]
			cf. \eqref{eq:proxIneq}, yields
			\begin{equation}\label{eq:inexactSD}
				\varphi(x^{k+1})
				{}+{}
				\tfrac14\|x^{k+1}-x^k\|_{\Sx^{-1}}^2
			{}\leq{}
				\varphi(x^{k+1})
				{}+{}
				\tfrac12\|x^{k+1}-x^k-e^k\|_{\Sx^{-1}}^2
				{}+{}
				\tfrac12\|e^k\|_{\Sx^{-1}}^2
			{}\leq{}
				\varphi(x^k)
				{}+{}
				\|e^k\|_{\Sx^{-1}}^2,
			\end{equation}
			proving assertions \ref{thm:PP:summable} and \ref{thm:PP:bounded}, and similarly \ref{thm:PP:cost} follows by invoking \cite[Lem. 2.2.2]{polyak1987introduction}.
			Next, let \(\seq{x^k}[k\in K]\) be a subsequence converging to a point \(x^\star\); then, it also holds that \(\seq{x^{k+1}}[k\in K]\) converges to \(x^\star\) owing to assertion \ref{thm:PP:summable}.
			From the proximal inequality \eqref{eq:proxIneq} we have
			\[
				\varphi(x^{k+1})
				{}+{}
				\tfrac12
				\|x^{k+1}-x^k-e^k\|_{\Sx^{-1}}^2
			{}\leq{}
				\varphi(x^\star)
				{}+{}
				\tfrac12
				\|x^\star-x^k-e^k\|_{\Sx^{-1}}^2,
			\]
			so that passing to the limit for \(K\ni k\to\infty\) we obtain that \(\limsup_{k\in K}\varphi(x^{k+1})\leq\varphi(x^\star)\).
			In fact, equality holds since \(\varphi\) is lsc, hence from assertion \ref{thm:PP:cost} we conclude that \(\varphi(x^{k+1})\to\varphi(x^\star)\) as \(k\to\infty\), and in turn from the arbitrarity of \(x^\star\) it follows that \(\varphi\) is constantly equal to this limit on the whole set of cluster points.
			To conclude the proof of assertion \ref{thm:PP:omega}, observe that the inclusion
			\(
				\Sx^{-1}(x^k+e^k-x^{k+1})
			{}\in{}
				\hat\partial\varphi(x^{k+1})
			\),
			cf. \eqref{eq:subgradMoreau}, implies that
			\begin{equation}\label{eq:dist}
				\dist\bigl(0,\partial\varphi(x^{k+1})\bigr)
			{}\leq{}
				\dist\bigl(0,\hat\partial\varphi(x^{k+1})\bigr)
			{}\leq{}
				\|\Sx^{-1}\|
				\left(
					\|x^k-x^{k+1}\|
					{}+{}
					\|e^k\|
				\right),
			\end{equation}
			and with limiting arguments (recall that \(\lim_{k\in K}\varphi(x^k)=\varphi(\lim_{k\in K}x^k)\)) the claimed stationarity of the cluster points is obtained.
		\grayout{%
			Suppose now that the hypotheses of assertion \ref{thm:PP:global} are satisfied.
			Then, the set of accumulation points \(\omega\) is nonempty and compact, with \(\dist(x^k,\omega)\to0\) as \(k\to\infty\) and \(\varphi\) is constant on \(\omega\) with value, say, \(\varphi_\star\).
			Let \(\func{\xi}{\R^n\times\R}{\R}\) be defined as
			\(
				\xi(x,t)\coloneqq\varphi(x)+\tfrac12t^2
			\),
			and observe that \(\partial\xi=\partial\varphi\times\id\).
			For \(k\in\N\), define
			\(
				t_k
			{}\coloneqq{}
				\sqrt{
					2
					\sum_{j\geq k}\|e^j\|_{\Sx^{-1}}^2
				}
			\)
			which, as its square is the tail of a convergent series, vanishes as \(k\to\infty\), and in particular \(\seq{x^k,t_k}\) accumulates at \(\omega\times\set0\), where \(\xi\) equals \(\varphi_\star\).
			Clearly, \(\xi\) is also a semialgebraic function and by virtue of \cite{attouch2010proximal} and \cite[Lem. 6]{bolte2014proximal} there thus exist \(\eta,\varepsilon>0\) together with a \emph{desingularizing function} \(\func{\psi}{[0,\eta]}{[0,\infty)}\), namely such that
			\begin{enumeratprop}
				\item\label{thm:UKL1}%
					\(\psi\) is concave, continuous, and increasing on \([0,\eta)\) with \(\psi(0) = 0\);
				\item\label{thm:UKL2}%
					\(\psi\) is \(C^1\) (with \(\psi'>0\)) on \((0,\eta)\);
				\item\label{thm:UKL3}%
					for all points \((x,t)\) such that \(\dist((x,t),\omega\times\set0)<\varepsilon\) and \(\varphi_\star < \xi(x,t) < \varphi_\star + \eta\) it holds that
					\begin{equation}\label{eq:UKL3}
						\psi'\bigl(\xi(x,t)-\varphi_\star\bigr)
						\dist\bigl(0,\partial\xi(x,t)\bigr)
					{}\geq{}
						1.
					\end{equation}
			\end{enumeratprop}
			Observe that
			\begin{align*}
				\xi(x^k,t_k)-\xi(x^{k+1},t_{k+1})
			{}={} &
				\varphi(x^k)-\varphi(x^{k+1})
				{}+{}
				\tfrac{t_k^2-t_{k+1}^2}{2}
			{}={}
				\varphi(x^k)-\varphi(x^{k+1})
				{}+{}
				\|e^k\|_{\Sx^{-1}}^2
			\\
			{}\overrel*[\geq]{\eqref{eq:inexactSD}}{} &
				\tfrac14\|x^k-x^{k+1}\|_{\Sx^{-1}}^2.
			\end{align*}
			In particular, \(\xi(x^k,t_k)\) converges strictly decreasing (to \(\varphi_\star\)).
			By possibly discarding the first iterates, we may assume that \(\xi(x^k,t_k)<\varphi_\star+\eta\) and \(\dist((x^k,t_k),\omega\times\set0)<\varepsilon\) hold for every \(k\in\N\).
			Denoting \(\Delta_k\coloneqq\xi(x^k,t_k)-\varphi_\star\), the convexity of \(\psi\) together with the inequality above and \eqref{eq:UKL3} implies
			\begin{align*}
				\Delta_k-\Delta_{k+1}
			{}\geq{} &
				\overbracket{
					\psi'\bigl(\xi(x^k,t_k)-\varphi_\star\bigr)
				}^{>0}
				\overbracket{
					(\xi(x^k,t_k)-\xi(x^{k+1},t_{k+1}))
				}^{>0}
			{}\geq{}
				\frac{
					\|x^{k+1}-x^k\|_{\Sx^{-1}}^2
				}{
					4\dist(0,\partial\xi(x^k,t_k))
				}
			\\
			{}\geq{} &
				c\frac{
					\|x^{k+1}-x^k\|^2
				}{
					\dist(0,\partial\varphi(x^k))+t_k
				}
			\end{align*}
			for some \(c>0\) (in light of the equivalence of \((x,t)\mapsto\|x\|+|t|\) and the Euclidean metric \((x,t)\mapsto\sqrt{\|x\|^2+t^2}\), and that of \(\|{}\cdot{}\|\) and \(\|{}\cdot{}\|_{\Sx^{-1}}\)).
		%
		
			From \eqref{eq:dist} we thus obtain
			\begin{align*}
				\|x^{k+1}-x^k\|
			{}\leq{} &
				\sqrt{
					\|\Sx^{-1}\|\bigl(
						\|x^k-x^{k-1}\|+\|e^{k-1}\|
					\bigr)
					{}+{}
					t_k
					\,
				}
				\sqrt{
					\vphantom{\Sx^{-1}}
					c^{-1}(\Delta_k-\Delta_{k+1})
					\,
				}
			\\
			{}\leq{} &
				\tfrac{1}{2\alpha}\Bigl(
					\|\Sx^{-1}\|\bigl(
						\|x^k-x^{k-1}\|
						{}+{}
						\underbracket{
							\|e^{k-1}\|
						}
					\bigr)
					{}+{}
					\underbracket{
						\,t_k\,
					}\,
				\Bigr)
				{}+{}
				\tfrac\alpha2
				c^{-1}(
					\underbracket{
						\Delta_k-\Delta_{k+1}
					}
				)
			\end{align*}
			for any \(\alpha>0\), owing to the Cauchy-Schwarz inequality \(\sqrt{uv}\leq\tfrac12(\nicefrac u\alpha+\alpha v)\) holding for any \(u,v>0\).
			Note that the under-bracketed terms are all summable over \(k\in\N\) (the former two by assumption and the latter by a telescoping argument and the fact that \(\Delta_k\geq0\)), hence by choosing \(\alpha=\|\Sx^{-1}\|\) we have that
			\[\textstyle
				\|x^{k+1}-x^k\|
			{}\leq{}
				\tfrac12\|x^k-x^{k-1}\|
				{}+{}
				\delta_k
			{}\leq{}
				(\nicefrac12)^k\|x^1-x^0\|
				{}+{}
				\sum_{j=0}^k(\nicefrac12)^{k-j}\delta_j
			\]
			for some summable sequence \(\seq{\delta_k}\).
			Therefore,
			\[
				\sum_{k\in\N}\|x^{k+1}-x^k\|
			{}\leq{}
				\sum_{k\in\N}{
					(\nicefrac12)^k\|x^1-x^0\|
				}
				{}+{}
				\sum_{k\in\N}{
					\sum_{j=0}^k(\nicefrac12)^{k-j}\delta_j
				}
			{}={}
				2\|x^1-x^0\|
				{}+{}
				2\sum_{k\in\N}\delta_k
			{}<{}
				\infty,
			\]
			hence \(\seq{x^k}\) has finite length and thus converges to a (unique) point, which is stationary owing to assertion \ref{thm:PP:omega}.
		}%
		\end{appendixproof}
	\end{appendix}


	\ifarxiv
		\bibliographystyle{plain}
	\else
		\phantomsection
		\addcontentsline{toc}{section}{References}
		\bibliographystyle{spmpsci}
	\fi
	\bibliography{Bibliography.bib}

\begin{thebibliography}{10}

\bibitem{absil2007newton}
Pierre-Antoine Absil and Andr\'e~L. Tits.
\newblock {N}ewton-{KKT} interior-point methods for indefinite quadratic
  programming.
\newblock {\em Computational Optimization and Applications}, 36(1):5--41, 2007.

\bibitem{amestoy2004algorithm}
Patrick~R. Amestoy, Timothy~A. Davis, and Iain~S. Duff.
\newblock Algorithm 837: {AMD}, an approximate minimum degree ordering
  algorithm.
\newblock {\em ACM Transactions on Mathematical Software (TOMS)},
  30(3):381--388, 2004.

\bibitem{banjac2019infeasibility}
Goran Banjac, Paul Goulart, Bartolomeo Stellato, and Stephen Boyd.
\newblock Infeasibility detection in the alternating direction method of
  multipliers for convex optimization.
\newblock {\em Journal of Optimization Theory and Applications},
  183(2):490--519, 2019.

\bibitem{benzi2005numerical}
Michele Benzi, Gene~H. Golub, and J\"org Liesen.
\newblock Numerical solution of saddle point problems.
\newblock {\em Acta numerica}, 14:1, 2005.

\bibitem{bertsekas1979convexification}
Dimitri~P. Bertsekas.
\newblock Convexification procedures and decomposition methods for nonconvex
  optimization problems.
\newblock {\em Journal of Optimization Theory and Applications},
  29(2):169--197, Oct 1979.

\bibitem{bertsekas1982constrained}
Dimitri~P. Bertsekas.
\newblock Constrained optimization and {L}agrange multiplier methods.
\newblock {\em Computer Science and Applied Mathematics, Boston: Academic
  Press, 1982}, 1982.

\bibitem{bertsekas2016nonlinear}
Dimitri~P. Bertsekas.
\newblock {\em Nonlinear Programming}.
\newblock Athena Scientific, 2016.

\bibitem{birgin2014practical}
Ernesto~G. Birgin and José~Mario Martínez.
\newblock {\em Practical Augmented {L}agrangian Methods for Constrained
  Optimization}.
\newblock Society for Industrial and Applied Mathematics, Philadelphia, PA,
  2014.

\bibitem{bolte2018nonconvex}
Jérôme Bolte, Shoham Sabach, and Marc Teboulle.
\newblock Nonconvex {L}agrangian-based optimization: Monitoring schemes and
  global convergence.
\newblock {\em Mathematics of Operations Research}, 43(4):1210--1232, 2018.

\bibitem{bot2020proximal}
Radu~Ioan Boţ and Dang-Khoa Nguyen.
\newblock The proximal alternating direction method of multipliers in the
  nonconvex setting: Convergence analysis and rates.
\newblock {\em Mathematics of Operations Research}, 45(2):682--712, 2020.

\bibitem{burer2008finite}
Samuel Burer and Dieter Vandenbussche.
\newblock A finite branch-and-bound algorithm for nonconvex quadratic
  programming via semidefinite relaxations.
\newblock {\em Mathematical Programming}, 113(2):259--282, 2008.

\bibitem{chen2012globally}
Jieqiu Chen and Samuel Burer.
\newblock Globally solving nonconvex quadratic programming problems via
  completely positive programming.
\newblock {\em Mathematical Programming Computation}, 4(1):33--52, 2012.

\bibitem{chen2008algorithm}
Yanqing Chen, Timothy~A. Davis, William~W. Hager, and Sivasankaran
  Rajamanickam.
\newblock Algorithm 887: {CHOLMOD}, supernodal sparse {C}holesky factorization
  and update/downdate.
\newblock {\em ACM Transactions on Mathematical Software (TOMS)}, 35(3):1--14,
  2008.

\bibitem{combettes2004proximal}
Patrick~L. Combettes and Teemu Pennanen.
\newblock Proximal methods for cohypomonotone operators.
\newblock {\em SIAM journal on control and optimization}, 43(2):731--742, 2004.

\bibitem{cottle1970classes}
Richard~W. Cottle, G.J. Habetler, and C.E. Lemke.
\newblock On classes of copositive matrices.
\newblock {\em Linear Algebra and Its Applications}, 3(3):295--310, 1970.

\bibitem{davis2005algorithm}
Timothy~A. Davis.
\newblock Algorithm 849: A concise sparse {C}holesky factorization package.
\newblock {\em ACM Transactions on Mathematical Software (TOMS)},
  31(4):587--591, 2005.

\bibitem{davis2006direct}
Timothy~A. Davis.
\newblock {\em Direct Methods for Sparse Linear Systems}.
\newblock Society for Industrial and Applied Mathematics, 2006.

\bibitem{davis1999modifying}
Timothy~A. Davis and William~W. Hager.
\newblock Modifying a sparse {C}holesky factorization.
\newblock {\em SIAM Journal on Matrix Analysis and Applications},
  20(3):606--627, 1999.

\bibitem{davis2001multiple}
Timothy~A. Davis and William~W. Hager.
\newblock Multiple-rank modifications of a sparse {C}holesky factorization.
\newblock {\em SIAM Journal on Matrix Analysis and Applications},
  22(4):997--1013, 2001.

\bibitem{davis2005row}
Timothy~A. Davis and William~W. Hager.
\newblock Row modifications of a sparse {C}holesky factorization.
\newblock {\em SIAM Journal on Matrix Analysis and Applications},
  26(3):621--639, 2005.

\bibitem{dolan2002benchmarking}
Elizabeth~D. Dolan and Jorge~J. Mor\'e.
\newblock Benchmarking optimization software with performance profiles.
\newblock {\em Mathematical programming}, 91(2):201--213, 2002.

\bibitem{dontchev2009implicit}
Asen~L. Dontchev and R.~Tyrrell Rockafellar.
\newblock {\em Implicit functions and solution mappings}, volume 208.
\newblock Springer, 2009.

\bibitem{facchinei2003finite}
Francisco Facchinei and Jong-Shi Pang.
\newblock {\em Finite-dimensional variational inequalities and complementarity
  problems}, volume~II.
\newblock Springer, 2003.

\bibitem{ferreau2014qpoases}
Hans~Joachim Ferreau, Christian Kirches, Andreas Potschka, Hans~Georg Bock, and
  Moritz Diehl.
\newblock {qpOASES}: A parametric active-set algorithm for quadratic
  programming.
\newblock {\em Mathematical Programming Computation}, 6(4):327--363, 2014.

\bibitem{frison2020hpipm}
Gianluca Frison and Moritz Diehl.
\newblock Hpipm: a high-performance quadratic programming framework for model
  predictive control.
\newblock {\em IFAC-PapersOnLine}, 53(2):6563--6569, 2020.

\bibitem{gertz2003object}
E.~Michael Gertz and Stephen~J. Wright.
\newblock Object-oriented software for quadratic programming.
\newblock {\em {ACM} Transactions on Mathematical Software ({TOMS})},
  29(1):58--81, 2003.

\bibitem{gill2015methods}
Philip~E. Gill and Elizabeth Wong.
\newblock Methods for convex and general quadratic programming.
\newblock {\em Mathematical programming computation}, 7(1):71--112, 2015.

\bibitem{golub2013matrix}
Gene~H. Golub and Charles~F. Van~Loan.
\newblock {\em Matrix Computations}.
\newblock Johns Hopkins Studies in the Mathematical Sciences. Johns Hopkins
  University Press, 2013.

\bibitem{gould2016note}
Nicholas Gould and Jennifer Scott.
\newblock A note on performance profiles for benchmarking software.
\newblock {\em ACM Transactions on Mathematical Software (TOMS)}, 43(2):1--5,
  2016.

\bibitem{gould2003galahad}
Nicholas~I.M. Gould, Dominique Orban, and Philippe~L. Toint.
\newblock {GALAHAD}, a library of thread-safe fortran 90 packages for
  large-scale nonlinear optimization.
\newblock {\em {ACM} Transactions on Mathematical Software ({TOMS})},
  29(4):353--372, 2003.

\bibitem{gould2015cutest}
Nicholas~I.M. Gould, Dominique Orban, and Philippe~L. Toint.
\newblock {CUTEst}: a constrained and unconstrained testing environment with
  safe threads for mathematical optimization.
\newblock {\em Computational Optimization and Applications}, 60(3):545--557,
  2015.

\bibitem{gurobi2018gurobi}
LLC Gurobi~Optimization.
\newblock Gurobi optimizer reference manual, 2018.

\bibitem{hermans2019qpalm}
Ben Hermans, Andreas Themelis, and Panagiotis Patrinos.
\newblock {{QP}ALM}: A {N}ewton-type proximal augmented {L}agrangian method for
  quadratic programs.
\newblock In {\em 2019 IEEE 58th Conference on Decision and Control (CDC)},
  pages 4325--4330, 2019.

\bibitem{iusem2003inexact}
Alfredo~N. Iusem, Teemu Pennanen, and Benar~F. Svaiter.
\newblock Inexact variants of the proximal point algorithm without
  monotonicity.
\newblock {\em SIAM Journal on Optimization}, 13(4):1080--1097, 2003.

\bibitem{knyazev2001toward}
Andrew~V. Knyazev.
\newblock Toward the optimal preconditioned eigensolver: Locally optimal block
  preconditioned conjugate gradient method.
\newblock {\em SIAM journal on scientific computing}, 23(2):517--541, 2001.

\bibitem{kong2019complexity}
Weiwei Kong, Jefferson~G. Melo, and Renato~D.C. Monteiro.
\newblock Complexity of a quadratic penalty accelerated inexact proximal point
  method for solving linearly constrained nonconvex composite programs.
\newblock {\em SIAM Journal on Optimization}, 29(4):2566--2593, 2019.

\bibitem{li2015global}
Guoyin Li and Ting~Kei Pong.
\newblock Global convergence of splitting methods for nonconvex composite
  optimization.
\newblock {\em SIAM Journal on Optimization}, 25(4):2434--2460, 2015.

\bibitem{lin2019inexact}
Qihang Lin, Runchao Ma, and Yangyang Xu.
\newblock Inexact proximal-point penalty methods for non-convex optimization
  with non-convex constraints.
\newblock {\em arXiv preprint arXiv:1908.11518}, 2019.

\bibitem{luo1993error}
Zhi-Quan Luo and Paul Tseng.
\newblock Error bounds and convergence analysis of feasible descent methods: a
  general approach.
\newblock {\em Annals of Operations Research}, 46(1):157--178, 1993.

\bibitem{maros1999repository}
Istvan Maros and Csaba M\'esz\'aros.
\newblock A repository of convex quadratic programming problems.
\newblock {\em Optimization Methods and Software}, 11(1-4):671--681, 1999.

\bibitem{meszaros1999bpmpd}
Csaba M\'esz\'aros.
\newblock The {BPMPD} interior point solver for convex quadratic problems.
\newblock {\em Optimization Methods and Software}, 11(1-4):431--449, 1999.

\bibitem{mosek2019mosek}
ApS {MOSEK}.
\newblock Mosek optimization toolbox for matlab.
\newblock {\em User's Guide and Reference Manual, Version}, 9.2.22, 2019.

\bibitem{nocedal2006numerical}
Jorge Nocedal and Stephen Wright.
\newblock {\em Numerical optimization}.
\newblock Springer Science \& Business Media, 2006.

\bibitem{patrinos2010new}
Panagiotis Patrinos and Haralambos Sarimveis.
\newblock A new algorithm for solving convex parametric quadratic programs
  based on graphical derivatives of solution mappings.
\newblock {\em Automatica}, 46(9):1405 -- 1418, 2010.

\bibitem{polyak1987introduction}
Boris~T. Polyak.
\newblock Introduction to optimization.
\newblock {\em Inc., Publications Division, New York}, 1, 1987.

\bibitem{rockafellar1976augmented}
R.~Tyrrell Rockafellar.
\newblock Augmented {L}agrangians and applications of the proximal point
  algorithm in convex programming.
\newblock {\em Mathematics of Operations Research}, 1(2):97--116, 1976.

\bibitem{rockafellar2011variational}
R.~Tyrrell Rockafellar and Roger~J.B. Wets.
\newblock {\em Variational analysis}, volume 317.
\newblock Springer Science \& Business Media, 2011.

\bibitem{ruiz2001scaling}
Daniel Ruiz.
\newblock A scaling algorithm to equilibrate both rows and columns norms in
  matrices.
\newblock Technical report, Rutherford Appleton Laboratorie, 2001.

\bibitem{sherali1995reformulation}
Hanif~D. Sherali and Cihan~H. Tuncbilek.
\newblock A reformulation-convexification approach for solving nonconvex
  quadratic programming problems.
\newblock {\em Journal of Global Optimization}, 7(1):1--31, 1995.

\bibitem{stellato2020osqp}
Bartolomeo Stellato, Goran Banjac, Paul Goulart, Alberto Bemporad, and Stephen
  Boyd.
\newblock {OSQP}: An operator splitting solver for quadratic programs.
\newblock {\em Mathematical Programming Computation}, 2020.

\bibitem{sun2017convergence}
Tao Sun, Hao Jiang, Lizhi Cheng, and Wei Zhu.
\newblock A convergence framework for inexact nonconvex and nonsmooth
  algorithms and its applications to several iterations.
\newblock {\em arXiv preprint arXiv:1709.04072}, 2017.

\bibitem{themelis2018acceleration}
Andreas Themelis, Masoud Ahookhosh, and Panagiotis Patrinos.
\newblock On the acceleration of forward-backward splitting via an inexact
  {N}ewton method.
\newblock In R.~Luke, H.~Bauschke, and R.~Burachik, editors, {\em Splitting
  Algorithms, Modern Operator Theory, and Applications}. Springer, 2019.

\bibitem{themelis2020douglas}
Andreas Themelis and Panagiotis Patrinos.
\newblock {D}ouglas--{R}achford splitting and {{ADMM}} for nonconvex
  optimization: Tight convergence results.
\newblock {\em SIAM Journal on Optimization}, 30(1):149--181, 2020.

\bibitem{vanderbei1995symmetric}
Robert~J. Vanderbei.
\newblock Symmetric quasidefinite matrices.
\newblock {\em SIAM Journal on Optimization}, 5(1):100--113, 1995.

\bibitem{wachter2006implementation}
Andreas W{\"a}chter and Lorenz~T. Biegler.
\newblock On the implementation of an interior-point filter line-search
  algorithm for large-scale nonlinear programming.
\newblock {\em Mathematical programming}, 106(1):25--57, 2006.

\bibitem{yannakakis1981computing}
Mihalis Yannakakis.
\newblock Computing the minimum fill-in is {NP}-complete.
\newblock {\em SIAM Journal on Algebraic Discrete Methods}, 2(1):77--79, 1981.

\bibitem{ye1992affine}
Yinyu Ye.
\newblock On affine scaling algorithms for nonconvex quadratic programming.
\newblock {\em Mathematical Programming}, 56(1-3):285--300, 1992.

\end{thebibliography}

\end{document}